\documentclass[preprint]{imsart}

\usepackage{amsmath,amssymb}
\usepackage{graphicx}
 \usepackage{yhmath}
\RequirePackage[numbers]{natbib}
\RequirePackage[colorlinks,citecolor=blue,urlcolor=blue]{hyperref}

\startlocaldefs

\input{diagrams.sty}

\graphicspath{{./graphics/}}

\renewcommand{\Bar}{\overline}
\newcommand{\bfr}{{\bf R}}
\newcommand{\dotprod}{{\mbox{
\hspace{-.06in}\raisebox{-.015in}{\bf {\Large $\cdot$}}}}}
\renewcommand{\Hat}{\widehat}
\newcommand{\id}{{\rm id.}}
\newcommand{\ident}{\equiv}
\newcommand{\intersect}{\mbox{\small \ $\bigcap$\ }}
\newcommand{\iso}{\cong}

\renewcommand{\Re}{{\rm Re}}
\newcommand{\plus}{\oplus}

\newcommand{\sgn}{{\rm sgn}}
\newcommand{\spn}{{\rm span}}
\newcommand{\Sym}{{\rm Sym}}
\renewcommand{\Tilde}{\widetilde}

\newcommand{\tr}{{\rm tr}}
\newcommand{\union}{\mbox{\small $\bigcup$}}
\newcommand{\Union}{\bigcup}
\newcommand{\be}{\begin{equation}}
\newcommand{\ee}{\end{equation}}
\newcommand{\bearray}{\begin{eqnarray}}
\newcommand{\eearray}{\end{eqnarray}}
\newcommand{\bestar}{\begin{eqnarray*}}
\newcommand{\eestar}{\end{eqnarray*}}
\newcommand{\ben}{\begin{displaymath}}
\newcommand{\een}{\end{displaymath}}

\newtheorem{theorem}{Theorem}[section]
\newtheorem{thm}[theorem]{Theorem}

\newtheorem{prop}[theorem]{Proposition}

\newtheorem{cor}[theorem]{Corollary}
\newtheorem{defn}[theorem]{Definition}

\newtheorem{remark}[theorem]{Remark}

\newtheorem{notation}[theorem]{Notation}
\newtheorem{result}[theorem]{Result}

\newtheorem{Convention}[theorem]{Convention}

\newenvironment{example}%
{\addtocounter{theorem}{1}%
\noindent {\bf
Example \arabic{section}.\arabic{theorem} }}{\\[2 ex]}

\def\qedns{\mbox{\hspace{1ex}}
\mbox{\hfill\vrule height10pt width10pt} \\ \vspace{0.5 ex}
}

\newcommand{\pf}{\noi{\bf Proof}: }

\renewcommand{\b}{\beta}
\renewcommand{\d}{\delta}
\newcommand{\e}{\epsilon}
\newcommand{\g}{\gamma}

\renewcommand{\l}{\lambda}
\renewcommand{\th}{\theta}
\newcommand{\z}{\zeta}

\newcommand{\noi}{\noindent}

\newcommand{\bsn}{\vspace{.2in}\noindent}

\newcommand{\indenum}{\mbox{\hspace{3ex}}}

\newcommand{\bfc}{{\bf C}}
\newcommand{\Sop}{SO(p)}
\newcommand{\So}{SO}
\newcommand{\sop}{{\mathfrak {so}}(p)}

\newcommand{\sympp}{{\rm Sym}^+(p)}

\newcommand{\ts}{{\tilde{S}}}
\newcommand{\tsp}{\ts_p}

\newcommand{\D}{{\rm Diag}}
\newcommand{\Dp}{{\rm Diag}(p)}
\newcommand{\Dpp}{{\rm Diag}^+(p)}
\newcommand{\E}{{\cal E}}

\newcommand{\I}{{\cal I}}
\newcommand{\J}{{\sf J}}
\newcommand{\K}{{\sf K}}
\renewcommand{\L}{{\Lambda}}
\newcommand{\M}{{\cal M}}

\newcommand{\Mthree}{({\rm SO} \times {\rm Diag}^+)(3)}

\renewcommand{\P}{{\cal P}}

\renewcommand{\S}{{\cal S}}

\newcommand{\dsr}{d_{\cal SR}}
\newcommand{\partpset}{{\rm Part}(\{1,\dots,p\})}
\newcommand{\partp}{{\rm Part}(p)}

\newcommand{\mbbi}{I}

\newcommand{\boldsig}{{\mbox{\boldmath $\sigma$ \unboldmath}
           \mbox{\hspace{-.05in}}}}
\newcommand{\miniboldsig}{{\mbox{\scriptsize \boldmath $\sigma$ \unboldmath}
           \mbox{\hspace{-.05in}}}}
\newcommand{\mbs}{\miniboldsig}
\newcommand{\ztwo}{{\bf Z}_2}
\renewcommand{\H}{{\bf H}}
\newcommand{\imh}{{\rm Im(\H)}}

\newcommand{\piid}{{\pi_{\rm id}}}
\renewcommand{\Im}{{\rm Im}}
\newcommand{\ET}{{\rm ET}}
\newcommand{\bfd}{{\bf d}}

\newcommand{\diag}{{\rm diag}}
\newcommand{\bone}{{\bf 1}}

\newcommand{\proj}{{\rm proj}}
\newcommand{\sproj}{\mbox{\small ${\rm proj}_2$}}
\newcommand{\bsproj}{\mbox{\small $\Bar{{\rm proj}_2}$}}
\newcommand{\mat}{{\sf mat}}
\newcommand{\Comp}{{\rm Comp}}
\newcommand{\cald}{{\cal D}}

\newcommand{\gdp}{{g_{\cald^+}}}

\newcommand{\dso}{d_{\So}}

\newcommand{\ddp}{d_{\cald^+}}
\newcommand{\esstop}{\S_{\rm top}}
\newcommand{\smid}{\S_{\rm mid}}
\newcommand{\sbot}{\S_{\rm bot}}
\newcommand{\jtop}{{\J_{\rm top}}}
\newcommand{\jbot}{{\J_{\rm bot}}}
\newcommand{\simc}{\sim_c}
\newcommand{\lbl}{{\rm lbl}}
\newcommand{\slbl}{\mbox{\small ${\rm lbl}$}}
\newcommand{\bslbl}{\mbox{\small $\Bar{{\rm lbl}}$}}
\newcommand{\quo}{{\rm quo}}

\newcommand{\Real}{\mathbf R}

\def\argmin{\mathop{\rm argmin}}

\def\tr{\mbox{trace}}
\def\half{\frac{1}{2}}
\def\proj{\mbox{Proj}}
\def\diag{\mbox{diag}}

\def\symp{\mbox{Sym}^+}
\def\sym{\mbox{Sym}}
\def\SO{\mbox{SO}}

\def\diagp{\mbox{Diag}^+}
\def\diag{\mbox{Diag}}

\def\Complex{\mathbf{C}}

\newcommand{\Dc}{\mathcal{D}}
\newcommand{\Ec}{\mathcal{E}}

\newcommand{\Mc}{\mathcal{M}}

\newcommand{\Rc}{\mathcal{R}}
\newcommand{\Sc}{\mathcal{S}}

\newcommand{\Nsf}{{\sf N}}

\def\diagg{\mbox{diag}}
\def\1v{\mathbf 1}
\def\0v{\mathbf 0}

\newcommand{\Stop}{\Sc_{\rm top}}
\newcommand{\Sbot}{\Sc_{\rm bot}}
\newcommand{\Smid}{\Sc_{\rm mid}}

\endlocaldefs

\begin{document}

\begin{frontmatter}

\title{
Geometric foundations for scaling-rotation statistics on symmetric positive definite matrices: minimal smooth scaling-rotation curves in low dimensions
 \thanksref{t1}}
\thankstext{t1}{This work was supported by NIH grant
  R21EB012177 and NSF grant DMS-1307178}
\runtitle{Eigenvalue stratification and MSSR curves}

\author{\fnms{David} \snm{Groisser}\corref{}\ead[label=e1]{groisser@ufl.edu}}
\address{Department of Mathematics, University of Florida, Gainesville, FL 32611, USA\\ \printead{e1}}
\author{\fnms{Sungkyu} \snm{Jung}\ead[label=e2]{sungkyu@pitt.edu}}
\address{Department of Statistics, University of Pittsburgh, Pittsburgh, PA 15260, USA\\ \printead{e2}}
\and
\author{\fnms{Armin} \snm{Schwartzman}\ead[label=e3]{armins@ucsd.edu}}
\address{Division of Biostatistics, University of California, San Diego, CA 92093, USA\\ \printead{e3}}

\runauthor{Groisser et al.}

\begin{abstract}
We investigate a  geometric computational framework, called the ``scaling-rotation framework'', on ${\rm Sym}^+(p)$, the set of $p \times p$ symmetric positive-definite (SPD) matrices.
 The purpose of our study is to lay geometric foundations for statistical analysis of SPD matrices,  in situations in which eigenstructure is of fundamental importance, for example diffusion-tensor imaging (DTI). Eigen-decomposition, upon which the scaling-rotation framework is based, determines both a stratification of ${\rm Sym}^+(p)$, defined by eigenvalue multiplicities, and fibers of the ``eigen-composition" map  $SO(p)\times{\rm Diag}^+(p)\to{\rm Sym}^+(p)$. This leads to the notion of {\em scaling-rotation distance} [Jung et al. (2015)], a measure of the minimal amount of scaling and rotation needed to transform an SPD matrix, $X,$ into another, $Y,$ by a smooth curve in ${\rm Sym}^+(p)$.
Our main goal in this paper is the systematic characterization and analysis of {\em minimal smooth scaling-rotation (MSSR) curves},
images in ${\rm Sym}^+(p)$ of minimal-length geodesics connecting two fibers  in the ``upstairs" space $SO(p)\times{\rm Diag}^+(p)$.
The length of such a geodesic connecting the fibers over $X$ and $Y$ is what we define to be the scaling-rotation distance from $X$ to $Y.$
For the important low-dimensional case $p = 3$ (the home of DTI), we find new explicit formulas for MSSR curves and for the scaling-rotation distance, and identify  ${\cal M}(X,Y)$  in all ``nontrivial" cases.  The quaternionic representation of $SO(3)$ is used in these computations. We also provide closed-form expressions for scaling-rotation  distance and MSSR curves for the case $p = 2$.
\end{abstract}





\end{frontmatter}

\section{Introduction}

In recent years there has been increased interest in stratified manifolds for
statistical applications.
For example, stratified manifolds have recently received attention in the study of
phylogenetic trees \cite{Billera2001,Hotz2013} and Kendall's 3D shape
space \cite{Kendall1999}. New analytic tools for such manifolds are
fast developing \cite{Damon2013, Bhattacharya2013}.
Our
work contributes to the development of such tools on both a theoretical and practical
level, providing a solid geometrical
foundation for development of statistical procedures on the stratified
manifold $\sympp$, the set of $p \times p$ symmetric positive-definite (SPD) matrices.

In this work, we investigate a geometric structure on  $\sympp$, resulting from the stratification defined by eigenvalue multiplicities. This stratification is tied inextricably to our main goal in this paper: the systematic characterization and analysis of minimal smooth scaling-rotation curves in low dimensions. Such curves were defined in \cite{JSG2015scarot} as smooth curves whose length minimizes the amount of scaling and rotation needed to transform an SPD matrix into another. The techniques developed in this paper, when applied to the case $p=3$, allow us to find new explicit formulas for such curves.
Our work builds fundamental mathematical and geometric grounds that facilitate developments of statistical procedures for SPD matrices, and is instrumental in understanding general stratified manifolds.

To elaborate how our work here relates to advancing statistical analysis of SPD matrices, we present below a rather long introduction. We first give some background on statistical analysis of SPD matrices, and more generally on analysis of data in stratified manifolds, followed by a brief discussion on the statistical motivation of studying scaling-rotation curves and distances. We then informally introduce the main results of the paper.

\subsection{Background }\label{sec:intro-background}

 \paragraph{Statistical analysis of SPD matrices} The statistical analysis of SPD matrices has several applications, especially in some biological problems, such as diffusion-tensor imaging (DTI).  A diffusion tensor may be viewed as an ellipsoid, represented by a $3\times 3$ SPD matrix. DTI researchers are interested in smoothing a raw noisy diffusion-tensor field  \cite{zhou2016regularisation}, registering fibers of tensor fields \cite{Alexander2005}, regression models \cite{Zhu2009,Yuan2012} and classification of `noisy' tensors into strata \cite{Zhu2007}.
Our eigenvalue-multiplicity stratification categorizes the ellipsoids associated with the SPD matrices into distinct shapes, which in the case $p=3$ are known as spherical, prolate/oblate, and tri-axial (scalene).  We believe that the scaling-rotation framework studied in this work and in \cite{JSG2015scarot,GJS2016mathpaper} will be highly useful in developing new methodologies of smoothing, registration and regression analysis of diffusion tensors.

A major hurdle in analyzing SPD matrix-valued data is that
 the data are best viewed as lying in a curved space, making the application of conventional statistical tools inappropriate.
To briefly discuss the drawback of using a naive approach (i.e., using the fact that the data lie in the vector space of all $p\times p$ symmetric matrices), take as an example the simplest case of $ 2 \times 2$ SPD matrices. In order for a $ 2 \times 2$  symmetric matrix $X$ to be positive-definite, the squared off-diagonal element $x_{12}$ must be absolutely smaller than    the product of two diagonal elements $x_{11}$ and $x_{22}$ (which themselves must be positive). This entails the set $\{(x_{11},x_{22},x_{12}) :  X = (x_{ij}) \in \rm{Sym}^+(2)\}$ being a proper subset of $\mathbf{R}^3$, the set of points inside of a convex cone (this is visualized in Fig. \ref{fig:Symp2stratification} in Section~\ref{sec:exsr2}.)
A naive approach to handle data in $\sympp$ is to use the usual metric defined in the ambient space, which gives rise to Euclidean metric $d_E(X,Y) = \|{X-Y}\|_F$ (Frobenius norm). There are several disadvantages of using Euclidean metric: the straight line given by the Euclidean framework  has undesirable features such as ``swelling'' \cite{Arsigny2007} and limited extrapolation.
Recently, several different geometric tools have been proposed to handle the data as lying in a curved space with the help of Riemannian geometry and Lie group theory \cite{Moakher2005,Arsigny2007,Pennec2006,Lenglet2006,Osborne2013,Schwartzman2006,Schwartzman2010} or by borrowing ideas from shape analysis \cite{Dryden2009,zhou2013procrustes,zhou2016regularisation}.
Among these, we point out three existing frameworks.

The log-Euclidean geometric framework \cite{Moakher2005,Arsigny2007} handles the data in a ``log-transformed space'', the set of symmetric matrices, $\sym(p) = \log(\symp(p))$. This gives rise to the log-Euclidean metric $d_L(X,Y) = \|\log(X) - \log(Y) \|_F$.
Effectively, the log-transform provides a ``local linearization" of $\symp(p)$ near the identity matrix;  the results it yields are less good for matrices farther from the identity.
A second framework, the ``affine-invariant Riemannian framework" \cite{Pennec2006}, provides a local linearization of $\symp(p)$ in a neighborhood of an arbitrary point $\mu \in \symp(p)$.
This framework makes use of the identification of $\sym(p)$ with the tangent space of $\symp(p)$ at $\mu$ to endow $\symp(p)$  with a $GL(p,R)$-invariant Riemannian metric.
 This gives rise to the metric $d_{AI}(X,Y) = \| \log (X^{-\half}YX^{-\half}) \|_F$. When $X$ and $Y$ are understood as covariance matrices of random vectors $x$ and $y$, the distance $d_{AI}(X,Y)$ is invariant under ``affine'' transformations applied to both $X,Y$; for any $p \times p $ invertible matrix $G$, $d_{AI}(X,Y) = d_{AI}(GXG^T,GYG^T)$.
From a third standpoint, the Procrustes size-and-shape framework of \cite{Dryden2009} turns the problem of analyzing SPD matrices into a problem of analyzing reflection size-and-shapes of $(p+1)$-landmark configurations in $p$ dimensions. Specifically, an SPD matrix $X$ is represented by an equivalent class $\{{LR} : R \in O(p)\}$, where the  lower triangular matrix $L$ satisfies $X = LL^T$. The size-and-shape metric is defined as $d_{S}(X_1,X_2) = \inf_{R\in O(p)}\| L_1 - L_2 R\|_F$, where $L_i$ satisfies $X_i = L_iL_i^T$. The size-and-shape framework can also be applied to symmetric non-negative definite matrices.

These three different measures of ``distance'' dictate the method of interpolation of two or more SPD matrices, and lead to different definitions of the population and sample mean. The results of smoothing a tensor field and registration of fiber tracts will also depend on the choice of geometric framework for computation.
These frameworks also provide methods for local linearization of data, methods that are useful for e.g. dimension-reduction, regression modeling, approximate multivariate-normal-based inference and large-sample asymptotic distributions. The log-transformation-based geometric frameworks, log-Euclidean and affine-invariant Riemannian frameworks, have been heavily used in statistical modeling and estimations \cite[cf.][]{Schwartzman2010,Zhu2007}, partly due to their simple geometric structures.
In  previous work \cite{JSG2015scarot}, we introduced a fourth framework, the ``scaling-rotation framework'', that is the subject of this paper.  In \cite[Section 5]{JSG2015scarot}, we presented evidence of advantages of this framework over the popular log-transformation-based frameworks for tensor interpolations.
In Section~\ref{sec:intro-importance} of the present paper, we briefly discuss some other advantages of the scaling-rotation framework in statistical analysis.

\paragraph{Statistical analysis of data on stratified spaces}
As we shall see in this paper, the scaling-rotation framework leads us to treat $\sympp$ as a stratified space.
Many statistical analyses now deal with data that naturally lie in non-Euclidean spaces.
In particular, stratified spaces have recently received attention in the study of, e.g., phylogenetic trees \citep{Billera2001} and Kendall's 3D shape space \citep{Kendall1999}.
A stratified space is a union of ``nice'' topological subspaces called \emph{strata}, with certain restrictions on the way the strata join. A simple example is a \emph{spider} (half-lines joined by a point) or an \emph{open book} (half-planes joined by a line) \citep{Hotz2013}. Another example is the phylogenetic tree space of Billera, Holmes and Vogtmann \citep{Billera2001},  the union of Euclidean positive orthants, each representing different topology of phylogenetic trees (see also \citep{nye2011principal}).
The space of SPD matrices is naturally stratified by eigenvalue multiplicities. For example if $p = 2$, there are two strata, one consisting of SPD matrices with distinct eigenvalues and the other consisting of matrices with equal eigenvalues.


For statistical analysis on  stratified spaces, it is crucial to devise appropriate notions of distance and  shortest path(s) between two points, together with associated computational algorithms.
These  tasks, in general, are challenging. For example, it is known that for computing a graph-edit distance between two geometric tree-like shapes is NP-complete \citep{bille2005survey}. To overcome these computational burdens, Feragen and her colleagues \citep{feragen2012complexity,feragen2013toward} have proposed and studied a quotient Euclidean distance on the space of tree-like shapes, which is a stratified space.
Wang and Marron \cite{wang2007object} defined a notion of ``average tree'' as well as a principal-component analysis of trees, and an efficient algorithm \cite{aydin2009principal} was needed to compute the principal components. For the phylogenetic-tree spaces, there has been an ongoing effort to advance efficient computations for distances \cite{owen2011fast}, mean and median \cite{bacak2014computing,miller2015polyhedral}, clustering \cite{chakerian2012computational}, and estimating principal components \cite{nye2011principal}. For stratified shape-spaces, Huckemann et al. \cite{huckemann2010intrinsic} have also developed a form of principal component analysis.

New analytic tools for these stratified spaces are fast developing. Hotz et al. \cite{Hotz2013} established a central limit theorem for the open-book space, and showed that the sample Fr\'{e}chet mean can be ``sticky'' to the one-dimensional stratum. For a special phylogenetic-tree space,  central limit theorems were derived in \cite{barden2013central} for each of three cases: when the population Fr\'{e}chet mean is in the top stratum, a co-dimension-one stratum, or the bottom stratum (a point). See \cite{barden2014limiting} for an extension.  Nye has defined diffusion processes for some simple stratified spaces \cite{nye2014diffusion}  and   for the phylogenetic-tree space \cite{nye2015convergence}. See \cite{feragen2014mini} and references therein for other recent developments.

In analogy to the literature on tree spaces, in this paper we develop the concepts of shortest paths and scaling-rotation distance, and provide closed-form formulas, as  first steps toward developing eigenstructure-based statistics on $\sympp$.
%
In the future, new concepts and analytical tools such as mean, principal component analysis, regression analysis, and inference procedures may be developed within the scaling-rotation framework. Our work contributes to the development of such tools on both theoretical and practical level, providing a solid geometrical foundation for development of eigenstructure-based statistical procedures on the stratified manifold $\sympp$.

\subsection{Scaling-rotation geometric framework and its statistical importance}
\label{sec:intro-importance}


Recall that every $X\in\sympp$ can be diagonalized by a rotation
matrix: $X=UDU^{-1}=UDU^T$ for some $U\in \Sop, D\in \D^+(p)$.  Here,
$\D^+(p)$ denotes the set of $p\times p$ diagonal matrices all of whose diagonal entries are positive. We refer to $(U,D)$ as an {\em eigen-decomposition}
of $X$.  Conversely, for all $U\in \Sop, D\in \D^+(p)$, the matrix
$UDU^T$ lies in $\sympp$.  Thus the {\em space of eigen-decompositions
  of $p\times p$ SPD matrices} is the manifold
\be
M:=M(p):=(SO\times\D)^+(p) := \Sop\times \D^+(p).
\ee

\noi This manifold comes to us naturally equipped with a smooth
surjective map $F: M\to\sympp$ defined by
\be
F(U,D) = UDU^T.
\ee

\noi {To name the set of eigen-decompositions corresponding to a single SPD matrix,} for each $X\in \sympp$, we define the {\em fiber over $X$} to be
the set
\ben
\E_X := F^{-1}(X) = \{(U,D)\in M : UDU^T=X\}.
\een
\noi The relation $\sim$ on $M$ defined by lying in the same
fiber---i.e. $(U,D)\sim (V,\L)$ if and only if $F(U,D)=F(V,\L)$---is
an equivalence relation. The quotient space $M/\sim$ (the set of
equivalence classes, endowed with the quotient topology) is
canonically identified with $\sympp$.  It should be noted that $F$ is
not a submersion (cf. \cite{Absil2009, CE1975}), and that $M$ is {\em not} a
fiber bundle over $\sympp$; as we will see explicitly later, the
fibers are not all mutually diffeomorphic (or even of the same
dimension).

The different structures of fibers  naturally lead to a
stratification of $\sympp$ and $M$.  The stratum to which an $X \in
\sympp$ belongs depends on the diffeomorphism type of $\E_X$.
As we shall see in Section \ref{sect:three-stratification}, this
stratification based on ``fiber types'' is equivalent to stratifications by
orbit-type and by eigenvalue-multiplicity type.

The strata of $\sympp$ and $M$ are determined by  patterns of
eigenvalue multiplicities, and are labeled by partitions of the
integer $p$ and the set $\{1,\ldots,p\}$. We will always assume $p>1$,
the case $p=1$ being uninteresting. For each $p$, one can obtain the
numbers of strata (of $\sympp$ and $M$), the dimension of each
stratum, and the diffeomorphism type of fibers belonging to each
stratum. Several group-actions are involved, and the deepest
understanding comes from identifying the relevant groups and the
various actions.

In \cite{Schwartzman2006}, Schwartzman introduced scaling-rotation
curves as a way of interpolating between SPD matrices in such a way
that eigenvectors and eigenvalues both change at uniform speed.
To provide a geometric framework for these curves, Section 2 is
devoted to systematic characterization of fibers and its connection to
the stratification of $\sympp$. This allows us to build upon the
scaling-rotation framework for SPD matrices proposed in
\cite{JSG2015scarot}, which provided a geometric interpretation for
the scaling-rotation curves in \cite{Schwartzman2006}. In particular,
our characterization of fibers is essential in understanding
differential topology and geometry of this framework.

In the scaling-rotation framework for SPD matrices, the ``distance"
$\dsr(X,Y)$ between any two matrices  $X,Y\in \sympp$  is
defined to be the distance between fibers $\E_X$ and $\E_Y$ in $M$, as determined by a suitable Riemannian structure on $M$.
We choose the Riemannian metric
on $M=\Sop\times \D^+(p)$ to be a product metric determined by
bi-invariant Riemannian metrics $g_{\So}, \gdp$ on the two factors (each
of which is a Lie group).
The corresponding squared distance function
$d_M^2$ is a sum of squares. The geodesics connecting two fibers
$\E_X$ and $\E_Y$ with the minimal length give rise to \textit{minimal
  smooth scaling-rotation curves (MSSR) curves}, ``efficient''
scaling-rotation curves that join $X$ and $Y$.

The scaling-rotation framework has the potential to improve statistical analysis of SPD matrices in situations in which eigenstructure is fundamental.  Take, for example, a regression analysis of SPD-matrix-valued data.
Using scaling-rotation curves, one can explicitly model the changes of SPD matrices separately in terms of eigenvalues or eigenvectors. In the setting of DTI, this means that diffusion intensities and diffusion directions can be modeled individually or jointly.
Thus the changes of diffusion tensor (either along the fibers of tensors, or as a function of time or covariates) may be interpreted more meaningfully than is the case with some alternative frameworks.
In particular, we found in \cite{JSG2015scarot} that MSSR curves oftentimes exhibit  deformations of ellipsoids (representing SPD matrices) that are more natural to the human eye than are the deformations determined by the interpolation methods of  \cite{Arsigny2007,Pennec2006}; the summary measures of diffusion tensors ($3 \times 3$ SPD matrices) such as fractional anisotropy and mean diffusivity evolve in a regular fashion. Moreover, in the scaling-rotation framework, exploratory statistics such as mean, median, and principal components may carry high interpretability, again due to separability of eigenvalues and eigenvectors.
 The scaling-rotation framework carries over to SPD-matrix-valued data of higher dimensions, such as in dynamic-factor models concerning covariance matrices varying over time \cite{forni2000generalized}. Our computational algorithms for low dimensions are still applicable through dimension reduction; we leave such developments for future work.

\subsection{Overview of main results}\label{sec:intro-results}

 We carefully characterize the eigenvalue-based stratification of $\sympp$ in Section \ref{sect:strat}. We begin with identifying
 all
the fibers of the eigen-composition map $F$ systematically in terms of partitions of the integer $p$ and the set $\{1,2,\dots,p\}$.  This  culminates in Section \ref{fibstr} with a very explicit description of all the fibers.
In Sections \ref{sect:orbit-type-strat}-\ref{sec:4ss} we show how these ideas
lead to stratifications of $\sympp$. In
Section \ref{sec:fibex}, we explicitly describe all the strata and all the fiber-types for the cases for $p=2$ and $p=3$.

Understanding the stratification enables us to analyze some
non-trivial features of the scaling-rotation framework.  For example,
$\dsr$ is a metric on the top stratum of $\sympp$, but is not a metric
on all of $\sympp$.  For any $p$, the analysis of $\dsr(X,Y)$ and MSSR
curves from $X$ and $Y$ depend on the strata to which $X$ and $Y$ belong, because
fibers are topologically and geometrically different for different
strata.
In Section \ref{sec:dsrmssr} we  review the geometry of scaling-rotation framework. In Section~\ref{sect3.1},  we first introduce our choice of Riemannian metric $g_M$ on $M$, and define scaling-rotation curves in $\sympp$ as images of geodesics in $(M,g_M)$.
While the geometry of the ``upstairs" Riemannian manifold $(M,g_M)$ is relatively simple, the problem of determining MSSR curves between arbitrary $X$, $Y$ in the quotient space $\sympp$  is highly nontrivial, as is determining how the set of all such curves depends on $X$ and $Y$.
%
%
In Section \ref{sec:dsrmssr2}, we define scaling-rotation distance and MSSR curves, and in Section~\ref{sect:uniqq} we summarize results from \cite{GJS2016mathpaper} on general tools used in computing these objects. These results are applied to the important $p = 3$ case in Sections \ref{sec:dsr3} and \ref{sec:mssr3}.

 As we shall see, for any $X,Y \in \sympp$, an  MSSR curve from $X$ to $Y$ always exists, but need not be unique.
 This paper also characterizes when such a curve {\em is} unique, very explicitly for the cases $p=2$ and $p=3$.   Precisely describing the conditions of uniqueness is vital in any probability statement on random objects on $\sympp$. For example, for any two random objects $X$ and $Y$ drawn from continuous distributions defined on $\sympp$, with probability 1 there exists a unique MSSR curve between them.

Because all strata of $\sympp$ other than the top stratum have positive codimension, any random object $X$ drawn from a continuous distribution defined on $\sympp$ will lie in the top stratum with probability 1. Nonetheless, we cannot assume that a population-mean or parameter $\mu \in \sympp$ for a continuous distribution lies in the top stratum.  Therefore, with the possibility that for $\mu \in \symp(p)$,  $\mu$ does not have distinct eigenvalues, a closed-form expression for $\dsr(\mu,X)$, and  a systematic characterization and analysis of MSSR curves from $\mu$ to $X$, are desirable.
In this paper, we focus on the cases $p = 2$ and $p = 3$.

In Section~\ref{sec:appendix-MSSR2}, for $p = 2$, we provide closed-form expressions for the scaling-rotation distance, provide conditions on $X,Y \in \symp(2)$ under  MSSR curves between $X$ and $Y$ are unique, and illustrate the cases of uniqueness and non-uniqueness. (When there is not a unique MSSR curve from $X$ to $Y$, there are several possibilities for the number of MSSR curves from $X$ to $Y$.)


Sections \ref{sec:dsr3}--\ref{sec:appendix-MSSR3} are devoted to the case $p=3$.
 In Section \ref{sec:dsr3}, we use the quaternionic parametrization of $SO(3)$ to help us
characterize scaling-rotation distances, to evaluate closed-form expressions for the distances, and to identify and parameterize MSSR curves between $X, Y \in \symp(3)$.
In this section we also reduce the combinatorial complexity of these problems depending  on  the strata to which $X$ and $Y$ belong. A catalog of the ``nontrivial" unique and non-unique cases
of MSSR curves is given   in Section \ref{sec:6.1}, and a detailed algorithm for computing scaling-rotation distance and the set of MSSR curves is given in Section~\ref{algo}. 
In Section \ref{sec:appendix-MSSR3}, we schematically illustrate the conditions on $X,Y \in \symp(3)$ in the catalog of Section~\ref{sec:mssr3}, and provide some pictorial examples of unique and non-unique MSSR curves (including cases in which both $X$ and $Y$ lie in the top stratum; these cases are omitted from the
catalog in Section \ref{sec:mssr3}).

Some of the material in Sections~\ref{sect:strat} and \ref{sec:dsrmssr}   summarizes \cite{JSG2015scarot}, and especially, \cite{GJS2016mathpaper}.
However, particularly in Section~\ref{sect:strat}, for some topics we greatly expand upon \cite{GJS2016mathpaper}, including giving detailed descriptions and illustrations of fibers and strata.

 Frequently used notations and symbols are listed in Table~\ref{table:notations}.
\begin{table}[t]\caption{Frequently used notations and symbols. \label{table:notations}}
\begin{tabular}{ll}
  Notation     & Definition or description\\
  \hline
  $M = \Sop \times \diagp(p)$ &  the space of eigen-decompositions of $p\times p$ SPD matrices            \\
  $d_M$  & the geodesic distance function on $M$\\
  $F:M \to \sympp$ & the eigen-composition map\\
  $\E_X = F^{-1}(X)$ & the set of eigen-decompositions of $X$; fiber over $X$\\
  $\Comp(\E_X)$ & the set of connected components of $\E_X$ \\
  $\dsr(X,Y)$ & the scaling-rotation distance between $X,Y \in \sympp$ \\
  $\chi$ & a scaling-rotation curve in $\sympp$ \\
  $\Mc(X,Y)$ & the set of MSSR curves between $X,Y \in \sympp$ \\
  $\partpset$  & the set of partitions of  $\{1,2,\dots,p\}$ \\
  $\partp$   & the set of partitions of $p$ \\
  $G_D$   & the stabilizer group of $D \in \Dp$ under the action of $SO(p)$ on $\sympp$ \\
  $G_D^0$   & the identity component of $G_D$ \\
  $\J_D $      & the partition of $\{1,2,\dots,p\}$ determined by  $D \in \Dp$ \\
  $S_p$   & the permutation group of the set $\{1,2,\dots, p\}$ \\
  $\I_p$  & the  group  of sign-change matrices  \\
  $\I_p^+$ & the group  of  even sign-change matrices\\
  $\tilde{S}_p$ & the group of signed-permutation matrices \\
  $\tilde{S}_p^+$ & the group of even signed-permutation matrices \\
  $\J $      & a typical element  in $ \partpset$ \\
  $[\J] $     & a typical element in $ \partp$; projection of $\J$ under natural map \\
  $\Sc_{\J}$ & the stratum of $M$ labeled by $\J$ \\
  $\Sc_{[\J]}$ & the stratum of $\sympp$ labeled by $[\J]$\\
  $ \Sc_{\rm{top}} = \Sc_{[\J_{\rm{top}}]}$ & the top stratum of $\sympp$\\
  $ \Sc_{\rm{bot}} = \Sc_{[\J_{\rm{bot}}]}$ & the bottom stratum of $\sympp$\\
  $\Sc_{\rm{mid}} = \S_{[\J_{\rm mid}]}$ & the ``middle'' stratum of $\symp(3)$\\
  $\mathbf{H}$ & the space of quaternions \\
  $S_{\mathbf{H}}^3$ & the unit sphere in $\mathbf{H}$ \\
  $S_{\mathbf{C}}^1$ & the unit circle in $\mathbf{C}$, the complex plane \\
  $S_{\mathbf{C}^2}^3$ & the unit sphere in $\mathbf{C}^2$ \\
  $\phi:S^3_\H\to\So(3)$& the natural two-to-one Lie-group homomorphism (see Section \ref{sec:5.1.1.relhso3}) \\
  $\So(3)_{<\pi}$ &  the set of non-involutions in $\So(3)$ \\
  $s:\So(3)_{<\pi}\to S^3_\H$ & a smooth right-inverse to $\phi$ on $\So(3)_{<\pi}$
\end{tabular}
\end{table}

\setcounter{equation}{0}
\section{ Stratification of $\sympp$ }
\label{sect:strat}
%

\subsection{Partitions of $p$ and $\{1,2,\dots,p\}$}
\label{sec:partn}

We will consider several stratified spaces in this paper.  The strata
we define will be labeled by two different types of {\em partitions}.
For the sake of efficiency we first review these partitions
and fix some related notation.

Recall that a partition of the {\em positive integer} $p$ is a
(necessarily finite) sequence of positive integers $k_1\geq k_2\geq
k_3\geq \dots$ with $\sum k_i = p$, while a partition of the {\em set}
$\{1, 2, \dots, p\}$ is a  finite collection $\J=\{J_1,
J_2, \dots\}$ of one or more nonempty, pairwise disjoint subsets $J_i$
whose union is $\{1, 2, \dots, p\}$.
Partitions of an integer are
commonly written using additive notation, e.g. $2+2+1$ (a partition of 5).
In a partition $k_1+k_2+\dots$ of $p$,
the terms $k_i$ of the sequence are called the {\em parts} of the
partition (and are counted with multiplicity; the parts of $2+2+1$ are
2, 2, and 1).
 In a
partition $\J=\{J_1, J_2, \dots\}$, the $J_i$ are called the {\em
  blocks} of $\J$.



\begin{notation}$~$

{\rm 1. We write $\partpset$ for the set of partitions of
  $\{1,2,\dots,p\}$, and $\partp$ for the set of partitions of $p$.

 2. We write $S_p$ for the symmetric group (permutation group) of the set $\{1,2,\dots, p\}$.

3. The natural left-action of  $S_p$ on $\{1,2,
\dots, p\}$ induces left-actions of $S_p$ on $\partpset$ and $\bfr^p$,
given by
\bearray
\pi \dotprod \{J_1, J_2, \dots, J_r\} &=& \{\pi(J_1), \pi(J_2), \dots,
\pi(J_r)\}, \\
\pi\dotprod (x_1, x_2, \dots, x_p)
&=& (x_{\pi^{-1}(1)}, x_{\pi^{-1}(2)}, \dots,x_{\pi^{-1}(p)}),
\label{sprp}
\eearray

\noi where $\pi\in S_p$ and $\J\in\partpset.$ For $\J\in \partpset$,
we write $[\J]$ for its image in the quotient space $\partpset/S_p$.
}
\end{notation}

There is an obvious $S_p$-invariant map $\partpset\to\partp$ that
assigns to $\J=\{J_1,\dots,J_r\}$ the sequence $|J_1|, \dots, |J_r|,$
rearranged in nonincreasing order.  This map induces a bijection
$\partpset/S_p\to \partp$. {\em Henceforth we will use this bijection
  implicitly and will regard $\partpset/S_p$ and $\partp$ as the
  \underline{same} \underline{set}}; e.g. we will generally write
$[\J]$ for a typical element of $\partp$.


The sets $\partpset$ and $\partp$ are partially ordered by the {\em
  refinement} relation.  For $\J,\K \in\partpset$, we say that {\em
  $\K$ is a refinement of $\J$}, or that {\em $\K$ refines $\J$}, if
every element of $\K$ is a subset of an element of $\J$ (remember that
an {\em element} of $\K$ or $\J$ is a {\em subset} of
$\{1,\dots,p\}$); equivalently, if $\K$ can be obtained by
partitioning the elements of $\J$.  We write $\J\leq \K$ if $\K$
refines $\J$; ``$\leq$'' is then a partial ordering on
$\partpset$.
 Similarly, for $[\J], [\K]\in\partp$, we say that {\em
  $[\K]$ is a refinement of $[\J]$}, or that {\em $[\K]$ refines
  $[\J]$}, if $[\K]$ can be obtained by partitioning the parts of
$[\J]$. (For example, $\{3,2,2\}$ refines $\{7\}, \{5,2\}$, and
$\{4,3\}$, but neither of $\{3,3,2\}, \{4,1,1,1\}$ refines the other.)
We write $[\J]\leq [\K]$ if $[\K]$ refines $[\J]$; this ``$\leq$'' is
a partial ordering on $\partp$. Note that the quotient map
$\partpset\to \partpset/S_p = \partp$ is
order-preserving. These relations are illustrated for $p=3$ in Table
\ref{tab:Partition_Example}.

\begin{table}[tb]
\caption{Part($\{1,2,3\}$) and Part(3). The relation $\J_2 \leftarrow
  \jtop$ stands for $\J_2 \le \jtop$ (i.e., $\jtop$ refines $\J_2$).
  None of $\J_1$, $\J_2$, $\J_3$, refines either of the others. $[\jtop]$
  refines $[\J_2]$. \label{tab:Partition_Example}}
\begin{tabular}{lcccc}
\hline
\multicolumn{3}{c}{Partitions of $\{1,2,3\}$} & &
\multicolumn{1}{c}{Partitions of 3} \\
\cline{1-3} \cline{5-5}
$\jtop = \{\{1\},\{2\},\{3\}\}$ & & $\jtop$ &
$\longmapsto$ & $[\jtop] = 1+1+1$\\
$\J_3 = \{\{1,2\},\{3\}\}$ &  &
\reflectbox{\rotatebox[origin=c]{140}{$\leftarrow$}}
\reflectbox{\rotatebox[origin=c]{90}{$\leftarrow$}}
\reflectbox{\rotatebox[origin=c]{40}{$\leftarrow$}} & &
\reflectbox{\rotatebox[origin=c]{90}{$\leftarrow$}} \\
$\J_2 = \{\{1,3\},\{2\}\}$ & &
$\J_1\quad \J_2\quad \J_3$ & $\longmapsto$  & $[\J_1] =[\J_2] =[\J_3]= 2+1$\\
$\J_1 = \{\{2,3\},\{1\}\}$ & &
\reflectbox{\rotatebox[origin=c]{40}{$\leftarrow$}}
\reflectbox{\rotatebox[origin=c]{90}{$\leftarrow$}}
\reflectbox{\rotatebox[origin=c]{140}{$\leftarrow$}}    & &
\reflectbox{\rotatebox[origin=c]{90}{$\leftarrow$}} \\
$\jbot = \{\{1,2,3\}\}$ &  & $\jbot$ & $\longmapsto$ & $[\jbot] = 3$\\
\hline
\end{tabular}
\end{table}

For all partial-order relations ``$\leq$'' in this paper, the meanings
of the symbols ``$<$'', ``$\geq$'', and ``$>$'' are defined from
``$\leq$'' the obvious way.
%
%
Note that there is a well-defined largest (also called highest) and
smallest (also called lowest) element of $\partpset$ and of $\partp$:
for all $\J\in\partpset$, we have
\bestar
\jbot : = \{ \{1,2, \dots,p\} \}
\leq \ \ \J &\leq& \{ \{1\}, \{2\}, \dots, \{p\}\} =:\jtop, \\
p \leq \ \ [\J] &\leq& 1 + 1+ \cdots + 1.
\eestar


\subsection{Relation of partitions to eigenstructure}
Let $\D(p)$ denote the set of $p\times p$ diagonal matrices, and
recall that $\Dpp:=\{\diag(d_1, \dots, d_p) : d_i>0, 1\leq i\leq
p\} \subset \D(p)$.

\begin{defn}\label{defjd}{\rm
For $D=\diag(d_1,\dots, d_p)\in {\rm Diag}(p)$, let $\J_D$ denote
the partition of $\{1, 2, \dots, p\}$ determined by the
equivalence relation $i\sim_D j \iff d_i=d_j$.
}
\end{defn}

Various objects we can define that depend on $D$
actually depend only on the partition $\J_D$. As $D$ runs over all of
$\Dpp$, the partitions $\J_D$ run over all of $\partpset$.  For this
reason we define certain objects, such as the groups $G_\J$ below, in
terms of general partitions $\J$ of $\partpset$.

\begin{defn}\label{defgj}{\rm
For $\emptyset \neq J \subset \{1, 2, \dots, p\}$, let $\bfr^J\subset
\bfr^p$ denote the subspace \linebreak $\{(x_1,\dots, x_p)\in
\bfr^p\mid x_j = 0 \ \forall j\notin J\}$. For a partition $\J=\{J_1,
\dots, J_r\}$ of $\{1, 2, \dots, p\}$, let $\{W_1, \dots, W_r\} =
\{W_1^\J, \dots, W_r^\J\} = \{\bfr^{J_1}, \dots,
\bfr^{J_r}\}$ denote the corresponding subspaces of $\bfr^p$; note
that we have an orthogonal decomposition $\bfr^p = \bfr^{J_1}\plus
\dots \plus \bfr^{J_r}$. Define the subgroup $G_\J\subset\Sop$
by
\be\label{blockstab}
G_\J=\{R\in \Sop\mid RW_i = W_i, 1\leq i\leq r\},
\ee

\noi a Lie group with (generally) more than one connected
component. We write $G_\J^0$ for the identity component of $G_\J$ (the
connected component of $G_\J$ containing the identity).  }
\end{defn}

If each block $J_i$ consists of consecutive integers, then the
elements of $G_\J$ are block-diagonal. For example, if $p=5$ and
$\J=\{ \{1,2\}, \{3,4,5\}\}$, then
\bestar
G_\J &=&\left\{ \left[
\begin{array}{cc} R_1 & 0 \\ 0 & R_2
\end{array}\right] : R_1\in O(2),\ R_2\in O(3),\ \det(R_1)\det(R_2)=1
\right\}\\
&\iso& S(O(2)\times O(3)).
\eestar

\noi
 (For any
subgroup $H\subset O(p)$, we write $S(H)$ for $H\intersect \So(p)$.) In this example, $G_\J$ has two connected components, one in
which $\det(R_1)=\det(R_2)=1$ (the component $G_\J^0$), and one in
which $\det(R_1)=\det(R_2)=-1$.

For general $\J = \{J_1,\ldots,J_r\}$, the elements of $G_\J$ have ``interleaved blocks''. Writing $k_i = | J_i | $, we have
\be
G_\J \iso S(O(k_1)\times O(k_2) \dots \times O(k_r)),
\label{gjiso}
\ee
and the identity component $G_\J^0$ is isomorphic to $SO(k_1)\times SO(k_2) \dots \times SO(k_r)$.
If the $k_i$ are non-decreasing then $[\J]=k_1+\dots+k_r$. For  concreteness we define
\be\label{eq:gj0}
G_{[\J]}^0 = SO(k_1)\times SO(k_2) \dots \times SO(k_r)
\ \ \ \mbox{if}\ [\J]=k_1+\dots+k_r.
\ee

The groups $G_\J$ are also partially-ordered. For $\J, \K \in \partpset$,
$$ \J\leq \K \iff G_\J \supset G_\K. $$
This partial-ordering will be reflected in the stratifications of $\sympp$ and $M$ discussed in Section \ref{sec:4ss}.

%

\begin{defn}\label{defgd}{\rm
For each $D\in \Dp$, we define the stabilizer group of $D$,
$$
G_D=\{R\in \So(p) : RD=DR\} =\{R\in \So(p) : RDR^{-1}=D\}.
$$
}
\end{defn}

Note that if $D_1,D_2 \in \Dp$ have each distinct diagonal entries, then $G_{D_1} = G_{D_2}$.
In general, $G_D$ does not depend on the
absolute or relative sizes of the diagonal entries of $D$, but only on
which entries are equal to which others.
The stabilizer group is closely related to eigenstructure:
if $(U,D)\in M$ is an eigen-decomposition of $X\in \sympp$, then for any $R \in G_D$, $(UR,D)$ is also an eigen-decomposition of $X$. But $G_D$ is precisely the group $G_{\J_D}$ defined using Definitions \ref{defjd} and \ref{defgj}, and the identity components are related similarly: $G_D^0 = G_{\J_D}^0$.

\subsection{The groups of signed-permutation matrices}
\label{gpthy}

In this subsection we define two groups, $\tsp$ and $\tsp^+$, related to the stabilizer group of $D \in \Dp$.  Both extend the symmetric group
$S_p$, and we interpret these groups in terms of matrices. 

\begin{notation}$~$

{\rm
1. We write $\I_p$ for the group $(\ztwo)^p = \ztwo \times \ztwo
\times \dots \times \ztwo$ ($p$ copies).
 Each $\ztwo$ is the {\em   group of signs}  with elements $\pm
1$.  We write typical elements of $\I_p$ by $\boldsig=(\sigma_1,
\dots, \sigma_p)$.  We call $\I_p$ the {\em group of sign-changes},
and write $\bone$ for its identity element.

2. For $\pi\in S_p, \boldsig=(\sigma_1, \dots, \sigma_p)\in \I_p$, in
accordance with \eqref{sprp} we set
\be\label{pidsig}
\pi\dotprod \boldsig = (\sigma_{\pi^{-1}(1)}, \dots, \sigma_{\pi^{-1}(p)}).
\ee

}\end{notation}

%

Observe that $\sgn: \I_p\to \ztwo$ is indeed a homomorphism, and is
$S_p$-invariant:
\be\label{spinv}
\sgn(\pi\dotprod \boldsig) =\sgn (\boldsig)\ \ \mbox{for all}\
\pi\in S_p,\ \boldsig\in \I_p\ .
\ee

We also write ``$\sgn$'' for the usual sign-homomorphism $S_p\to
\ztwo$. Both of these sign-homomorphisms determine index-two
subgroups, the sets of elements of sign 1. For $\I_p$, our notation
for this subgroup will be
\ben
\I_p^+ :=\{\boldsig\in \I_p : \sgn(\boldsig)=1\}.
\een

\noi For $S_p$, of course, the corresponding subgroup is the
group of even permutations. 
By analogy, we call $\I_p^+$ the group of {\em even sign-changes}.

One may easily check that \eqref{pidsig} defines a left action of the
symmetric group on $\I_p$, and that $\mbox{``$\pi\dotprod$''}:
\I_p\to\I_p$ is an automorphism. Hence this action determines a {\em
  semidirect product} group: a group
\be\label{semi2}
\tsp = \I_p\rtimes S_p
\ee

\noi whose underlying set is $\I_p\times S_p$, and which contains
subgroups $\I_p\times\{\id\}$ and $\{\bone\}\times S_p$ isomorphic to
$\I_p, S_p$, respectively, but in which the group operation is given
by $
(\boldsig_1, \pi_1)(\boldsig_2, \pi_2) = (\boldsig_1 (\pi_1\dotprod
\boldsig_2),  \pi_1\pi_2).
$

%

Because of \eqref{spinv}, the sign-homomorphisms $\I_p\to \ztwo,
S_p\to\ztwo$ determine a third sign-homomorphism $\sgn: \tsp\to
\ztwo$, defined by
$
\sgn(\boldsig,\pi) = \sgn(\boldsig)\,\sgn(\pi).
$

\begin{defn}\label{deftspp}{\rm We write $\tsp^+$ for
$\ker(\sgn:\tsp\to \ztwo)$, an index-two subgroup of
    $\tsp$. Equivalently, $\tsp^+ = \{(\boldsig, \pi)\in\tsp :
    \sgn(\boldsig)=\sgn(\pi)\}.$}
\end{defn}

%

%

For later use, we record the orders (cardinalities) of the groups
$\tsp$ and $\tsp^+$:

\begin{result} The orders (cardinalities) of the groups $\tsp$ and
  $\tsp^+$ are as follows:
\be\label{gporder}
|\tsp| = 2^p p!\ , \ \ \ |\tsp^+| = 2^{p-1}p! \ .
\ee

\end{result}

\pf Immediate from \eqref{semi2} and the fact
that $\tsp^+$ has index 2 in $\tsp$.
\qedns

\begin{remark}\label{rem:gjcomps}{\rm
For $\sigma \in \ztwo=\{\pm 1\}$, let $O_\sigma(k)\subset O(k)$ denote
the set of orthogonal transformations with determinant $\sigma$. In the
setting of \eqref{gjiso}, the connected components of $G_\J$ are
$O_{\sigma_1}(k_1)\times O_{\sigma_2}(k_2)\times \dots \times O_{\sigma_r}(k_r)$,
subject to the restriction $\prod_i\sigma_i=1$.
Thus for each partition $\J$ with $r$ blocks, there is a 1-1 correspondence between the set of connected components of $G_\J$ and $\I_r^+$ (in which $(\sigma_1,\ldots,\sigma_r)$ lies).
This fact leads that the number of connected components is $2^{r-1}$, which is used in describing the fibers of $F$; see Proposition~\ref{prop:biject}.}
\end{remark}

The group $\tsp$ has a natural representation on $\bfr^p$, the map
$\mat: \tsp\to O(p)$ defined by
\be\label{ispp}
\mat(\boldsig,\pi) = I_\mbs P_\pi \ ,
\ee

\noi where $I_\mbs =\diag(\sigma_1, \dots, \sigma_p)$ and $P_\pi$ is
the matrix of the linear map
$\mbox{``$\pi\dotprod$''}:\bfr^p\to\bfr^p$ in \eqref{sprp}.   The
entries of the permutation matrix $P_\pi$ are $(P_\pi)_{ij} =
\d_{i,\pi(j)}$. (We will see shortly that $\mat$ is a homomorphism,
justifying the term ``representation on $\bfr^p$''.) It is easily seen
that $\mat$ is injective.

\begin{defn}{\rm
We call a $p\times p$ matrix $P$ a {\em signed-permutation matrix} if
for some (necessarily unique) $\pi\in S_p$ the entries of $P$ satisfy
$P_{ij}=\pm \d_{i,\pi(j)}$.  We call such $P$ {\em even} if
$\det(P)=1$ and {\em odd} if $\det(P)=-1$. (Note that evenness of $P$
is not the same as evenness of the associated permutation $\pi$.)
The set of signed $p\times p$ permutation matrices is exactly $\mat(\tsp)\subset O(p)$; the subset of even elements is exactly $\mat(\tsp^+)\subset \Sop$.

}
\end{defn}

It is easy to see that $\mat(\tsp)$ is actually a sub{\em group} of $O(p)$. (This also follows from the fact, shown below, that $\mat$  is a homomorphism $\tsp\to O(p)$.)  Furthermore, at the level of matrices, the sign-homomorphism $\tsp\to\ztwo$ is
simply determinant:
\be
\sgn(\boldsig,\pi)
=\det(\mat(\boldsig,\pi))
= \det(I_\mbs P_\pi).
\ee
\noi It follows that $\mat(\tsp^+)$ is a
subgroup of $\So(p)$.




Identifying $\Dpp$ with $(\bfr_+)^p\subset \bfr^p$, the action
  \eqref{sprp} yields an action of $S_p$ on $\Dpp$, given by
\bearray\label{act2}
\pi\dotprod \diag(\bfd) &=&\diag(\pi\dotprod \bfd\,)\\
\nonumber&=& \diag(d_{\pi^{-1}(1)}, \dots,d_{\pi^{-1}(p)})
\ \ \ \mbox{if}\ \bfd=(d_1, \dots, d_p).
\eearray

\begin{notation}{\rm
 For $D\in\Dpp$, we write $[D]$ for its image in the quotient
space $\Dpp/S_p$. }
\end{notation}

One may easily check that for any $\pi\in S_p$, $D\in\D(p)$, we have
\be\label{ppidotd}
\pi\dotprod D = P_\pi D (P_\pi)^T = P_\pi D (P_\pi)^{-1},
\ee

\noi 
and that the restrictions of the map $\mat$ to the subgroups
$\I_p\times\{\id\} \iso \I_p$ and $\{\bone\}\times S_p\iso S_p$ are
homomorphisms.  It follows easily that the map $\mat:\tsp\to
O(p)\subset GL(p,\bfr)$ is a homomorphism (hence a representation on
$\bfr^p$, as asserted earlier):
\bestar
\mat(\boldsig_1,\pi_1)\ \mat(\boldsig_2,\pi_2) &=&
I_{\mbs_1} \ \, P_{\pi_1} I_{\mbs_2} (P_{\pi_1})^{-1}\ \,
P_{\pi_1} P_{\pi_2}\\
&=&
I_{\mbs_1} (\pi_1\dotprod I_{\mbs_2}) P_{\pi_1\pi_2}\\
&=&
I_{\mbs_1} I_{\pi_1\dotprod \mbs_2} P_{\pi_1\pi_2}\\
&=&
I_{\mbs_1 (\pi_1\dotprod \mbs_2)} P_{\pi_1\pi_2}\\
&=&
\mat((\boldsig_1,\pi_1) (\boldsig_2,\pi_2)).
\eestar

\noi Since $\mat$ is an injective homomorphism, it is an isomorphism
onto its image, the subgroup $\mat(\tsp)\subset O(p)$.

%
%

%


As in \cite{JSG2015scarot}, we call the elements of
{ $\mat(\I_p)$}
{\em
  sign-change matrices}, even or odd according to their determinants.

For any subgroup $H$ of $\tsp$, $\mat$ restricts to an isomorphism
$H\to \mat(H)$.  Therefore to simplify notation, {\em henceforth in
  most expressions we will not write the map $\mat$
explicitly;} rather, we will use (for example) the
notation $\tsp$ for both $\tsp$ and
{$\mat(\tsp)$}. It should always be
clear from context whether our notation refers to an element (or
subgroup) of $\tsp$, or the corresponding matrix (or finite group of
matrices) under the map $\mat$.
However, to
avoid some odd-looking formulas we
will use the following notation:

\begin{notation}\label{defpig} {\rm We write typical elements
of the (abstract) signed-permutation group $\tsp$ as $g$, and define
the matrix $P_g=\mat(g)\in\tsp$. Thus $P_{(\mbs,\pi)}= I_\mbs P_\pi$. The
image of $g$ under the projection $\proj_2:\tsp\to S_p$ will be
denoted $\pi_g$.}
\end{notation}

We remark that if $\tsp$ is interpreted as
{$\mat(\tsp)$}, $\proj_2$ is the
map $\I_{\mbs}P_\pi\mapsto \pi$ (well-defined, since every
element of
{$\mat(\tsp)$} can be written {\em uniquely} in the form
$\I_{\mbs}P_\pi$).

Note that the action of $S_p$ on $\Dpp$ lifts to an action of $\tsp$
on $\Dpp$:
\be
g\dotprod D : = \pi_g\dotprod D.
\ee

\noi In terms of matrices, this is just the conjugation action:
\be
P_{(\mbs,\pi)}\dotprod D = I_\mbs P_\pi D P_\pi^{-1}
  I_\mbs^{-1}= P_\pi D P_\pi^{-1},
\ee

\noi the latter equality holding since sign-change matrices are
diagonal (and therefore commute with diagonal matrices).

\subsection{Structure of the fibers}\label{fibstr}

  We are now ready to provide a systematic description of the
fibers of $F$. We start with a result from  \cite{GJS2016mathpaper}:

\begin{prop}[{\cite[Corollary 2.6]{GJS2016mathpaper}}]
\label{cor:fibform4}
Let $X\in\sympp$ and $(U,D)\in \E_X$. Then
\be\label{fibform4}
\E_X= \{ (UR(P_g)^{-1}, \pi_g\dotprod D) : R\in G_D^0, g\in \tsp^+\}.
\ee
\end{prop}

 The fiber $\E_X$ generally has more than one connected component. The ``shape'' of the fiber $\E_X$ depends on the partition $[\J_D]$.

\begin{defn}\label{defkj}$~$

{\rm 1.
For $\J\in\partpset$,   $\Gamma_\J = \tsp^+\intersect G_\J$,
and $\Gamma_\J^0 = \Gamma_\J\intersect G_\J^0 = \tsp^+\intersect
G_\J^0$ .

2. For any $X\in\sympp$ and $(U,D)\in \E_X$, define
\be\label{defcompud}
[(U,D)]=\{(UR,D) : R\in G_D^0\},
\ee

\noi the connected component of $\E_X$ containing $(U,D)$. We write
$\Comp(\E_X)$ for the set of connected components of $\E_X$.

3. For any Lie group $G$ and closed subgroup $K$, we write $G/K$ and
$K\backslash G$ for the spaces of left- and right-cosets,
respectively, of $K$ in $G$. (In particular, we use this notation when
$G$ is a finite group.)
}
\end{defn}

The group $\tsp^+$ acts on $M = \So(p) \times \Dp$ via setting
$g\dotprod (U,D) = (UP_g^{-1}, g\dotprod D)$.
 This action preserves every fiber of $F$. Thus for each $X\in
\sympp$ there is an induced action of $\tsp^+$ on $\Comp(\E_X)$, given
by
\be\label{act3}
g\dotprod[(U,D)]=[g\dotprod(U,D)].
\ee
Each $g\in\tsp^+$, acting as above, permutes the connected components
of $\E_X$; the subgroup $\Gamma^0_{\J_D}$ is the stabilizer of $[(U,D)]\in\Comp(\E_X)$ under this action. 


\begin{prop}\label{prop:biject}\label{fibdescrip}
Let $X\in\sympp$.
\begin{enumerate}
 \item[(i)]
Then every $(U,D)\in\E_X$ determines a
  bijection between $\Comp(\E_X)$ and the set
  $\tsp^+/\Gamma^0_{\J_D}$.
  \item[(ii)] Let $(U,D)\in \E_X$, and $[\J_D] = k_1 + \cdots k_r$.
Then $\E_X$ is diffeomorphic to a disjoint union of
$2^{r-1}\frac{p!}{k_1! k_2! \dots k_r!}$ copies of $SO(k_1)\times
SO(k_2)\times \dots \times SO(k_r)$.
  \end{enumerate}
\end{prop}

 The proposition above is proved in \cite{GJS2016mathpaper}.

An important special case of Proposition \ref{prop:biject} is the case
in which all eigenvalues of $X$ are distinct.
In this case, $\J_D=\jtop = \{\{1\}, \{2\}, \dots, \{p\} \}$,
 $G_{\J_D} = S(O(1) \times O(1) \times \cdots \times O(1)) = \I_p^+$ and
 $G_{\J_D}^0 = SO(1) \times SO(1) \times \cdots \times SO(1)) = \{ I\}$.
 Thus $\Gamma^0_{\J_D} =\{\id\}$ and action of $\tsp^+$ on $\Comp(\E_X)$ is free as well as transitive.
Since $G_D^0=G_{\J_D}^0 = \{I\}$, each connected component of $\E_X$ is
a single point; $\Comp(\E_X)=\E_X$. Thus, by part (i) of the Proposition, any choice of $(U,D)\in \E_X$ yields a bijection $\tsp^+\to \E_X$, $g\mapsto g\dotprod (U,D)$.
Furthermore, $[\J_D] = \{1,1,\ldots,1\}$. Thus applying part (ii) of the Proposition, $\E_X$ is diffeomorphic to a disjoint union of $2^{p-1} p!$ copies of $SO(1)\times
SO(1)\times \dots \times SO(1)$, which is a point.

Examples of the fibers for $p =2,3$ can be found in Section~\ref{sec:fibex}.

\subsection{Orbit-type stratification of $\sympp$}\label{sect:orbit-type-strat}

 The compact Lie group $G=\Sop$ acts from the left on
the manifold $\sympp$ via
\be\label{actonsympp}
(U,X)\mapsto U\dotprod X= UXU^T.
\ee

\noi For each $X\in\sympp$, the orbit $G\dotprod X$ of $X$ is
diffeomorphic to $G/G_X$, where $G_X\subset G$ is the stabilizer
subgroup of $X$:
\be\label{stabx}
G_X:=\{ U\in G : UXU^T = X\}.
\ee

\noi If $Y\in G\dotprod X$ then $G_Y = UG_XU^{-1}$ for any
$U$ for which $Y=U\dotprod X$; hence $G_Y$ is conjugate to $G_X$.
More generally, whether or not $X,Y\in \sympp$ lie in the same orbit,
we say that $X$ and $Y$ have the {\em same orbit type} if the
stabilizers $G_X, G_Y$ are conjugate subgroups of $G$ (i.e. if $G_Y =
UG_XU^{-1}$ for some $U\in G$, an equivalence relation we will write
as $G_X\simc G_Y$).  If $X$ and $Y$ have the same orbit type then the
orbits $G\dotprod X, G\dotprod Y$ are diffeomorphic.  Define the {\em
  orbit-type stratum} of $\sympp$ associated with a given orbit-type
to be the union of all orbits of that type; we refer to the collection
${\bf S}$ of these strata as the {\em orbit-type stratification} of
$\sympp$.

The pair $(\sympp,{\bf S})$ is an example of a {\em Whitney stratified
  manifold}, one of several notions of ``stratified space'' in the
literature. In all such notions, a stratification of a topological
space $Z$ is a collection ${\bf S}$ of pairwise disjoint subsets of
$Z$, called {\em strata}, whose union is $Z$ and which are required to
satisfy certain conditions that depend on which notion of ``stratified
space'' is being used. The ``nicest'' type of stratification of a
manifold is a {\em Whitney stratification} \cite[Section
  1.1]{GorMac1988}. It is known that, for any compact Lie group acting on
a smooth manifold, the orbit-type stratification is a Whitney
stratification (\cite[p. 21]{Gibsonetal1976}).


However, not all the criteria for a Whitney stratification are
relevant to this paper.  Slightly modifying the terminology of
\cite{GorMac1988}, the notion of greatest relevance here is that of a {\em
  $\P$-decomposed space}, where $(\P,\leq)$ is a partially ordered
set. A {\em $\P$-decomposition} of a closed subset $Z$ of a manifold
$N$ is a locally finite collection
${\bf S}=\{\S_i\}_{i\in \P}$ of pairwise disjoint submanifolds of $N$
whose union is $Z$ and for which $\S_i\intersect
\Bar{\S_j}\neq\emptyset \iff \S_i\subset \Bar{\S_j} \iff i\leq j$,
where ``overbar'' denotes closure.  For the purposes of this paper, we
allow
``{\em stratified space}" to mean simply a $\P$-decomposition ${\bf S}$ of a
closed subset $Z$ of some manifold, where $\P$ is any partially
ordered set; the submanifolds $S_i$ are called the {\em strata}
of this stratification.
We will make pervasive use of
the ``$\P$-decomposition'' notion. Our $\P$ will always be either
$\partp$ or $\partpset$, and we will refer to it as a {\em label set}.

%
%

\subsection{Three equivalent stratifications of $\sympp$}
\label{sect:three-stratification}

There are three ``types'' that we will associate to each $X\in\sympp$.
The first, already defined, is the {\em orbit type} of $X$ under the
action \eqref{actonsympp}. The other two types, {\em fiber type} and
{\em eigenvalue-multiplicity type}, will be defined below.  For any
of these types, ``$X$ has the same type as $Y$'' is an equivalence
relation.  We will see that all three relations are identical.  Thus
the orbit-type stratification may be thought of just as well as a
fiber-type stratification or as an eigenvalue-multiplicity-type
stratification.

\begin{defn} {\rm

1.
We say that
$X,Y\in\sympp$ {\em have the same fiber type} if $[\J_D]=[\J_\L]$. In this case, the fibers $\E_X, \E_Y$ are diffeomorphic (cf. Proposition~\ref{prop:biject}).

2.
For $X\in \sympp$, we define the {\em eigenvalue-multiplicity
  type} of $X$, which we will denote $\ET(X)$, to be the multi-set of
multiplicities of eigenvalues of $X$ (the collection of eigenvalues of
$X$, enumerated with their multiplicities), an element of $\partp$.
}
\end{defn}

For example, if $p=3$, then for any $R_1,R_2,R_3\in\So(3)$, the
matrices
\be\label{ETexamp}
X_1= R_1 \left[ \begin{array}{ccc} 1 & 0 & 0 \\ 0& 6 & 0 \\ 0 & 0 & 6
\end{array} \right] R_1^{-1}, \ \
X_2=R_2 \left[ \begin{array}{ccc} 4 & 0 & 0 \\ 0& 1 & 0 \\ 0 & 0 & 1
\end{array} \right] R_2^{-1}, \ \
X_3=R_3 \left[ \begin{array}{ccc} 5 & 0 & 0 \\ 0& 7 & 0 \\ 0 & 0 & 5
\end{array} \right] R_3^{-1}
\ee
 in $\Sym^+(3)$ all have the same eigenvalue-multiplicity type,
the partition $2+1$ of 3.  The relative {\em sizes} of the eigenvalues
of $X\in\sympp$ have no bearing on the eigenvalue-multiplicity type of
$X$; all that matters are the eigenvalue {\em multiplicities}. As we shall see
later, in our stratification of $\D^+(p)$ the three diagonal matrices in \eqref{ETexamp} represent
two different strata.

The three ``types'' we have defined are conceptually different: For
$X\in \sympp$ and $D\in\D^+(p)$ for which $(U,D)\in\E_X$ for some
$U\in\Sop$,
\begin{itemize}

\item[(i)] the concept of {\em orbit-type} is based on (though not
  necessarily equivalent to) diffeomorphism type of the {\em orbit}
  $G\dotprod X$, a 
  submanifold of $\sympp$ diffeomorphic to
  $\Sop/G_D$; 

\item[(ii)] the concept of {\em fiber-type of $X$} is based on the diffeomorphism type of the {\em
    fiber} $\E_X$, a possibly non-connected submanifold of
  $\Sop\times\D^+(p)$ diffeomorphic to finitely many copies of $G_D$
  (the number of copies being the multinomial coefficient
  $\frac{p!}{k_1! \dots k_r!}$ appearing in Proposition
  \ref{prop:biject}); and 

\item[(iii)] the concept of {\em eigenvalue-multiplicity type} is based
  directly on discrete information: the partition $[\J_D]$ of $p$
  determined by the eigenvalues of $X$.

\end{itemize}
Even though the three kinds of ``types'' are conceptually
different, they are equivalent.

\begin{prop} \label{prop:types-strat}
{\rm
 For all $X,Y \in \sympp$,
\be\label{equivtypes}
 \mbox{same orbit-type  = same fiber-type = same
  eigenvalue-multiplicity type}.
\ee
}
\end{prop}

\pf
Let $X,Y\in \sympp$, $(U,D)\in\E_X$, and $(V,\L)\in\E_Y$. Note
that $G_X\simc G_D$ and $G_Y\simc G_\L$. Hence
\bestar
\mbox{$X,Y$ have the same orbit-type} &\iff & G_X\simc G_Y\\
&\iff& G_D\simc G_\L\\
&\iff& \J_D=\pi\cdot J_\L  \ \ \ \mbox{for some $\pi\in S_p$}\\
&\iff& [\J_D]=[\J_\L]\\
&\iff& \ET(X)=\ET(Y). \qedns
\eestar

%

Because of Proposition~\ref{prop:types-strat}, we are free to view the orbit-type
stratification as an eigenvalue-multiplicity-type stratification, and
to label strata accordingly.  We will do this in
Section \ref{sec:4ss}.

%
%
%

\subsection{Four stratified spaces}\label{sec:4ss}

As it is clear that the stratifications of $\sympp$ and $M = SO(p) \times \D^+(p)$ only depend on the eigenvalues, we also define stratifications of the spaces of eigenvalues: $\Dpp$  and $\Dpp/S_p$. Typical elements of $\Dpp/S_p$ will be denoted by $[D] = \{\pi\dotprod D : \pi \in S_p\}$. Strata of $\sympp$ (thus $\Dpp/S_p$) will be labeled by $\partp$; strata of $M$ and $\Dpp$ will be labeled by  $\partpset$.
The commutative
diagram in Figure~\ref{comdiag}, with notation as defined in Definition~\ref{defn:2.17}, indicates the relationships among these spaces and label-sets.
\begin{figure}[h]
\begin{diagram}
M = \So(p)\times\D^+(p) & \rTo^{\ \sproj} &\D^+(p) & \rTo^{\slbl\ \ \ } &
\partpset\\
\dTo >F  &                &             \dTo >{\quo_1} &  &
\dTo >{\quo_2}\\
\sympp                  & \rTo_{\bsproj\ \ \ \ \ \ } & \D^+(p)/S_p &
\rTo_{\ \ \ \ \bslbl} & \partp\\
\end{diagram}
\caption{Commutative diagram for the stratifications of $\sympp$ and related spaces. (This figure also appears on \cite{GJS2016mathpaper}.)}
\label{comdiag}
\end{figure}

\begin{defn}$~${\rm
\begin{itemize}
\item[(i)] $\proj_2:\So(p)\times\D^+(p)\to\D^+(p)$ is projection onto the
  second factor.

\item[(ii)] For $X\in\sympp$, if $(U,D)\in\E_X$ we define
  $\Bar{\proj_2}(X)=[D]\in\sympp/S_p$.

\item[(iii)] $\lbl: \D^+(p)\to\partpset$ is defined
  by $\lbl(D)=\J_D$.

\item[(iv)]  $\Bar{\lbl}: \D^+(p)/S_p\to\partp$ is
  defined by $\lbl([D])=[\J_D]$.

\item[(v)] $\quo_1$ and $\quo_2$ are the quotient maps $\D^+(p)\to
  \D^+(p)/S_p$ and \linebreak $\partpset\to \partpset/S_p=\partp$
  respectively.
\end{itemize}
}
\label{defn:2.17}
\end{defn}

The diagram suggests a natural definition of strata of the four spaces.

\begin{defn}\label{defn:4ss}
The four spaces $\Dpp, \Dpp/S_p, M$, and
 $\sympp$ are each stratified by strata labeled by stratum-labeling maps $\J\in\partpset$ and
$[\K]\in\partp$:
\bestar
\cald_\J := \lbl^{-1}(\J) &=&\{ D\in \Dpp : \J_D = \J\}\subset \Dpp, \\
\cald_{[\K]} := \Bar{\lbl}^{\,-1}([\K]) &=&\{ [D]\in \Dpp/S_p : \pi\dotprod D\in
\cald_{\K}\ \mbox{\rm for some}\ \pi\in S_p\}\\
&&\subset \Dpp/S_p \ ,
\\
\S_\J
:= {\rm proj}_2^{-1}(\cald_\J) &=& \Sop\times \cald_\J\\
&=& \{ (U,D)\in M : \J_D =\J\} \subset M, \\
\S_{[\K]} :=
\Bar{{\rm proj}_2}^{\,-1}(\cald_{[\K]})
&=& \{X\in \sympp \ \mbox{with  eigenvalue-multiplicity type}\ [\K]\}\\
&=& \{ X\in \sympp
: X=F(U,D)\ \mbox{for some}\\
&& \phantom{\{ X\in \sympp
: X} U\in \Sop, D\in \cald_\K\}\subset \sympp.
\eestar
\end{defn}

\begin{example}
The matrices $X_1, X_2, X_3$ in \eqref{ETexamp} all
 lie in the same stratum $\S_{[\J]}$ of $\sympp$, the one labeled by
 the partition $[\J]=2+1$ of 3. The diagonal matrices $D_1, D_2, D_3$
 appearing in the formulas in \eqref{ETexamp} for $X_1, X_2, X_3$,
 respectively, lie in two different strata of $M$: the first two lie
 in $\cald_{\J_1}$ while the third lies in $\cald_{\J_2}$, where
 $\J_1=\{\{2,3\},\{1\}\}$ and $\J_2=\{\{1,3\},\{2\}\}$. Note that
 strata need not be connected. For example, the stratum $\S_{2+1}$ in
 $\Sym^+(3)$ has two connected components, one in which the
 double-eigenvalue is the larger of the two distinct eigenvalues, and
 one in which it is the smaller. The matrix $X_1$ in \eqref{ETexamp}
 lies in the first of these components, while $X_2$ and $X_3$ lie in
 the second. The diagonal matrices $D_1$ and $D_2$ lie in different
 connected components of $\cald_{\J_1}$.
\end{example}

The map $\diag(d_1,\dots, d_p)\to (d_1,\dots,
  d_p)$ identifies $\D^+(p)$ diffeomorphically with $(\bfr_+)^p$. Under this identification,
for each $\J\in\partpset$ the stratum $\cald_\J$ is the intersection of a
linear subspace of $\bfr^p$ with the open subset $(\bfr_+)^p\subset \bfr^p$, hence is a submanifold of $(\bfr_+)^p$.  The stratum $\S_\J=\Sop\times \cald_J$ is therefore a submanifold of $M$.  The quotient $\D^+(p)/S_p$ is simply the
$p$-fold symmetric product of $\bfr_+$, which can be identified homeomorphically with
$Z:=\{(x_1,\dots, x_p)\in (\bfr_+)^p : x_1\leq x_2\leq \dots \leq x_p\}$, a closed subset
of $(\bfr_+)^p$.  This homeomorphism identifies the stratum $\cald_{[\J]}$ of $\D^+(p)/S_p$ with a submanifold of
$(\bfr_+)^p$ (diffeomorphic to a connected component of $\cald_\J$).  Thus our collections of strata of
$\D^+(p), \D^+(p)/S_p$, and $M$ meet our definition of ``stratified space".  As noted earlier, our stratification of $\sympp$ is an orbit-type stratification, hence automatically a Whitney stratification. Thus, each of $\Dpp,  M,$ $\Dpp/S_p$, and
$\sympp$, equipped with the strata defined above, is a stratified
space.
 Note also that for any $\J \in \partpset$,
$F(\S_\J) = \S_{[\J]}$ and $\quo_1 ( \cald_\J) = \cald_{[\J]}$.

If $\J, \K\in\partpset$ and $\K$ is a strict refinement of $\J$
(i.e. $\K$ refines $\J$ but $\K\neq \J$; equivalently, $\J<\K$), it is
easy to see that every element of the stratum $\S_{\J}$ in $M$ can
be obtained as a limit of a sequence lying in $\S_{\K}$, but that no
element of $\S_{\K}$ can be obtained as a limit of a sequence lying
in $\S_{\J}$ (in the limit of a sequence of matrices, distinct
eigenvalues can coalesce but equal eigenvalues cannot separate). Thus
\be\label{stratclo1}
\Bar{\S_\K} = \Union_{\J\leq \K} \S_\J,
\ee
where  $\Bar{\S}$ denotes the closure of a stratum $\S$.
A similar comment applies to strata $\cald_\J, \cald_\K$ in $\Dpp$; to strata
$\S_{[\J]}, \S_{[\K]}$ in $\sympp$; and to strata $\cald_{[\J]},
\cald_{[\K]}$ in $\Dpp/S_p$.
%



For any of the stratified spaces defined in Definition~\ref{defn:4ss}, the set of strata has a natural partial ordering, given by
 \bearray\label{ordstrat1}
\S_\J \subset \Bar{\S_\K} \iff \J\leq \K \iff \cald_\J \subset \Bar{\cald_\K} \ , \\
\S_{[\J]} \subset \Bar{\S_{[\K]}} \iff [\J]\leq [\K] \iff \cald_{[\J]} \subset \Bar{\cald_{[\K]}}\ .
\label{ordstrat2}
\eearray
  In each of the stratified spaces above, there is a highest stratum, corresponding to $\jtop$ and $[\jtop]$, and a lowest stratum, labeled by $\jbot$ and $[\jbot]$.
Note that for $\J, \K\in\partpset$, $\J\leq \K$ implies $[\J]\leq
[\K]$, but the converse is false for $p>2$.  In view of
\eqref{ordstrat1}--\eqref{ordstrat2}, a similar comment applies to $M$ and $\sympp$: $\S_\J\subset \Bar{\S_\K}$ implies $\S_{[\J]}\subset \Bar{\S_{[\K]}}$, but the converse  is false
for $p>2$. As a counterexample, set $\J = \{\{1,2\},\{3\}\}$ and $\K = \{\{1,3\},\{2\}\}$ for $p = 3$.

%
%
%

\begin{remark}[Number of strata]\label{numstrat}
{\rm The number of strata of $\sympp$ is the number of partitions of
  $p$, while the number of strata of $M$ is the number of partitions
  of $\{1,\dots, p\}$. In number theory, the {\em partition function}
  is the function that assigns to each positive integer $n$ the number
  partitions of $n$.
The number
  of partitions of $\{1,\dots,n\}$ is known as the {\em $n^{th}$ Bell
    number}.
Both the partition function and the Bell numbers
  have a long history and have been extensively studied; see
  \cite[Chapter XIX]{HarWri1968}, \cite{Stanley1997}.

}
\end{remark}

\begin{remark}[Dimensions of strata]\label{dimstrat} {\rm The dimensions of the strata
    in each of the four stratified spaces in diagram \eqref{comdiag} can
    easily be worked out; we will simply state the answers. If
    $\J\in\partpset$ and $[\J]=(k_1,\dots, k_r)$, then
\bestar
\dim(\S_\J) &=& r+ \dim(\Sop) = r+\frac{p(p-1)}{2}\ ,\\
\dim(\S_{[\J]}) &=& r+ \left(\dim(\Sop/G_\J)\right) =
r+ \left(\dim(\Sop)-\dim(G_\J^0)\right)\\
 &=& r+\frac{p(p-1)}{2} - \sum_{i=1}^r\frac{k_i(k_i-1)}{2}\ , \ \ \ \mbox{and}\\
\dim(\cald_\J) &=& \dim(\cald_{[\J]}) = r.
\eestar

}
\end{remark}


\subsection{Examples}\label{sec:fibex}

Using Proposition \ref{prop:biject}, Definition~\ref{defn:4ss} and Remarks \ref{dimstrat} and \ref{numstrat},
for any given $p$ we can, in principle, describe all the fibers of $F$ and all the strata of $M$ and $\sympp$ very explicitly.  As $p$ grows, the number of strata and the number of diffeomorphism-types of fibers grows rapidly, so below we do this
exercise only for the cases $p=2$ and $p=3$.

 \subsubsection{Example: $\Sym^+(2)$}\label{sec:exsr2}

\begin{figure}[tb]

  \centering

  \includegraphics[width=1\textwidth]{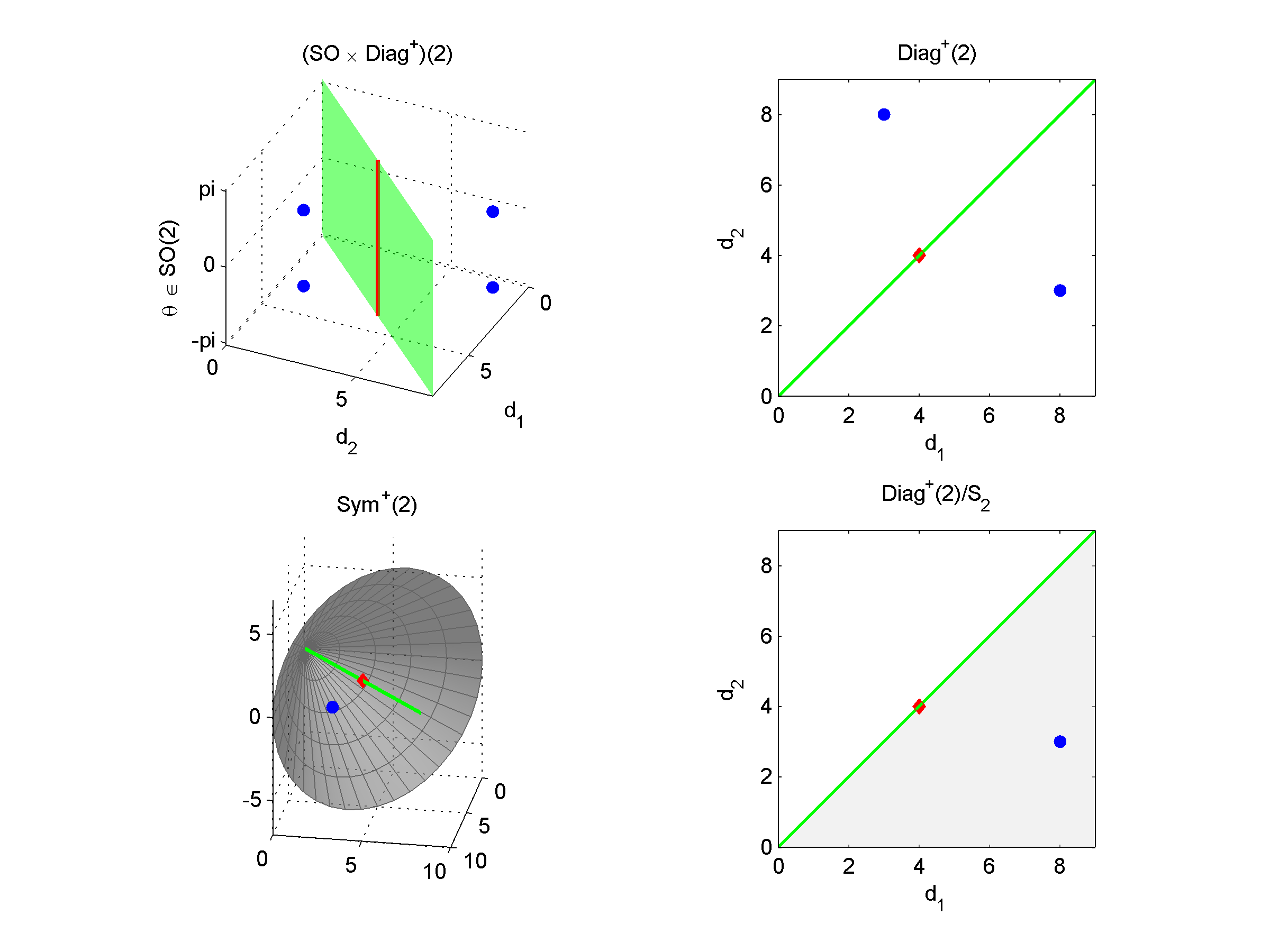}

   \caption{{\small Stratification of $M
    =\So(2) \times \mbox{Diag}^+(2)$ (top left), $\mbox{Sym}^+(2)$ (bottom left),
    $\mbox{Diag}^+(2)$ (top right) and $ \mbox{Diag}^+(2) / S_2$
    (shaded area in bottom right). $\Sym^+(2)$ is embedded in $\Sym(2)\iso \bfr^3$
as the cone $\{(a_{11},a_{22},\sqrt{2}a_{12}) : a_{11}>0, a_{22}>0, a_{11}a_{22}  - a_{12}^2 >0\}$.
    The space $ \D^+(2) / S_2$ is represented as a {\em fundamental domain} for the action of $S_2$ on $\D^+(2)$ (a subset $A\subset \D^+(2)$ containing exactly one point of each orbit, and such that  $\quo_1: \D^+(2)\to \D^+(2)/S_2$ restricts to a homeomorphism $\quo_1^{-1}(A\setminus \partial A)\to A\setminus\partial A$).
Also shown are $X = \mbox{diag}(8,3) \in \S_{[\J_{\rm top}]}
   \subset \mbox{Sym}^+(2)$ (blue dot) and
$Y = {\rm diag}(4,4) \in
    \S_{[\J_{\rm bot}]}$
(red dot), as well as their pre-images in $M$.
The projections to $\D^+(2)$ and $\D^+(2)/S_2$ of these subsets of $M$ and
$\Sym^+(2)$ are illustrated as correspondingly-colored dots in the right-hand
panels.
    }\label{fig:Symp2stratification}}

    \end{figure}

 \textbf{Stratification of $M$ and $\D^+(2)$}. There are two strata of
 $M =\linebreak (\So \times \D^+)(2)$ (and of $\D^+(2)$), labeled by the two
 partitions of $\{1,2\}$: $\J_{\rm top} = \{\{1\},\{2\}\}$ and
 $\J_{\rm bot } = \{\{1,2\}\}$.
\begin{itemize}
\item[(a)] The two-dimensional stratum $\cald_{\J_{\rm top}}$ consists of two
  connected components, $\{\mbox{diag}(d_1,d_2) : d_1> d_2 > 0\}$ and
  $\{ \mbox{diag}(d_1,d_2) : d_2> d_1 > 0\}$. Correspondingly, the
  three-dimensional stratum $\S_{\J_{\rm top}} = \So(2) \times \cald_{\J_{\rm top}}$
  also has two connected components.
\item[(b)] The one-dimensional stratum $\cald_{\J_{\rm bot }}$ is the connected set
  $\{ \mbox{diag}(d_1,d_2) : d_1 = d_2 > 0\}$. Therefore the two-dimensional stratum
  $\S_{\J_{\rm bot}} = \So(2) \times \cald_{\J_{\rm bot}}$ is also
  connected.
\end{itemize}
In the top panels of Fig. \ref{fig:Symp2stratification},
$\S_{\J_{\rm bot }}$ (respectively, $\cald_{\J_{\rm bot }}$) is
schematically depicted as the green plane (resp., line), which
separates the two connected components of $\S_{\J_{\rm top }}$ (resp.,
$\cald_{\J_{\rm top }}$).

 \textbf{Stratification of $\Sym^+(2)$ and $\D^+(2) / S_2$}. There are
 two strata of $\Sym^+(2)$ (and of $\D^+(2) / S_2$), corresponding to
 the two partitions of 2: $[\J_{\rm top }] = 1+1$, and $[\J_{\rm bot
 }] =2$. It is easily checked that for any $\J \in \partpset$,
 $F(\S_\J) = \S_{[\J]}$, $\quo_1 ( \cald_\J) = \cald_{[\J]}$.
 \begin{itemize}
\item[(a)] The top stratum $\S_{[\J_{\rm top}]} = \S_{1+1}$ is three-dimensional
and consists
  of SPD matrices with two distinct eigenvalues. Unlike $\S_{\J_{\rm
      top}}$ in $M$, the stratum $S_{[\J_{\rm top}]}$ in $\Sym^+(2)$
  is connected. In the bottom left panel of Figure
  \ref{fig:Symp2stratification}, $\S_{[\J_{\rm top }]}$ corresponds to
  the inside of the cone, minus the green line.
\item[(b)] The bottom stratum $\S_{[\J_{\rm bot}]} = \S_{2}$ is
one-dimensional and consists
  of SPD matrices with only one distinct eigenvalue.  This stratum is depicted
  as the green line in Fig. \ref{fig:Symp2stratification}.
\end{itemize}

\textbf{Fibers of $X \in \Sym^+(2)$}. The fibers are characterized by
Corollary~\ref{fibdescrip}.
 \begin{itemize}
\item[(a)] For any $ X \in \S_{[\J_{\rm top}]} = \S_{1+1}$, the fiber
  $\E_X \subset\S_\jtop\subset M$ consists of four points.
\item[(b)] For any $ X \in \S_{[\J_{\rm bot}]} = \S_{2}$, the fiber
  $\E_X \subset \S_{\jbot} \subset M$ is diffeomorphic to a circle.
An example of this circle is depicted schematically as the red line
segment in the top left panel of Fig.~\ref{fig:Symp2stratification}.
\end{itemize}

 \subsubsection{Example: $\Sym^+(3)$} \label{sec:exsr3}

\begin{figure}[tb]

  \centering

  \includegraphics[width=1\textwidth]{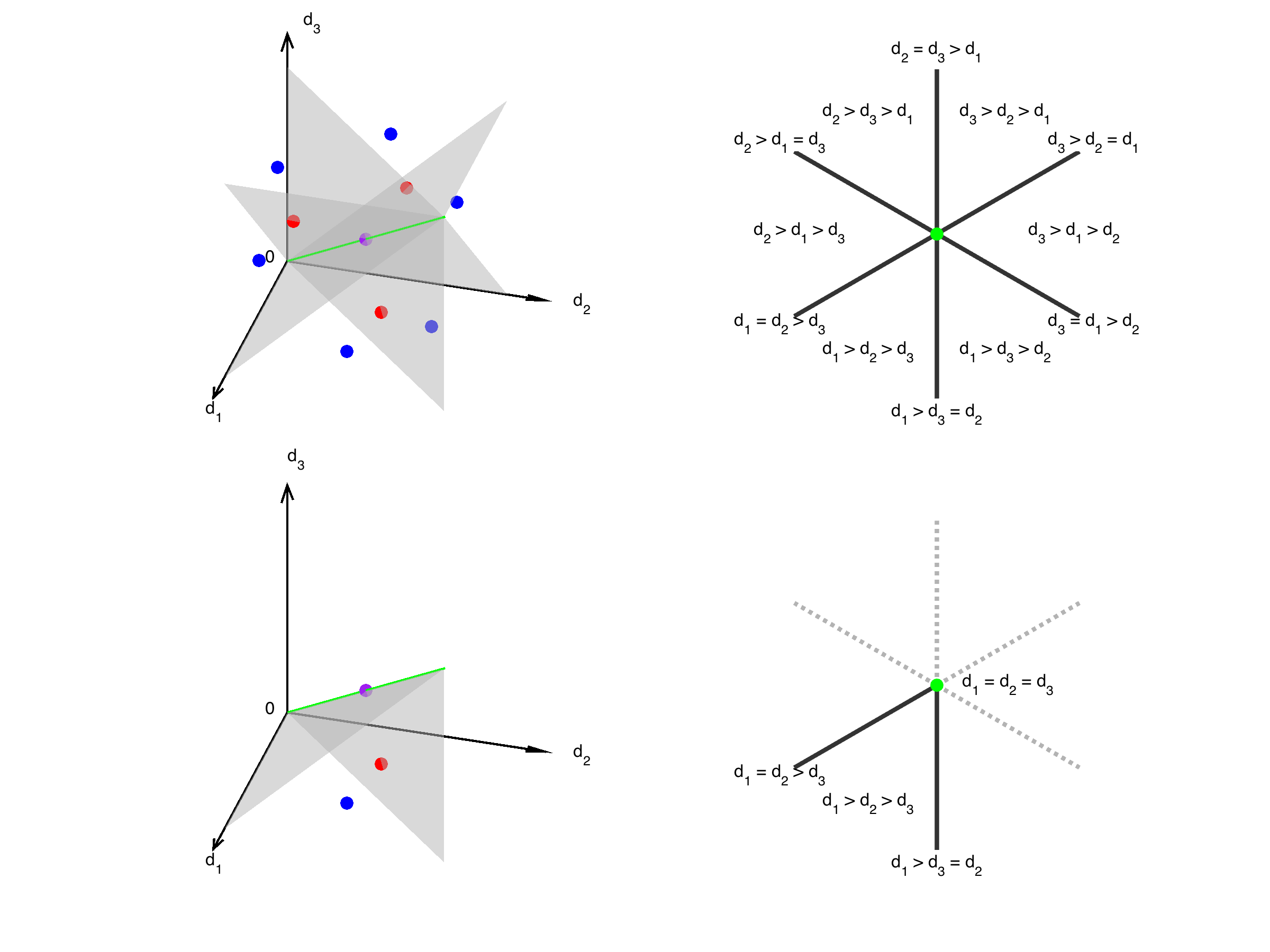}

  \caption{{\small Stratification of $\D^+(3)$ (top left) and $
      \mbox{Diag}^+(3) / S_3$ (bottom left). The star-figure (top right)
      represents the intersection of the top left figure and a hyperplane
      orthogonal to the green line.  In the bottom left panel, the
      space $ \D^+(3) / S_3$ is represented as a fundamental domain
      for the action of $S_3$ on $\D^+(3)$, a  convex cone bounded by below
by the positive quadrant of the $d_1d_2$-plane and on the sides by
      the two indicated gray planar sectors. Since the spaces
      $\Sym^+(3)$ and $M =(\So \times \D^+)(3)$ are six-dimensional,
      there are no simple visualizations of them.  Also shown are
      projections of $X = \mbox{diag}(8,5,1) \in \S_{[\jtop]} \subset
      \Sym^+(3)$ (blue dot), $Y = \mbox{diag}(6,6,2) \in \S_{\J_{\rm mid}}$ (red
      dot), and $Z = \mbox{diag}(5,5,5) \in \S_{[\jbot]}$ (green dot) to
      $\D^+(3) / S_3$ (bottom) and the pre-images of the quotient map
      in $\D^+(3)$ (top).}\label{fig:Symp3stratification}}

\end{figure}

 \textbf{Stratification of $M$ and $\D^+(3)$}. There are six strata of
 $M = (\So \times \D^+)(3)$ (and of $\D^+(3)$), labeled by the six
 partitions of $\{1,2,3\}$; see Table~\ref{tab:Partition_Example}. The
 features of the stratum $\cald_{\J}$ we discuss below apply also to the
 corresponding stratum $\S_{\J} = \So(3) \times \cald_{\J}$.
 \begin{itemize}
 \item[(a)] The stratum $\cald_{\J_{\rm bot}}$ is the connected
   component $\{\diag(d_1,d_2,d_3) : d_1 = d_2 = d_3\}$. In the top
   panel of Fig.~\ref{fig:Symp3stratification}, $\cald_{\J_{\rm
       bot}}$ corresponds to the green line.
 \item[(b)] The stratum $\cald_{\J_1}$ consists of two connected
   components: $\cald_{\J_1}^{\rm pro} = \linebreak
   \{\diag(d_1,d_2,d_3) : d_1 > d_2 = d_3\}$ and $\cald_{\J_1}^{\rm
     ob} = \{\diag(d_1,d_2,d_3) : d_1 < d_2 = d_3\}$. (The superscripts
     ``pro'' and ``ob'' stand for ``prolate'' and ``oblate'', respectively;
     see below.)
     The closures of
   these two connected components intersect in $\cald_{\J_{\rm
       bot}}$. In Fig.~\ref{fig:Symp3stratification}, $\cald_{\J_1}$
   corresponds to one of the three shaded planes except the green
   line.  The stratum $\S_{\J_1}$ also consists of two connected
   components: $\S_{\J_1}^{\rm pro} = \So(3) \times \cald_{\J_1}^{\rm
     pro}$ and $\S_{\J_1}^{\rm ob} = \So(3) \times \cald_{\J_1}^{\rm
     ob}$.  The strata $\cald_{\J_i}$ and $\S_{\J_i}$ for $i = 2$ and
   $3$ are similarly characterized.
\item[(c)] The stratum $\cald_{\J_{\rm top}}$ consists of six
  connected components, which can be labeled by permutations of
  $\{1,2,3\}$. Precisely, $\cald_{\J_{\rm top}} = \bigcup_{\pi \in
    S_3} \cald_{\J_{\rm top}}^\pi$, where $\cald_{\J_{\rm top}}^\pi =
  \{ \diag(d_1,d_2,d_3) : d_{\pi^{-1}(1)} > d_{\pi^{-1}(2)} >
  d_{\pi^{-1}(3)}\}$ for $\pi \in S_3$.
 \end{itemize}

 \textbf{Stratification of $\Sym^+(3)$ and $\D^+(3) / S_3$}.  There
 are three strata of $\Sym^+(3)$ (and of $\D^+(3) / S_3$),
 corresponding to the three partitions of 3: $[\J_{\rm top}]=1+1+1,
 [\J_{\rm mid}]=2+1,$ and $[\J_{\rm bot}]=3$.  These stratifications
 are closely related to an ellipsoid classification. An SPD matrix $X
 \in \Sym^+(3)$ with eigenvalues $a\ge b\ge c$ corresponds to the
 ellipsoid given by the equation $x^TX^{-1}x = 1$, and has the shape
 of a sphere if $a = b = c$, an oblate spheroid if $a = b > c$, a
 prolate spheroid if $a > b = c$, or a tri-axial ellipsoid if $a > b >
 c$. We will say that $X\in\Sym^+(3)$ is prolate (respectively oblate,
 triaxial) if the corresponding ellipsoid is a prolate spheroid (resp.
 oblate spheroid, triaxial ellipsoid).

 \begin{itemize}
\item[(a)] The stratum $\S_{[\J_{\rm top}]} = \S_{1+1+1}$ the
set of all SPD matrices with three distinct
  eigenvalues. Every $X \in \S_{[\J_{\rm top}]}$ is tri-axial. In the
  bottom panel of Fig.~\ref{fig:Symp3stratification}, the
  corresponding stratum $\cald_{[\J_{\rm top }]}\subset \D^+(3)$
  is depicted as an open convex
  cone. $\S_{[\J_{\rm top}]} $ is connected.
\item[(b)] The stratum $\S_{[\J_{\rm mid}]} = \S_{2+1}$ consists of
  SPD matrices with just two distinct eigenvalues, and is a disjoint
  union of two connected components: $\S_{2+1}^{\rm pro} =
  F(\S_{\J_1}^{\rm pro})$ and $\S_{2+1}^{\rm ob} = F(\S_{\J_1}^{\rm
    ob})$. If $X \in \S_{2+1}^{\rm pro}$, then $X$ is prolate; if $X
  \in \S_{2+1}^{\rm ob}$, then $X$ is oblate. Likewise, the stratum
  $\cald_{2+1}$ is a disjoint union of two connected components:
  $\cald_{2+1}^{\rm ob}$ and $\cald_{2+1}^{\rm pro}$.  In the bottom
  left panel of Figure~\ref{fig:Symp3stratification}, the two gray open
  planar sectors represent these two connected components of $\cald_{2+1}$.
\item[(c)] The  stratum $\S_{[\J_{\rm bot}]} = \S_{3}$ is the
  set of all SPD matrices with only one distinct
  eigenvalue.  The corresponding ellipsoids have the shape of
  a sphere. $\S_{[\J_{\rm bot}]}$ is connected.
\end{itemize}

\textbf{Fibers of $X \in \Sym^+(3)$}.
 \begin{itemize}
\item[(a)] For any $ X \in \S_{[\J_{\rm top}]} = \S_{1+1+1}$, the
  fiber $\E_X \subset M$ consists of 24 points, all of which lie in
  $\S_{\J_{\rm top}} \in M$.
\item[(b)] For any $ X \in \S_{[\J_{\rm mid}]} = \S_{2+1}$, the fiber
  $\E_X \subset M$ is diffeomorphic to 6 copies of the circle.
\item[(c)] For any $ X \in \S_{[\J_{\rm bot}]} = \S_{3}$, the fiber
  $\E_X \subset \S_{\J_{\rm bot}} \subset M$ is diffeomorphic to (one
  copy of) $\So(3)$ (and thus to $\bfr P^3$).
\end{itemize}
In Figure~\ref{fig:Symp3stratification}, examples of the three types
of fibers are provided. To help visualize, we show the projected
fibers, $\proj_2(\E_X) \subset \D^+(3)$. (For any $X \in
\S_{k_1+\cdots+k_r} \subset \Sym^+(p)$, $\proj_2(\E_X)$ is a discrete
set of cardinality $\frac{p!}{k_1! k_2! \ldots k_r!}$.)

\setcounter{equation}{0}
\section{ Scaling-rotation framework for curves and distances on $\symp(p)$}\label{sec:dsrmssr}



\subsection{Smooth scaling-rotation curves}\label{sect3.1}

 The space of eigen-decompositions $M = \Sop\times \D^+(p)$ is a Riemannian manifold. We define the Riemannian metric $g_M$ as a product Riemannian metric determined by metrics on $\Sop$ and $\D^+$ as follows.

The Lie algebra $\sop=T_I(\Sop)$ is the space of $p\times p$
antisymmetric matrices.
In this paper, for $U\in \Sop$ we identify the tangent space $T_U(\Sop)$ with the right-translate of $\sop$ by $U$:
\be\label{eq:Tu_sop}
T_U(\Sop) = \{AU : A\in\sop\}. 
\ee
The space $\D^+(p)$ is also a Lie group, but since it
is an open subset of a vector space, namely $\D(p)$, we make the identification
$T_D(\D^+(p))=\D(p)$ for all $D\in\D^+(p)$.
The tangent space of $M$ at $(U,D) \in M$ is $$ T_{(U,D)}M = T_U(\Sop)\plus
T_D(\D^+(p)).$$

 Using (\ref{eq:Tu_sop}), the standard bi-invariant Riemannian metric
$g_{\So}$ on $\Sop$ is defined by
\be\label{gsop}
\left.g_{\So}\right|_U(A_1,A_2)=-\frac{1}{2}\tr(A_1 U^{-1}A_2U^{-1}),
\ee
where $U\in \Sop$ and $A_1,A_2\in T_U(\Sop).$
 A (bi-)invariant Riemannian metric $\gdp$ is defined by setting
\be\label{gdp}
\left.\gdp\right|_D(L_1,L_2)=\tr(D^{-1}L_1 D^{-1}L_2)
\ee
 for $D\in \D^+(p)$ and $L_1,L_2\in T_D(\D^+(p)).$  Up to a constant factor, $\gdp$ is the only bi-invariant Riemannian metric on $\D^+(p)$ that is also invariant under the action of the symmetric group $S_p$.
  The product Riemannian metric is determined by the metrics above. Specifically, for $(U,D) \in M$ and $(A_i,L_i)\in  T_{(U,D)}M$
 we set
\be\label{gm}
\left.g_M\right|_{(U,D)}((A_1,L_1), (A_2,L_2)) =
k\left.g_{\So}\right|_U(A_1,A_2)
+
\left.\gdp\right|_D(L_1,L_2),
\ee
where  $k > 0$ is an
arbitrary parameter that can be adjusted as desired for applications.
Since the metrics $g_{\So}$ and $\gdp$ are
bi-invariant, the geodesics in $M$ can be obtained as either
left-translates or right-translates of geodesics through the identity
$(I,I)$.  In this paper, the right-translates are more convenient, which is why we have chosen the identification (\ref{eq:Tu_sop}) of the tangent spaces of $\Sop$.

\begin{defn}\label{defssr}
 {\rm A {\em smooth scaling-rotation (SSR) curve} is a
    curve $\chi$ in $\sympp$ of the form $F\circ\g$, where $\g:I\to M$
    is a geodesic defined on some interval $I$. }
\end{defn}


\begin{notation}{\rm
For $(U,D)\in M$, $A\in\sop$, and $L\in\D(p)$,
we define $\g_{U,D,A,L}:\bfr\to M$ and $\chi_{U,D,A,L}:\bfr\to \sympp$
by
\be\label{geod}
\g_{U,D,A,L}(t)=(\exp(tA)U, \exp(tL)D)
\ee

\noi and
\be\label{defchi}
\chi_{U,D,A,L}=F\circ\g_{U,D,A,L} .
\ee

We use the same notation $\g_{U,D,A,L},$ $\chi_{U,D,A,L}$ for the
restrictions of the curves above to any interval.
}
\end{notation}

The curve $\g_{U,D,A,L}:\bfr\to M$ is the geodesic in $M$ with
initial conditions
$\g(0)=(U,D)$, $\g'(0)=(AU,DL)\in T_{(U,D)}M$. The   curve
$\chi_{U,D,A,L}$ in $\symp(p)$, is the corresponding SSR
curve.

\subsection{
Scaling-rotation distance and MSSR curves}\label{sec:dsrmssr2}

Recall that in any group $G$, an
element $g$ is called an {\em involution} if $g^2$ is the identity element $e$  but $g\neq e$. Thus $R\in \Sop$ is an involution if $R^2=I\neq R$. That is, involutions in $SO(p)$ are reflections.
The
cut-locus of the identity in $\Sop$ is precisely the set of all involutions.  For
every non-involution $R\in\Sop$, there is a unique $A\in\sop$ of
smallest norm such that $\exp(A)=R$;
we define $\log(R)=A$.
If
$R$ is an involution, there is more than one smallest-norm $A\in\sop$
such that $\exp(A)=R$, and we allow $\log(R)$ to denote the {\em set}
of all such $A$'s. However, all elements $A$ in this set have the same
norm, which we write as $\|\log(R)\|$, where  $\| \ \|$
denotes the Frobenius norm on matrices: $\|A\|^2=\|A\|_F^2=\tr(A^TA)$.
  Thus $\|\log(R)\|$ is a
well-defined real number for all $R\in\Sop$, even when $\log(R)$ is
not a unique element of $\sop$. The geodesic-distance function $d_M$ on $M$ is then
\bearray
\label{defdM}
d_M^2\left((U,D),(V,\Lambda)  \right)
& = & k\, \dso(U,V)^2 + \ddp(D, \Lambda)^2
\\
&=& \frac{k}{2}\left\|\log (U^{-1}V)\right\|^2 +
\left\|\log( D^{-1}\L)\right\|^2.
 \nonumber 
\eearray

\begin{defn}[{\cite[Definition 3.10]{JSG2015scarot}}]{\rm For $X,Y\in \sympp$,
        the {\em scaling-rotation distance} $\dsr(X,Y)$ between $X$
        and $Y$ is defined by

\be\label{defdsr} \dsr(X,Y) := \inf_{ \substack{(U,D) \in \E_X,
    \\ (V,\Lambda) \in \E_Y }} d_M( (U,D), (V,\L) ).
\ee
}
\end{defn}

In \cite{JSG2015scarot}, $\dsr(X,Y)$ is interpreted as  ``the minimum amount of rotation and scaling needed to deform $X$ into $Y$.'' In the following, we provide an equivalent definition of $\dsr(X,Y)$ as the minimum length of SSR curves from $X$ to $Y$. However, we have not defined a Riemannian metric on $\sympp$, so there is no ``automatic" meaning attached to the phrase {\em length of a smooth curve} in $\sympp$.

\begin{defn} \label{defmssr}
{\rm Let $\g$ be a piecewise-smooth curve in $M$ and let
    $\ell(\g)$ denote the length of $\g$.
    \begin{itemize}
    \item[(i)] For $X,Y\in\sympp$,
    $\g:[0,1]\to M$ is called an {\em $(\E_X,\E_Y)$-minimal geodesic} if
    $\g(0)\in \E_X, \g(1)\in \E_Y$, and $\ell(\g)=\dsr(X,Y)$.
    \item[(ii)]  A pair of points $((U,D),(V,\L))\in \E_X\times \E_Y$ is called a {\em minimal pair} if $(U,D)=\g(0)$ and $(V,\L)=\g(1)$ for some
    $(\E_X,\E_Y)$-minimal geodesic $\g$.
    \item[(iii)] A {\em minimal smooth  scaling-rotation} (MSSR) curve from $X$ to $Y$ is a curve $\chi$ in $\sympp$ of the form $F\circ\g$ where $\g$ is an
    $(\E_X,\E_Y)$-minimal geodesic.
    We say that the MSSR curve $\chi=F\circ\g$ {\em corresponds to} the minimal pair formed by the endpoints of $\g$.
   \item[(iv)] The set of (not necessarily unique) MSSR curves from $X$ to $Y$ is denoted by $\M(X,Y)$.
    \item[(v)] For an SSR  curve $\chi$ in $\sympp$ we define the {\em length of $\chi$} to be $\ell(\chi):=\inf\{\ell(\g): \g \ \mbox{is a geodesic in $M$ and}\ F\circ\g =\chi\}$.
    \end{itemize} }
\end{defn}

Definition~\ref{defmssr}(i) also suggests the obvious fact that an $(\E_X,\E_Y)$-minimal geodesic is a minimal geodesic in
the usual sense: it is a curve of shortest
length among all piecewise-smooth curves with the same endpoints. From
the general theory of geodesics (see e.g. \cite{CE1975}), any such curve $\g$ is actually
{\em smooth}, and, when parametrized at constant speed, satisfies the
geodesic equation $\nabla_{\g'}\g'\ident 0$.  Thus \eqref{defdsr} is equivalent to
\be\label{defdsr2}
\dsr(X,Y) = \inf\left\{\ell(\g) \mid \g:[0,1]\to
M\ \mbox{is a geodesic with}\  \g(0)\in\E_X\ ,\ \g(1)\in
\E_Y\right\}.
\ee
Now with Definition~\ref{defmssr}(v), \eqref{defdsr2} becomes
\bearray\nonumber
\dsr(X,Y) &=& \inf\left\{\ell(\chi) \mid \chi:[0,1]\to
\sympp\ \mbox{is an SSR curve with}\right.\\
&&  \left.\phantom{\inf\left\{\ell(\chi) \mid \right.}
\  \chi(0)=X,\,\chi(1)=Y\right\}.
\label{defdsr3}
\eearray

\begin{remark}\label{whysmooth}{\rm
As noted in \cite{JSG2015scarot}, the ``scaling-rotation distance" $\dsr$ is
{\em not} a metric on
$\sympp$; it does not satisfy the triangle inequality.
(However, its restriction to the top stratum of
$\sympp$ {\em is} a metric; see \cite[Theorem 3.12]{JSG2015scarot}.)
}
\end{remark}

Computing $\dsr(X,Y)$ amounts to optimizing over the fibers of $X$ and $Y$.
Choosing $(U,D)\in\E_X, (V,\L)\in \E_Y$, it first appears from \eqref{fibform4} that this
requires optimizing over $(G_D^0\times\tsp^+)\times (G_\L^0\times \tsp^+)$,
thus a ``continuous" optimization over $G_D^0\times G_\L^0$ for each of the
$|\tsp^+|^2$ elements of $\tsp^+\times\tsp^+$.  However, there is
quite a bit of redundancy; clearly it suffices to do a continuous
optimization over each pair of connected components (an element of
$\Comp(\E_X)\times \Comp(\E_Y)$) and then a combinatorial optimization over
the finite set $\Comp(\E_X)\times \Comp(\E_Y)$.
When both $X$ and $Y$ are in the top stratum, the optimization (\ref{defdsr}) is purely combinatorial.
More generally, Proposition
  \ref{prop:biject}(i) implies that $|\Comp(\E_X)|
=|\tsp^+/\Gamma^0_{J_D}|=|\tsp^+|/|\Gamma^0_{\J_D}|$
and $|\Comp(\E_Y)|=|\tsp^+/\Gamma^0_{J_\L}| = |\tsp^+|/|\Gamma^0_{J_\L}|$,
so the product of these
two numbers is an upper bound on the number of continuous optimizations
needed. Even this bound is quite crude: using invariances of the metric on $M$,
it is not hard
to understand that the number of continuous optimizations
needed should not exceed $\min\{ |\Comp(\E_X)|, |\Comp(\E_Y)|\}$.   (In \cite{JSG2015scarot}, this idea is used in Theorems 4.2 and 4.3.) However,
Proposition \ref{distprop} below reduces this number further when neither $X$ nor $Y$ lies in the top or bottom stratum. To state the proposition, first recall
that given any group $G$ and subgroups $H_1,H_2$, an
$(H_1,H_2)$ {\em double-coset} is an equivalence class under the
equivalence relation $\sim$ on $G$ defined by declaring $g_1\sim g_2$
if there exist $h_1\in H_1, h_2\in H_2$ such that $g_2=h_1g_1h_2$. The
set of equivalence classes under this relation is denoted
$H_1\backslash G/H_2$.  By a {\em set of representatives} of
$H_1\backslash G/H_2$ we mean a subset of $G$ consisting of exactly
one element from each $(H_1,H_2)$ double-coset.

\begin{prop}[{\cite[Proposition 4.10]{GJS2016mathpaper}}]
\label{distprop}
Let $X,Y\in\sympp$ and let $(U,D)\in \E_X, (V,\L)\in \E_Y$.  Let $Z$
be any set of representatives of $\Gamma_{\J_D}^0\backslash \tsp^+/
\Gamma_{\J_\L}^0$.  Then the scaling-rotation distance from $X$ to $Y$ is
given by
\begin{align}
\dsr(X,Y)^2  = \min_{g\in Z}
  \left
    \{k\left(
          \wideparen{d}(g ; (U,D) , (V,\Lambda) )
        \right)^2
       +\|\log\left(D^{-1}(\pi_g\dotprod \L)\right)\|^2
      \right\}, \label{fibdist8a}
\end{align}
where
\begin{equation}
\wideparen{d}(g ; (U,D) , (V,\Lambda) )= \min_{R_U\in
                    G_D^0, R_V\in G_\L^0} \left\{ d_{\So}(UR_U, VR_VP_g^{-1}) \right\}. \label{fibdist8a-cts-opt}
\end{equation}

Every minimal smooth scaling-rotation curve from $X$ to $Y$
corresponds to some minimal pair whose first element lies in the
connected component $[(U,D)]$ of $\E_X$.
\end{prop}

%
%

To illustrate the reduction in the number of required continuous
optimizations \eqref{fibdist8a-cts-opt}
is reduced in the computation of $\dsr(X,Y)$, take for example $p = 3$ and  $X,Y \in \S_{[\J_{\rm mid}]},$ the middle stratum of $\symp(3)$, defined in Section~\ref{sec:exsr3}.
In this case we have $|\Comp(\E_X)|=|\Comp(\E_Y)|=6$, but, as we shall
see in Section~\ref{sec:tz}, the set $Z$ in \eqref{fibdist8a} has cardinality 3. Thus Proposition \ref{distprop} reduces the number of continuous optimizations needed down to 3.

\subsection{Existence and uniqueness of MSSR curves}\label{sect:uniqq}

From Proposition \ref{prop:biject}, every fiber of $F$ is compact, so
the infimum in \eqref{defdsr} is always achieved. Hence for all
$X,Y\in \sympp$, there always exists an
$(\E_X,\E_Y)$-minimal geodesic, a minimal
pair in $\E_X\times \E_Y$, and an MSSR curve from $X$ to $Y$.

Such an MSSR curve may not be unique. In \cite{JSG2015scarot}, a sufficient condition for uniqueness is given, and an example for $p = 2$ is provided. With statistical analysis in mind, it is natural to ask:  For which $X$ and $Y$ is there a unique MSSR curve from $X$ to $Y$? We address this question more generally by characterizing $\M(X,Y)$ for all $X,Y \in \sympp$.
In Sections~\ref{sec:appendix-MSSR2}, \ref{sec:mssr3} and \ref{sec:appendix-MSSR3}, we do this explicitly for low-dimensional cases: $p =2$ and $3$.
As preparation for this work, we briefly discuss here how non-uniqueness can occur and introduce a tool used to characterize $\M(X,Y)$  in low dimensions. A general treatment of this topic can be found in \cite{GJS2016mathpaper}.

Different $(\E_X,\E_Y)$-minimal
geodesics may or may not project to the same MSSR curve.
For given $X,Y$, for uniqueness of an MSSR curve from $X$ to $Y$ to fail, there must be distinct $(\E_X,\E_Y)$-minimal geodesics $\g_i:[0,1]\to M$, whose endpoints are minimal pairs $((U_i,D_i),(V_i,\L_i)) \in  \E_X \times \E_Y ,$ $i=1,2$, such that $F\circ\g_1\neq F\circ \g_2$. There are two possible ways in which this failure can occur: There exist such $\g_i$ whose endpoints are {\em distinct} minimal pairs $((U_i,D_i),(V_i,\L_i))$ (``Type I non-uniqueness''),  or {\em the same} minimal pair $((U,D),(V,\L))$ (``Type II non-uniqueness'').
%
%
%

Since for any $D,\L\in \D^+(p)$ the minimal geodesic from $D$ to $\L$ is unique, Type II non-uniqueness with minimal pair $((U,D),(V,\L))$ is equivalent to the existence of two or more minimal geodesics from $U$ to $V$,  which is equivalent to $U^{-1}V$ being an involution.
For $p =2,3$ it is shown in \cite{JSG2015scarot} that Type II non-uniqueness never occurs. This is because that, for $p \le 3$, for any pair $U,V \in \Sop$ such that $U^{-1}V$ is an involution, there exists a $\boldsig\in \I_p^+$ such that  $\dso(UI_\mbs, V)<\dso(U,V)$. In \cite{GJS2016mathpaper}, it is further shown that for small enough values of $p$, Type II non-uniqueness never occurs; for large enough $p$, it always occurs. In particular, for $p\leq 4$,  for all $X,Y\in \sympp$ for which $\M(X,Y)>1$, the non-uniqueness is purely of Type I.

Our strategy for understanding $\M(X,Y)$ for $p = 2, 3$ and all $X,Y \in \sympp$ is to list all MSSR curves from $X$ to $Y$.
Proposition \ref{distprop} assures us that,  for
any $(U,D)\in \E_X$, every MSSR curve from $X$ to $Y$ corresponds to
some minimal pair whose first element lies in the connected component
$[(U,D)]$ of $\E_X$.
We need a way to tell whether MSSR curves corresponding
to two minimal pairs with first point in
$[(U,D)]$ are the same.
 The following proposition, a special case of Proposition 4.19 of \cite{GJS2016mathpaper}, provides such a tool. We apply this result to the $p=3$ case in Section~\ref{sec:mssr3}.

\begin{prop} \label{prop:MSRequal}
Let $p \le 4$, $X,Y \in\sympp, X\neq Y$. For $i=1,2$ assume that
$\chi_i=F\circ\g_i$ is an MSSR curve from
$X$ to $Y$ corresponding to the minimal pair
$((UR_{U,i},D),(VR_{V,i}\,P_{g_i}^{-1}, \L_i)),$ where $R_{U,i}\in
G_D^0, R_{V,i}\in G_\L^0$, $g_i\in \tsp^+$, $\L_i=\pi_{g_i}\dotprod
\L$, and $\g_i:[0,1]\to M$ is a geodesic. Then $\chi_1=\chi_2$ if and only if
the following two conditions hold.

\begin{itemize}

\item[(i)] $
R_{V,2}P_{g_2}^{-1}R_{U,2}^{-1} =R_{V,1}P_{g_1}^{-1}R_{U,1}^{-1}\,;$
%
%

\item[(ii)] There exist $g\in\tsp^+, R\in G_{D,\L_1}^0$ such that
\begin{eqnarray}
\label{d2pdd1}
D &=& \pi_g\dotprod D, \\
\L_2 &=& \pi_g\dotprod \L_1, \\
\mbox{\rm and} \ \ \ R_{U,1}^{-1}R_{U,2} &=& R P_g^{-1}\,.
\label{rppp}
\end{eqnarray}

\end{itemize}

\end{prop}

\setcounter{equation}{0}

\section{Scaling--Rotation distance and MSSR curves in $\symp(2)$}\label{sec:appendix-MSSR2}
The space $\symp(2)$ has only two strata: $\Sbot := \Sc_{[\J_{\rm bot}]} = \{\lambda I : \lambda > 0 \}$ and $\Stop := \Sc_{[\J_{\rm top}]} = \symp(2) \backslash \Sbot$. To characterize all unique and non-unique cases of minimal smooth scaling-rotation (MSSR) curves in $\symp(2)$, it suffices to consider two possibilities for the strata in which $X,Y$ lie:
\begin{enumerate}
\item[(i)]  $X,Y \in \Stop$.
\item[(ii)] $X \in \Sbot$ ($Y \in \Stop \cup \Sbot = \symp(2)$).
\end{enumerate}

For any $X,Y \in \symp(2)$, with $a \ge b$, $c \ge d$ and $0 \le \theta < \pi$, one can write
\begin{equation} \label{eq:XandY}
X =  U \begin{pmatrix}
          e^a & 0 \\
          0    & e^b
\end{pmatrix} U^T , \
Y = U R_\theta
\begin{pmatrix}
          e^c & 0 \\
          0    & e^d
\end{pmatrix}
R_\theta^T U^T
\end{equation}
where $U \in \SO(2)$ and $R_\theta = \exp (A_\theta),$ $A_\theta = \begin{pmatrix}
          0 & -\theta \\
          \theta    & 0
\end{pmatrix}$.
Denote the apparent eigen-decompositions of $X$ and $Y$ appearing in (\ref{eq:XandY}) by $(U,D)$ and $(V,\Lambda)$. Then the scaling--rotation curve $\chi$ with parameters
$(U,D,A,L)$, where $A = \log(VU^T)$ and $L = \exp(D^{-1}\Lambda)$, is
\begin{equation}\label{eq:scarotcurve}
\chi(t) = (U R_{\theta t})
  \exp   \begin{pmatrix}
               {(1-t)a + t c } & 0 \\
               0    & {(1-t)b + t d }
     \end{pmatrix}
     (UR_{\theta t} )^T,
\end{equation}
satisfying $\chi(0) = X$, $\chi(1) = Y$.

\textbf{Case (i)} ($a > b$, $c > d$). Let $(V_i,\Lambda_i$), $i = 1,\ldots,4 $ be the four eigen-decompositions of $Y$. Specifically, these four eigen-decompositions are
\begin{align*}
 (V_1,\Lambda_1) &= (UR_\theta, \diagg(e^c, e^d)), \\
  (V_2,\Lambda_2) &= (UR_{\theta + \pi}, \diagg(e^c, e^d)), \\
   (V_3,\Lambda_3) &= (UR_{\theta + \pi/2}, \diagg(e^d, e^c)), \\
    (V_4,\Lambda_4) &= (UR_{\theta - \pi/2}, \diagg(e^d, e^c)).
\end{align*}
Then $d_{\Sc\Rc}(X,Y) = \min_{i = 1,2,3,4} d_i$, where $d_i =  d((U,D),(V_i,\Lambda_i))$.
For $0 \le \theta \le \pi/2$,
$$ d_1^2 = k\theta^2 + (a-c)^2 + (b-d)^2 \le d_2^2, $$
and  equality holds if and only if $\theta = \pi/2$.
On the other hand, if $\pi/2 < \theta < \pi$,
$$ d_2^2 = k(\pi - \theta)^2 + (a-c)^2 + (b-d)^2 < d_1^2, $$
For  $0 \le \theta < \pi$,
$$ d_3^2 = k(\theta - \pi/2)^2 + (a-d)^2 + (b-c)^2 \le d_4^2, $$
and  equality holds if and only if $\theta = 0 $.
 Furthermore, we have
\begin{align*}
 d_1^2 < d_3^2 \quad &\Longleftrightarrow \quad
   \theta < \frac{\pi}{4} + \frac{2(a-b)(c-d)}{k\pi}, \\
 d_2^2 < d_3^2 \quad &\Longleftrightarrow \quad
   \theta > \frac{3\pi}{4} - \frac{2(a-b)(c-d)}{k\pi}.
\end{align*}
These inequalities will be used later in the characterization of all MSSR curves for Case (i).

Case (i) has seven subcases: three in which there is a
unique MSSR curve, three in which there are  non-unique MSSR curves with multiplicity 2, and one in which there are non-unique MSSR curves with multiplicity 3. We denote these subcases  ``$d_i$'', ``$d_i = d_j$'', and ``$d_1 = d_2 = d_3$'' respectively.
In the subcase denoted by ``$d_i$'', the MSSR curve from $X$ to $Y$ is unique, has length $d_{\Sc\Rc}(X,Y) = d_i$, and corresponds to the minimal pair $((U,D),(V_i,\Lambda_i))$ using (\ref{eq:scarotcurve}).
In the subcase denoted ``$d_i = d_j$'', there are  exactly two MSSR curves from $X$ to $Y$, of length  $d_{\Sc\Rc}(X,Y) = d_i =d_j$, and corresponding to
 the minimal pairs $((U,D),(V_i,\Lambda_i))$ and $((U,D),(V_j,\Lambda_j))$. The notation for the last subcase with three MSSR curves is similarly understood.

\begin{figure}[t]
  \centering
\includegraphics[width=0.7\textwidth]{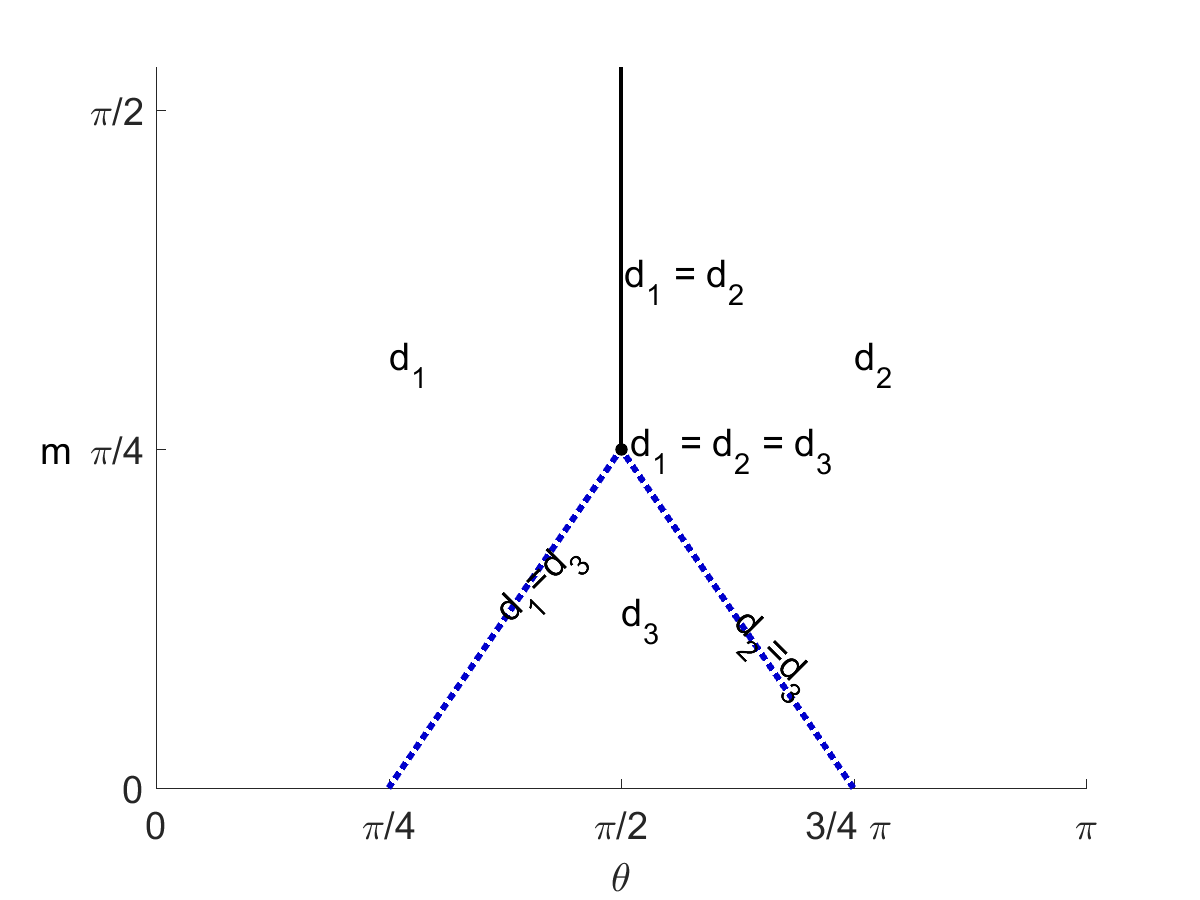}
  \caption{Unique and non-unique MSSR curves in $\symp(2)$: Schematic illustration for the seven subcases of Case (i) in which both $X$ and $Y$ are in $\Stop$.   \label{fig:1HDGeom}}
\end{figure}

 The seven subcases are distinguished by the relationship of the quantity $m := \frac{2(a-b)(c-d)}{k\pi}$ and the angle $\theta$.
If
$ m  >  \min(\theta,\pi-\theta) - \pi/4$, then
$$
d_{\Sc\Rc}(X,Y) = \left\{
             \begin{array}{ll}
             d_1, & \theta < \pi/2,  \\
            d_1 = d_2 ,& \theta = \pi/2,  \\
             d_2, & \theta > \pi/2.  \\
             \end{array}
             \right.
$$
If
$ m  =  \min(\theta,\pi-\theta) - \pi/4$, then
$$
d_{\Sc\Rc}(X,Y) = \left\{
             \begin{array}{ll}
             d_1=d_3, & \theta < \pi/2,  \\
            d_1 = d_2 = d_3 ,& \theta = \pi/2,  \\
             d_2 = d_3, & \theta > \pi/2.  \\
             \end{array}
             \right.
$$
Finally, if
$ m  < \min(\theta,\pi-\theta) - \pi/4$, then
$
d_{\Sc\Rc}(X,Y) = d_3.
$
The conditions leading to these seven subcases are graphically summarized in Fig.~\ref{fig:1HDGeom}.

The MSSR curves corresponding to subcases ``$d_1$'' and ``$d_2$'' can be understood as the rotation and scaling of the ellipse corresponding to $X$ to the ellipse corresponding to $Y$, where the rotation is either counterclockwise (case ``$d_1$''), or clockwise (case ``$d_2$''). (There is no rotation in subcase ``$d_1$'' if $\theta = 0$.) In these two subcases, the whole MSSR curve never leaves the distinct-eigenvalue subset of $\symp(2)$ (the top stratum). On the other hand, the MSSR curves corresponding to ``$d_3$'' always pass through the equal-eigenvalue subset of $\symp(2)$ (the bottom stratum). The direction of rotation for ``$d_3$'' depends on $\theta$: counterclockwise if $\theta < \pi/2$, clockwise if $\theta> \pi/2$. (There is no rotation if $\theta = \pi/2$.)

\begin{figure}[tb]
  \centering
\includegraphics[width=0.70\textwidth]{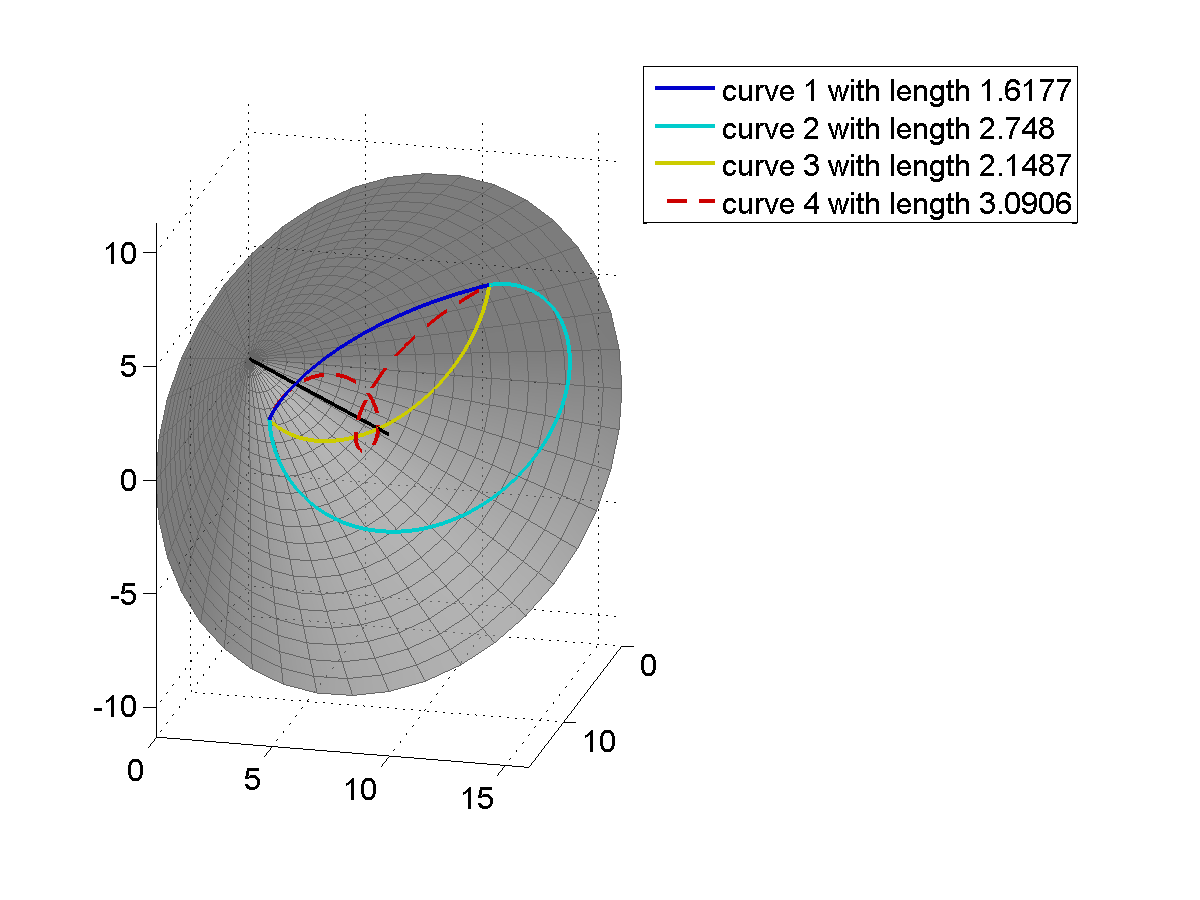}
\hspace{-1.5in}
\includegraphics[width=0.50\textwidth]{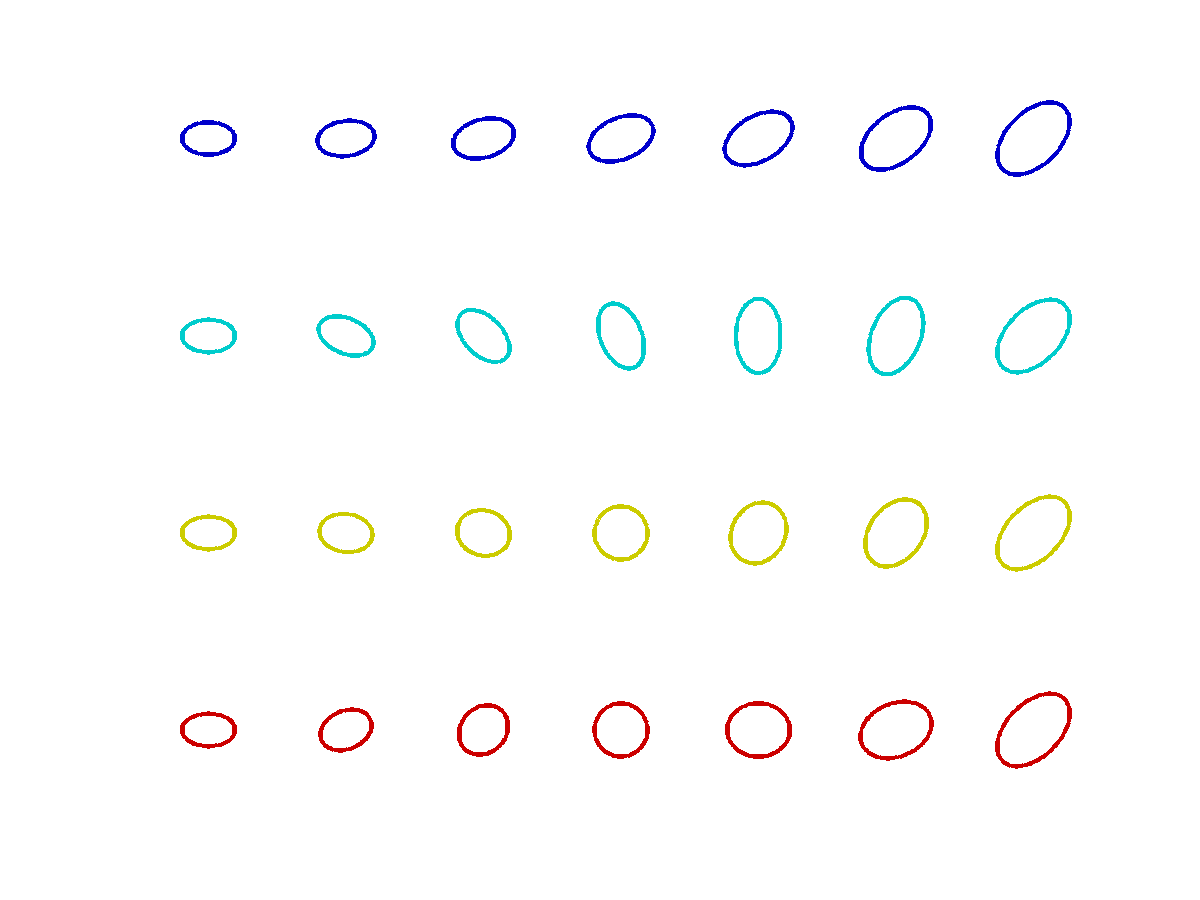}
\vspace{-.2in}
  \caption{An example for subcase ``$d_1$''. \label{fig:d1} For each $i = 1,2,3,4$, ``curve $i$'' represents the scaling-rotation curve corresponding to the pair $((U,D),(V_i,\Lambda_i))$, and whose length is $d_i$.
$\symp(2)$ is an open cone in the three-dimensional space $\sym(2)$. The black line is the axis of the cone.
Each curve, labeled by the same color in the left and right panels, is illustrated as the space curve $\chi(t)$, contained in this cone $\symp(2)$ (left), or as the sequence of corresponding ellipses (right). In this and all other examples,  the scaling factor $k$ is set equal to  $1$.}
\end{figure}

\begin{figure}[tb]
  \centering
\includegraphics[width=0.70\textwidth]{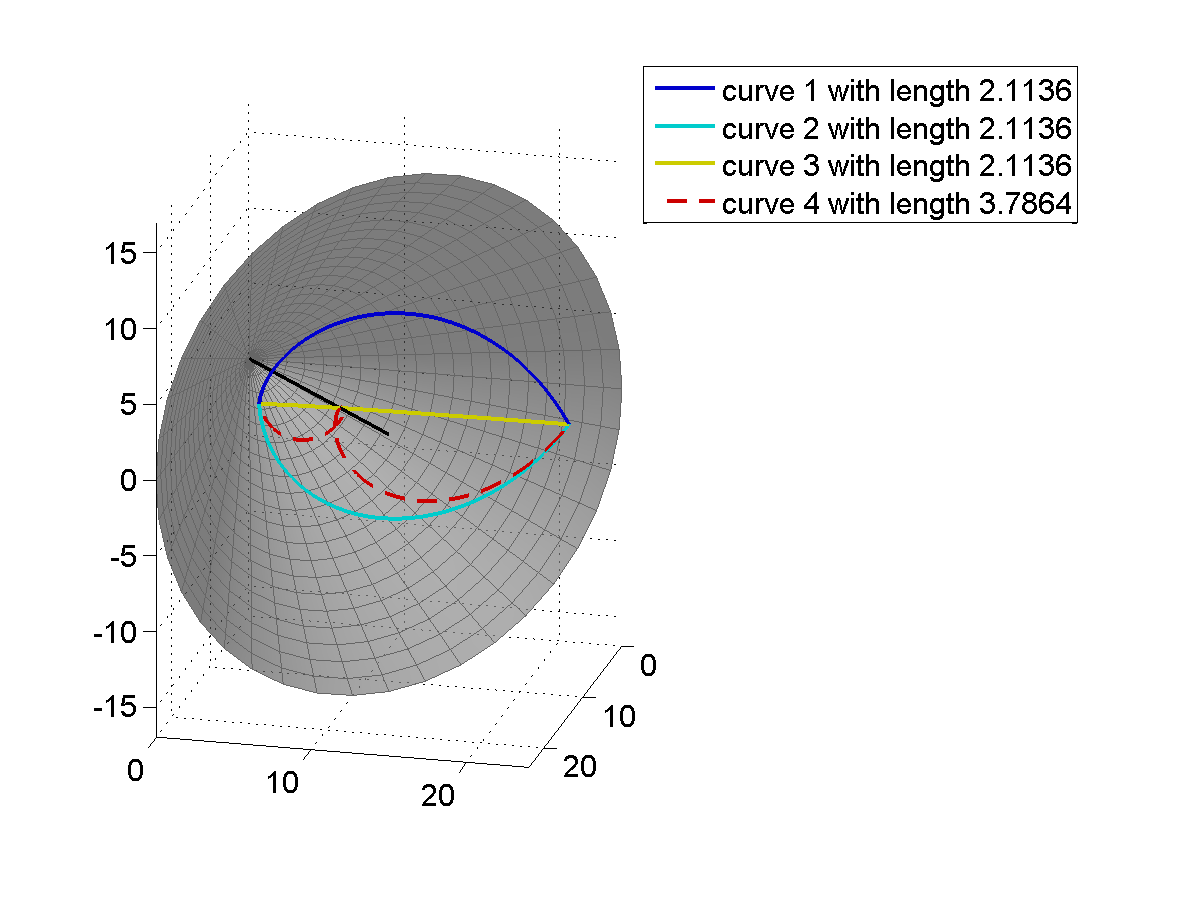}
\hspace{-1.5in}
\includegraphics[width=0.50\textwidth]{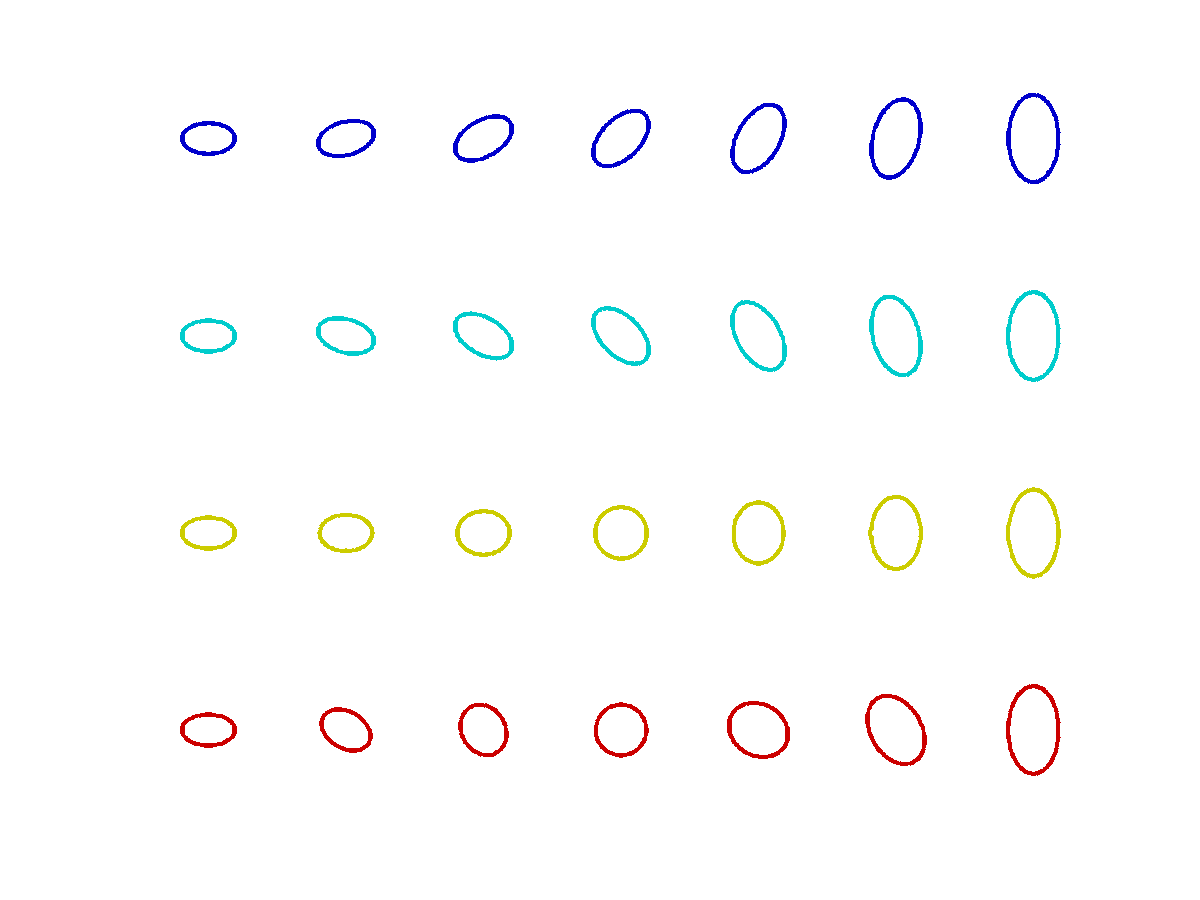}
 \vspace{-.2in} \caption{An example for subcase ``$d_1 = d_2 =  d_3$''. \label{fig:d1=d2=d3}}
\end{figure}

A few of these seven subcases are illustrated by  representative examples in Figs. \ref{fig:d1} and \ref{fig:d1=d2=d3}. If ``$d_1$'' and ``$d_2$'' are thought of as the same ``type'', then there are five different types of (non-)uniqueness behavior of MSSR curves in Case (i), as follows:
\begin{enumerate}
\item Unique MSSR curve (completely contained in the distinct-eigenvalue subset), if $ m  >  \min(\theta,\pi-\theta) - \pi/4$ and $\theta \neq \pi/2$. Subcases ``$d_1$'' (Fig.~\ref{fig:d1}) and ``$d_2$'' are of this type.
\item Unique MSSR curve (leaving   the distinct-eigenvalue subset and passing through the bottom stratum), if $ m  < \min(\theta,\pi-\theta) - \pi/4$. Subcase ``$d_3$''  is of this type.
\item Two MSSR curves with rotation angle $\pi/2$ (still completely contained in the distinct-eigenvalue subset), if  $ m  >  \min(\theta,\pi-\theta) - \pi/4$ and $\theta  =  \pi/2$. Case ``$d_1 = d_2$'' is of this type.
\item Two MSSR curves (one in the distinct eigenvalue subset, the other passing through the bottom stratum), if $ m  =  \min(\theta,\pi-\theta) - \pi/4$ and $\theta \neq \pi/2$. Subcases ``$d_1 = d_3$'' and ``$d_2 = d_3$''  are of this type.
\item Three MSSR curves (two with rotation-angle $\pi/2$, completely contained in the distinct-eigenvalue subset, and the other involving no rotation but passing through the bottom stratum), if $ m  =  \min(\theta,\pi-\theta) - \pi/4$ and $\theta  =  \pi/2$. Subcase ``$d_1=d_2=d_3$'' (Fig.~\ref{fig:d1=d2=d3}) is of this type.
\end{enumerate}

For each given $X$ and $Y$, if one takes $k$ small enough that $m > \pi/4$, then MSSR curves from $X$ to $Y$ are always of type ``$d_1$'' or ``$d_2$''. In other words, if $k$ is small enough (for fixed $X,Y$), the MSSR curve(s) are completely contained in the distinct-eigenvalue subset.

%
%
%
%
%
%
%
%
%

\textbf{Case (ii)} ($a=b$ and $c \ge d$). By Theorem 4.1 of  \cite{JSG2015scarot}, the MSSR curve from $X$ to $Y$ is unique, and is
$$\chi(t) = (U R_{\theta})
  \exp   \begin{pmatrix}
               {(1-t)a + t c } & 0 \\
               0    & {(1-t)a + t d }
     \end{pmatrix}
     (UR_{\theta} )^T,$$ and
$d_{\Sc\Rc}(X,Y) = \sqrt{ (a-c)^2 + (a-d)^2}$.
If in addition $c=d$, then $\chi(t) = e^{(1-t)a + t c } I$.

\setcounter{equation}{0}
\section{Scaling-rotation distances on $\Sym^+(3)$}\label{sec:dsr3}

For the case $p=3$, we will obtain explicit formulas for all MSSR curves and
scaling-rotation distances by using the quaternionic parametrization of $\So(3)$.
In Section \ref{sec:mssr3}, we use this to give explicit descriptions
of the set $\M(X,Y)$ of MSSR curves between two points $X,Y\in\Sym^+(3)$ in all ``nontrivial"
cases (as defined later in this section).

   \subsection{Characterization of SR distance using quaternions}\label{sec:5.1.quaternion}
    \subsubsection{Relation of quaternions to ${\rm SO} (3)$}\label{sec:5.1.1.relhso3}

The space $\H$ of quaternions, with its usual real basis
$\{1,i,j,k\}$ identified with the standard basis of $\bfr^4$, and with
$\{i,j,k\}$ identified with the standard basis of $\bfr^3$, provides a
convenient parametrization of $\So(3)$. Specifically, writing
$S^3_\H=$
$\{x_0+x_1i+x_2j+x_3k\ :\ x_0^2+x_1^2+x_2^2+x_3^2=1\}$, there is a
natural two-to-one Lie-group homomorphism $\phi:S^3_\H\to\So(3)$,
defined as follows.  Using the basis $\{i,j,k\}$ to identify $\bfr^3$
with $\imh$, the space of purely imaginary quaternions, for $q\in
S^3_\H$ and $x\in\imh$ we set $\phi(q)(x)=qx\bar{q}$, which lies in
$\imh$. For $q_1,q_2\in S^3_\H$ we have $\phi(q_2)=\phi(q_1)$ if and
only if $q_2=\pm q_1.$ Thus, for any $U\in \So(3)$, if $q_U\in
\phi^{-1}(U)$ then
\be\label{invimg} \phi^{-1}(U)= \{\pm q_U\}.  \ee

Let $S^2_\imh = \{\tilde{a}\in\imh \ :\ \|\tilde{a}\| =1\}$. For
$\tilde{a}\in S^2_\imh$ and $\th\in [0,\pi]$ let $R_{\th,\tilde{a}}$
denotes counterclockwise rotation by angle $\th$ about the axis
$\tilde{a}$ (``counterclockwise'' as determined by $\tilde{a}$ using the
right-hand rule).  Let
\ben \So(3)_{<\pi}:= \{R_{\th,\tilde{a}} \ :\ \tilde{a}\in S^2_\imh,
\th\in [0,\pi)\}, \een

\noi the set of non-involutions in $\So(3)$.  The map
$s:\So(3)_{<\pi}\to S^3_\H$ defined by
\be\label{halfangle} s(R_{\th,\tilde{a}})= \cos\frac{\th}{2} +
(\sin\frac{\th}{2})\tilde{a}  \ee
 is a smooth right-inverse to $\phi$ on $\So(3)_{<\pi}$
(i.e., $\phi\circ s$ is the identity map on this domain), but $s$ is
not a homomorphism and cannot be extended continuously to all of
$\So(3)$.

Distances between elements of $\So(3)$ are related very simply to
geodesic distances in $S^3_\H$ with respect to the standard Riemannian
metric on $S^3=S^3_\H$.  For $U,V\in\So(3)$,
\bestar d_{\So}(U,V)&=&2\min\{d_{S^3}(r_U, r_V)\ :\ r_U\in
\phi^{-1}(U), r_V\in
\phi^{-1}(V)\}\\ &=&2\cos^{-1}|\Re(\Bar{r_U}\ r_V)| \ \mbox{\rm where
  $r_U, r_V$ are as in previous line.}  \eestar

\noi Alternatively, $d_{\So}(U,V) = d_{\So}(\mbbi, U^{-1}V)
=2\cos^{-1}|\Re(q_{U^{-1}V})|$, where $q_{U^{-1}V}$ is either of the two
elements of $\phi^{-1}(U^{-1}V)$.  Thus
\be\label{master-b} d_{\So}(U,V) = 2\cos^{-1}|\Re(\Bar{q_U}\ q_V)|
= 2\cos^{-1}|\Re(q_{U^{-1}V})| \ee

\noi where each of $q_U, q_V, q_{U^{-1}V}$ is either of the two elements
in $S^3$ mapped by $\phi$ to $U, V,$ and $U^{-1}V$ respectively.

\subsubsection{Quaternionic pre-images of
signed-permutation   matrices}\label{sect:5.1.2}

For every subgroup $H\subset \So(3)$, let $\Hat{H}$ denote the
pre-image $\phi^{-1}(H)\subset S^3_\H$.  Writing
$\Gamma:=\ts_3^+\subset\So(3)$, the 24-element group of even signed-permutation matrices, the 48 elements of $\Hat{\Gamma}$ are the
following:
\bearray
\label{line1}
&\pm 1, \pm i, \pm j, \pm k,&\\
\label{line2}
&\frac{1}{\sqrt{2}}(\pm e_m\pm e_n), \mbox{where}\ e_0=1, e_2=i,
e_3=j, e_4=k,\ \mbox{and}\ m<n,& \\
\label{line3}
&\frac{1}{2}(\pm 1\pm i\pm j\pm k)&
\eearray

\noi (all sign-combinations allowed in all sums).  Under $\phi$, the
eight elements \eqref{line1} are mapped to the four even sign-change
matrices, the 24 elements \eqref{line2} are mapped to the 12
positive-determinant ``signed transposition matrices''
(permutation-matrices corresponding to transpositions, with an odd
number of 1's replaced by $-1$'s), and the 16 elements \eqref{line3}
are mapped to the 8 positive-determinant ``signed cyclic-permutation
matrices'' (permutation matrices corresponding to cyclic permutations,
with an even number of 1's replaced by $-1$'s).

   \subsubsection{Parameters corresponding to different strata}\label{sect:5.1.3}


For any subgroups $H_1, H_2$ of $\So(3)$, the map $\phi$ induces a
bijection $\Hat{H_1}\backslash \Hat{\Gamma}/ \Hat{H_2} \to
H_1\backslash \Gamma / H_2$. Thus if $\widehat{Z}$ is a set of
representatives of $\Hat{\Gamma_D^0}\backslash \Hat{\Gamma}
/\Hat{\Gamma_\L^0}$ in $\Hat{\Gamma}$, then $\phi(\widehat{Z})$ is a
set $Z$ of representatives of
$\Gamma_D^0\backslash \Gamma / \Gamma_\Lambda^0$ in $\Gamma$.

Now let $X,Y, (U,D),$ and $(V,\L)$ be as in Proposition
\ref{distprop}, and let $\widehat{Z}$ be a set of representatives of
$\Hat{\Gamma_D^0}\backslash \Hat{\Gamma} /\Hat{\Gamma_\L^0}$.
For $\z\in \Hat{\Gamma}$, define $\pi_\z = \pi_{\phi(\z)} \in S_3$
(see Notation \ref{defpig}).

From equation \eqref{master-b} and the fact
that $\phi$ is a homomorphism, it follows that in the setting of
equation \eqref{fibdist8a} (with $p=3$), for all $r_U\in
\phi^{-1}(R_U), r_V\in \phi^{-1}(R_V), g\in \ts_3^+,$ and $\z_g\in
\phi^{-1}(g)$,
\be\label{fibdist2}
d_{\So}(UR_U, VR_VP_g^{-1})^2
=4\left(\cos^{-1}\left|\Re\left(r_V\ \Bar{\z_g}\ \Bar{r_U}\ q_{U^{-1}V}
\right)\right|\right)^2.
\ee

\noi From Proposition \ref{distprop}, we therefore have
%
\begin{align} \label{fibdist8-new}
\dsr(X,Y)^2    &  = \min_{\z\in \widehat{Z}}
  \left
    \{k \,
          \hat{d}(\zeta)^2
       +\|\log(\L_{\pi_\z}D^{-1})\|^2
      \right\},
\end{align}
where
\begin{align} \label{fibdist8-new2}
\hat{d}(\z) :=\breve{d}(\zeta ; \Hat{G_D^0}, \Hat{G_\L^0}, q_{U^{-1}V}) := \min_{r_U \in \Hat{G_D^0}, r_V\in \Hat{G_\L^0} } \left\{ \cos^{-1} | \Re\left(r_V\ \bar{\z}\ \Bar{r_U}\ q_{U^{-1}V}\right) | \right\}.
\end{align}

As equations \eqref{fibdist8-new}-\eqref{fibdist8-new2}  suggest, computing $\dsr(X,Y)$ is a
minimization problem that breaks into two parts, one over the discrete
parameter-set $\widehat{Z}$ and the other over the (potentially)
``continuous'' parameter-set $\Hat{G_D^0}\times \Hat{G_\L^0}$. Both
parameter-sets depend on $X$ and $Y$. 


If $X$ or $Y$ lies in the bottom stratum of $\symp(3)$ (i.e., has only one distinct
eigenvalue),
then the set $\widehat{Z}$ has only one element $\z$, which we can take
to be $1\in \H$, and at least one of the groups
$\Hat{G_D^0},\Hat{G_\L^0}$ is all of $S^3_\H$. The set
$\left\{r_V\ \bar{\z}\ \Bar{r_U}: (r_U,r_V)\in \Hat{G_D^0}\times
\Hat{G_\L^0}\right\}$ is simply $S^3_\H$, the inner minimum $\hat{d}(\z)$ in
\eqref{fibdist8-new} is 0, and we immediately obtain $\dsr(X,Y)=\|\log (\L
D^{-1})\|$. We do not need to use quaternions to obtain this
result; it follows just as quickly from \eqref{fibdist8a}.

At the other extreme, if $X$ and $Y$ both lie in the top stratum of
$\sympp$ (i.e. both have three distinct eigenvalues) then
$\Hat{G_D^0}=\Hat{G_\L^0}=\{\pm 1\}$,  and the inner minimum is trivial to
compute ($\hat{d}(\zeta) = \cos^{-1}| \Re (\bar{\z} q_{U^{-1}V}) |$), 
so we are reduced immediately to a single minimization
over $\widehat{Z}$.  As in the previous case, we do not need the
quaternionic reframing of the distance formula at all: already in
\eqref{fibdist8a} we have $G_D^0=G_\L^0=\{\mbbi\}$, so the distance can
be found simply by minimizing over the discrete variable $g\in Z=\tilde{S}_3^+$.
 We need only have a computer calculate $d_M((U,D),(VP_g^{-1}, g\dotprod \L)$ for each of the
24 $g$'s and return the corresponding minimal pairs and MSSR curves.
For combinatorial reasons, a complete algebraic classification of the
pairs $(X,Y)$ (with $X,Y$ both in the top stratum of $\Sym^+(3)$) for which $\M(X,Y)$ has a given cardinality would be very complex, and   we do not attempt this.  

For the above reasons, for the remainder of this section we focus
on the cases in which $X$ and $Y$ do not both lie in the
top stratum of $\Sym^+(3)$, and neither lies in the bottom stratum.  Thus we restrict attention
to the cases in which one of the
matrices $X,Y$ has exactly two distinct eigenvalues, and the other has
either two or three.
We refer to these cases as the ``nontrivial''
cases (because the set of distances between elements of $\E_X$ and elements of
$\E_Y$ is not a finite set). To analyze them we introduce the following notation:
 \begin{align*}
\esstop &:=  \S_{[\J_{\rm top}]} = \S_{1+1+1} \subset \Sym^+(3),\\
\smid &:=  \S_{[\J_{\rm mid}]}  = \S_{2+1}  \subset  \Sym^+(3),\\
\cald_{\J_1} &:=\cald_{\{\{1\}, \{2,3\}\}}
 = \{{\rm diag}(d_1, d_2, d_2)
 :\ d_1, d_2>0, \ d_1\neq d_2\} \subset \D^+(3).
\end{align*}


   \subsection{Scaling-rotation distances for $\Sym^+(3)$ in the nontrivial cases}\label{sect:5.2}



From now through  Section \ref{sec:mssr3} we
assume that $X\in \smid$ and that either $Y\in \esstop$ or $Y\in \smid$.
Then $X$ has an eigen-decomposition $(U,D)$ with $D\in \cald_{\J_1}$, and if $Y\in \smid$
then $Y$ has a eigen-decomposition $(V,\L)$ with $\L\in \cald_{\J_1}$. We will always
assume that our pairs $(U,D), (V,\L)$ have been chosen this way. Then
we have
\ben
G_D^0= G_{\cald_{\J_1}}^0:=\left\{\left[ \begin{array} {cl} 1 & {\bf 0}^T \\ {\bf 0}
    & R\end{array}\right] \ :\ R\in \So(2)\right\},
\ \ G_\L^0= \left\{\begin{array}{ll} \{\mbbi\} & \mbox{if $Y\in
  \esstop$,}
\\
G_{\cald_{\J_1}}^0& \mbox{if $Y\in\smid$.}
\end{array}\right.
\een

\noi It is not hard to check that
\bearray\label{cinh}
\Hat{G_D^0}&=&S^1_\bfc := \{\exp(ti)\ :\  t\in\bfr\}
\subset\bfc\subset\H,   \\
\label{pm1}
\Hat{G_\L^0}&=& \left\{ \begin{array}{ll} \{\pm 1\} & \mbox{if
    $Y\in\esstop$,}\\  S^1_\bfc& \mbox{if $Y\in\smid$.}
\end{array}\right.
\eearray

\noi The inner minimum in \eqref{fibdist8-new} is then
\begin{equation} \label{fibdist9-new}
\hat{d}(\zeta) =  \left\{
             \begin{array}{ll}
                       \min_{r_U\in S^1_\bfc} \left\{\cos^{-1}\left|\Re\left(\bar{\z}\ \Bar{r_U}\ q_{U^{-1}V}\right)\right|\right\}  & \mbox{if $Y\in\esstop$,}\\
                       \min_{r_U,r_V\in S^1_\bfc} \left\{\cos^{-1}\left|\Re\left(r_V\ \bar{\z}\ \Bar{r_U}\ q_{U^{-1}V}\right)\right|\right\} & \mbox{if $Y\in\smid$.}
             \end{array}
                     \right.
\end{equation}


\noi Obviously, minimizing the arc-cosines above is equivalent to
maximizing
\bearray
\left|\Re\left(\bar{\z}\ \Bar{r_U}\ q_{U^{-1}V}\right)\right|
&\mbox{if}& Y\in \esstop\,, \phantom{\int_a^b}
\label{fibdist7}
\\ \left|\Re\left(r_V\ \bar{\z}\ \Bar{r_U}\ q_{U^{-1}V}\right)\right|
&\mbox{if}& Y\in \smid\,,
\label{fibdist3}
\eearray

\noi where the maximum is taken over $(r_U,r_V)\in S^1_\bfc\times
S^1_\bfc$ in \eqref{fibdist3}, and over just $r_U\in S^1_\bfc$ in
\eqref{fibdist7}.

\subsubsection{The discrete parameter-sets $\widehat{Z}$ in the
nontrivial cases}\label{sec:tz}

To compute the outer minimum (over $\z$) in  \eqref{fibdist8-new}, we will
need to select sets $\widehat{Z}$ of representatives of
$\Hat{\Gamma_D^0}\backslash \Hat{\Gamma} /\Hat{\Gamma_\L^0}$ in two
cases: (i) $\J_D=\J_1:= \{1,\{2,3\}\}, \J_\L=\jtop:=\{\{1\}, \{2\},
\{3\}\}$, and (ii) $\J_D=\J_1=\J_\L$.  Let us write
$\Gamma_1:=\Gamma_{\J_1}^0$. For case (i), since $\Hat{\Gamma_{\jtop}^0}
=\{\pm 1\}$, which commutes with every element of $\Hat{\Gamma}$ and
is a subgroup of $\Hat{\Gamma_\J^0}$ for every $\J$, the double-coset
space $\Hat{\Gamma_1}\backslash \Hat{\Gamma} /\Hat{\Gamma_{\jtop}^0}$
is simply the set $\Hat{\Gamma_1}\backslash\Hat{\Gamma}$ of right
$\Hat{\Gamma_1}$-cosets in $\Hat{\Gamma}$.
Observe that $\Hat{\Gamma_1}
  =\Hat{G_{\J_1}^0}\ \intersect \ \Hat{\Gamma} =
  S^1_\bfc\ \intersect\Hat{\Gamma} =\{\pm 1, \pm i, \frac{\pm 1\pm
    i}{\sqrt{2}}\}=\{e^{2n\pi i/8}: n=0,1,\dots, 7\}$. Thus the
  cardinality of $\Hat{\Gamma_1}\backslash\Hat{\Gamma}$ is 48/8 =
  6.  One can check that the following set $\widehat{Z}_{1,*}$ contains
  a representative of each of the six right $\Hat{\Gamma_1}$-cosets:

\be\label{cosetreps} \widehat{Z}_{1,*}:= \left\{1, \ j, \frac{1\pm j}{\sqrt{2}},
\ \frac{1\pm k}{\sqrt{2}}\right\}.  \ee

\noi For case (ii), the double-coset space can be viewed as the set of
orbits under the action of $\Hat{\Gamma_1}$ on
$\Hat{\Gamma_1}\backslash\Hat{\Gamma}$ (the coset-space in case
(i)) by right-multiplication. Thus a set of representatives can be
found by imposing the double-coset equivalence relation on the set
$\widehat{Z}_{1,*}$ above. The four elements $\frac{1\pm j}{\sqrt{2}},
\ \frac{1\pm k}{\sqrt{2}}$ all lie in the same $\Hat{\Gamma_1}$
  double-coset, since
$$-i\frac{1+j}{\sqrt{2}} i =
  \frac{1-j}{\sqrt{2}},\ \ -i\frac{1+k}{\sqrt{2}} i =
  \frac{1-k}{\sqrt{2}},\ {\rm and}\ \frac{1+i}{\sqrt{2}}\,
  \frac{1+j}{\sqrt{2}}\, \frac{1-i}{\sqrt{2}}=
  \frac{1+k}{\sqrt{2}}\ .$$ It is easily checked that no two of $1,
  j$, and $\frac{1+j}{\sqrt{2}}$ lie in the same
  $(\Hat{\Gamma_1},\Hat{\Gamma_1})$ double-coset.  Hence $Z_{1,1}:=\{1,j,
  \frac{1+j}{\sqrt{2}}\}$ is a set of representatives of
$\Hat{\Gamma_1}\backslash
\Hat{\Gamma} /\Hat{\Gamma_1}$.

The elements $\z$ of $Z_{1,*}$ are listed in Table~\ref{zetas}, along
with the images $\phi(\z)\in \So(3)$ and $\pi_\z\in S_3.$ Since
$Z_{1,1}\subset Z_{1,*}$, a separate listing for $Z_{1,1}$ is not
needed. In Table~\ref{zetas} and henceforth, we write $\piid$ for the
identity permutation, and, for distinct $a,b\in \{1,2,3\}$, we write
$\pi_{ab}$ for the transposition $(ab)$, the permutation that just
interchanges $a$ and $b$.

\begin{table}
\caption{\label{zetas}{\footnotesize Representatives of the
    double-coset space $\Hat{\Gamma_1}\backslash\Hat{\Gamma}/
    \Hat{\Gamma_{\jtop}^0} = \Hat{\Gamma_1}\backslash\Hat{\Gamma}$,
    where $\Hat{\Gamma_1}=\Hat{\Gamma_{\J_1}^0}$.  A set of
    representatives of $\Hat{\Gamma_1}\backslash \Hat{\Gamma}/
    \Hat{\Gamma_1}$ is $\widehat{Z}_{1,1}:=\{1, j, \z_{j,+}\}$.}}
\begin{tabular}{c|c|c|c|c}
\hline
$\z\in\widehat{Z}_{1,*}$ & 1 & $j$ & $\z_{j,\e}:=\frac{1+\e  j}{\sqrt{2}}
\,,\e\in\{\pm 1\}$
& $\z_{k,\e}:=\frac{1+\e  k}{\sqrt{2}}\,,\e\in\{\pm 1\}$
\\ \hline
\raisebox{-1ex}{$\phi(\z)\in Z_{1,*}$} & \raisebox{-1ex}{$\mbbi$} &
\raisebox{-2ex}{$\left[ \begin{array}{rrr}
-1 & 0 & 0 \\
0 & 1 & 0 \\
0 & 0 & -1
\end{array}\right]$} &
\raisebox{-2ex}{$\left[ \begin{array}{rrr}
0 & 0 & \e \\
0 & 1 & 0 \\
-\e & 0 & 0
\end{array}\right]$} &
\raisebox{-2ex}{$\left[ \begin{array}{rrr}
0 & -\e & 0 \\
\e & 0 & 0 \\
0 & 0 & 1
\end{array}\right]$}\\
\hline
\raisebox{-1ex}{$\pi_\z\in S_3$} & \raisebox{-1ex}{$\piid$}
&\raisebox{-1ex}{$\piid$}
&\raisebox{-1ex}{$\pi_{13}$}
&\raisebox{-1ex}{$\pi_{12}$}\\
\hline
\end{tabular}
\end{table}

\begin{remark}
  {\rm As seen in Section \ref{sec:exsr3}, $\E_X$ has six connected
    components, each diffeomorphic to the circle $G_{\J_1}^0$. In
    Proposition \ref{prop:biject}, for general $p$ and $X$ we
    exhibited a bijection between $\Comp(\E_X)$ and the left-coset
    space $\tsp^+/\Gamma_{\J_D}^0$. For any group $G$ and subgroup
    $H$, the inversion map $G\to G$ induces a 1-1 correspondence
    between left $H$-cosets and right $H$-cosets, so (for general $p$
    and $X$), $\Comp(\E_X)$ is also in bijection with
    $\Gamma_{\J_D}^0\backslash \tsp^+$.  In our current $p=3,
    X\in\smid$ setting, the set
    $Z_{1,*}:=\{\phi(\z)\ :\ \z\in\widehat{Z}_{1,*}\}$ is a set of
    representatives of $\Gamma_1\backslash \ts_3^+/\Gamma_{\jtop} =
    \Gamma_1\backslash \ts_3^+$. The fact that right $\Gamma_1$-cosets
    appear here instead of left cosets is an artifact of our having
    chosen $X$, rather than $Y$, to lie in $\smid$.}
\end{remark}

\subsubsection{Hypercomplex reformulation of the continuous-parameter
minimization}\label{sec:hyper}

We now have
%
\begin{align}
&\dsr(X,Y)^2 = d_M(\E_X, \E_Y)^2 = \nonumber \\
&\left\{\begin{array}{l}
\displaystyle\min_{\z\in \widehat{Z}_{1,*}}\left(
4k\left[\min_{r_U\in S^1_\bfc}
\left\{\cos^{-1}\left|\Re\left(\bar{\z}\
\Bar{r_U}\ q_{U^{-1}V}\right)\right|\right\}\right]^2
+\|\log(\L_{\pi_\z}D^{-1})\|^2\right)
\\
\mbox{\hspace{1in} if  $Y\in\esstop$},
\\
\displaystyle\min_{\z\in \widehat{Z}_{1,1}}
\left(
4k\left[\min_{r_U, r_V\in S^1_\bfc}
\left\{\cos^{-1}\left|\Re\left(r_V\ \bar{\z}\
\Bar{r_U}\ q_{U^{-1}V}\right)\right|\right\}\right]^2
+\|\log(\L_{\pi_\z}D^{-1})\|^2\right)
\\ \label{fibdist10}
\mbox{\hspace{1in} if  $Y\in\smid$}.
\end{array}
\right.
\end{align}
To allow us to refer efficiently to the minimization-parameters in
\eqref{fibdist10} without too much separate notation for the two cases
$Y\in\esstop, Y\in \smid$, for both cases we will refer to the triple
$(\z,r_U,r_V)$, with the understanding that we always take $r_V=1$
when $Y\in \esstop$.

Recall that quaternions can be written in ``hypercomplex'' form: we
regard the complex numbers $\bfc$ as the subset $\{a+bi\}\subset \H$,
and write $x_0+x_1i +x_2j+x_3k = z + wj= z + j\bar{w},$ where
$z=x_0+x_1i, w=x_2+x_3i$. This gives us a natural identification
\be\label{natident} S^3_\H \longleftrightarrow S^3_{\bfc^2} :=
\{(z,w)\in\bfc^2\ :\ |z|^2+|w|^2=1\}.  \ee

\noi To perform the maximization of \eqref{fibdist7} and
\eqref{fibdist3} (in order to minimize the arc-cosines in \eqref{fibdist9-new}), we will write $q_{U^{-1}V}$
in hypercomplex form:
\be\label{hyper}
q_{U^{-1}V}= z+wj, \ \ (z,w)\in S^3_{\bfc^2}\ .
\ee

\noi Henceforth whenever we refer to the
quantities $z$ and $w$, they are regarded as functions of the pair
$(U,V)$, satisfying \eqref{hyper}, and with the pair $(z,w)$
determined only up to an overall sign.

Because the parameters $r_U, r_V$ in \eqref{fibdist7} and
\eqref{fibdist3} run over the unit circle in $\bfc$, it is easy to
maximize \eqref{fibdist7} and \eqref{fibdist3} explicitly for each
$\z$ in $\widehat{Z}_{1,*}$ and $\widehat{Z}_{1,1}$, respectively, and
then to maximize over $\z$. To express some of our answers, we define
the following quantities, which we may regard as functions of
the pair $(U,V)$:
\be\label{varphi}
\varphi:=\cos^{-1}(\max\{|z|,|w|\}),
\ee
%

\be\label{betas}
\b:= \frac{1}{2} \cos^{-1}(2 |\Re(\bar{z}w)|), \ \ 
\b':=  \frac{1}{2} \cos^{-1}(2 |\Im(\bar{z}w)|). 
\ee
\noi Note that $\varphi,\b$, and $\b'$ all lie in the interval $[0,
  \frac{\pi}{4}]$.

\begin{remark}\label{rem:zweq0}
{\rm Let $(U,D)\in\E_X, (V,\L)\in \E_Y$ be eigen-decompositions
 of $X,Y$ respectively. When $X$
  and $Y$ both lie in $\smid$, the ellipsoids corresponding to the
  matrices $X$ and $Y$ are surfaces of revolution. When $D$ and $\L$
  both lie in $\cald_{\J_1}$, the case in which
\be\label{cond:zweq0}
w=0 \ \mbox{\rm or}\  z=0
\ee
\noi in \eqref{hyper} has a simple geometric interpretation, and will
have special significance later in Theorem \ref{allMSR}. Observe that
$w=0$ if and only if $V=UR$ for some $R\in G_{\cald_{\J_1}}^0$, while
$z=0$ if and only if $V=UR$ for some $R\in G_{\cald_{\J_1}}^0j:=
\{R\phi(j)\ :\ R\in G_{\cald_{\J_1}}^0\}$. But $G_{\cald_{\J_1}}^0$
and $G_{\cald_{\J_1}}^0j$ are exactly the two connected components of
$G_{\cald_{\J_1}}$. Hence \eqref{cond:zweq0} is equivalent to $V=UR$
for some $R\in G_{\cald_{\J_1}}$, which in turn is equivalent to $V\L
V'=U\L U^{-1}$. Thus (when $D,\L\in \cald_{\J_1}$) the following are
all equivalent: (i) $\varphi=0$; (ii) $w=0$ or $z=0$; (iii)
simultaneously,
\be\label{simul}
X=UDU^{-1}\ \mbox{\rm and}\ Y=U\L U^{-1}\,;
\ee
\noi (iv) the ellipsoids of revolution corresponding to the matrices
$X$ and $Y$ have the same axis of symmetry.
The latter condition is
obviously intrinsic to the pair $(X,Y)$, independent of any choices of
eigen-decompositions.  Note also that when we want to find all MSSR curves from $X$ to $Y$, we do not need to express
these in terms of {\em arbitrary} eigen-decompositions with $D,\L\in
\cald_{\J_1}$; it suffices to use any that we find
convenient. Thus, given the eigen-decomposition $(U,D)$ of $X$, if
\eqref{cond:zweq0} (hence \eqref{simul}) holds we are free to replace
$V$ with $U$, in which case $U^{-1}V=\mbbi$ and $(z,w)=(\pm 1,0)$. We
will adopt the following convention:

\begin{Convention}\label{convention5.3}
If $D,\L\in \cald_{\J_1}$ and $U,V$ are such that $w=0$ or $z=0$, we
  replace $V$ with $U$, and replace $(z,w)$ with $(1,0)$.
We do not change $U$.
\end{Convention}
}
\end{remark}


\subsubsection{Closed-form formulas for $\Sym^+(3)$
  distances}\label{sec:distform3}

\begin{thm}\label{thm:dsr3}
Let $X,Y\in \Sym^+(3)$, $(U,D)\in\E_X$, and $(V,\L)\in \E_Y$. If $X\in
\smid$ assume $D\in\cald_{\J_1}$; if $Y\in
\smid$ assume $\L\in\cald_{\J_1}$. The distance $\dsr(X,Y)$ is given as
  follows.

(i) If $X, Y\in \esstop$, then
\bearray \nonumber \dsr(X,Y)^2&=&
\min_{g\in \ts_3^+}
  \left\{\frac{k}{2}\left\| \log(U^{-1}VP_g^{-1})\right\|^2
       +\left\|\log\left(D^{-1}(\pi_g\dotprod \L)\right)\right\|^2
      \right\}.
\label{fibdist4-3to3}
\eearray

 (ii) If $X\in \smid, Y\in \esstop$, then
\be\label{fibdist4-2to3}
\dsr(X,Y) =
  \min\left\{\ell_{\rm id}, \ell_{(13)}, \ell_{(12)}\right\},
\ee

\noi where
\begin{equation}
\label{fibdist4-2to3.piid}
\ell_{\rm id} =
\sqrt{4k\varphi^2  +\|\log(\L
D^{-1})\|^2}\ ,
\end{equation}
\begin{equation}
\label{fibdist4-2to3.pi13}
\ell_{(13)} =
\sqrt{4k\b^2
+\|\log(\L_{\pi_{13}} D^{-1})\|^2}\ ,
\end{equation}
\begin{equation}
\label{fibdist4-2to3.pi12}
\ell_{(12)} =
\sqrt{4k(\b')^2
+\|\log(\L_{\pi_{12}} D^{-1})\|^2}\ ,
\end{equation}


\noi and where $\phi,\b,\b' \in [0,\frac{\pi}{4}]$ are defined by
\eqref{hyper} and \eqref{varphi}--\eqref{betas}. Writing $D$ as ${\rm
  diag}(d_1,d_2,d_2)$ and $\L$ as ${\rm diag}(\l_1,\l_2,\l_3)$, we
also have the following comparisons of $\ell_{\rm id}, \ell_{(13)},$
and $\ell_{(12)}$:
\bearray
\label{ellcompar1}
\ell_{{\rm id}}^2 - \ell_{(13)}^2 &=& 4k(\varphi^2-\b^2)
-2\log\left(\frac{d_1}{d_2}\right) \log\left(\frac{\l_1}{\l_3}\right), \\
\label{ellcompar2}
\ell_{{\rm id}}^2 - \ell_{(12)}^2 &=& 4k(\varphi^2-(\b')^2)
-2\log\left(\frac{d_1}{d_2}\right) \log\left(\frac{\l_1}{\l_2}\right), \\
\label{ellcompar3}
\ell_{(13)}^2 - \ell_{(12)}^2 &=& 4k(\b^2-(\b')^2)
+2\log\left(\frac{d_1}{d_2}\right)\log\left(\frac{\l_2}{\l_3}\right).
\eearray

 (iii) If $X,Y$ both lie in $\smid$, then
\be\label{fibdist4}
\dsr(X,Y) = \min\{\ell_{\rm id}, \ell_{(13)}\},
\ee

\noi where
\bearray
\label{lidiii}
\ell_{{\rm id}}&=&\sqrt{4k\varphi^2 +\|\log(\L D^{-1})\|^2}\ ,
\\ \ell_{(13)} &=& \sqrt{4k(\frac{\pi}{4}-\varphi)^2 +\|\log(\L_{\pi_{13}}
  D^{-1})\|^2}\ .
\label{l12iii}
\eearray

\noi Writing $D$ as ${\rm diag}(d_1,d_2,d_2)$ and $\L$ as ${\rm
  diag}(\l_1,\l_2,\l_2)$, we also have the following comparison of
$\ell_{\rm id}$ and $\ell_{(13)}$:
\be\label{l1l2-alt1b}
\ell_{\rm id}^2-\ell_{(13)}^2 =  2k\pi(\varphi - \frac{\pi}{8})
-2\log\left(\frac{d_1}{d_2}\right)
\log\left(\frac{\l_1}{\l_2}\right).
\ee

 (iv) If $X\in \sbot:=\S_{\jbot}$ then $\dsr(X,Y)=\ddp(D,\L) =
\|\log(D^{-1}\L)\|$, regardless of which stratum $Y$ lies in.

\end{thm}

\pf
(i) Since $G_D^0=G_\L^0=\{I\}$ in this case, \eqref{fibdist4-3to3}
follows immediately from \eqref{fibdist8a} in Proposition
\ref{distprop}.

(ii) We use \eqref{fibdist10} to compute $\dsr(X,Y)$.
 We proceed by determining the ``inner'' minimum for
each $\z\in \Tilde{Z}_{1,*}$, and comparing the answers for the
different $\z$'s. For a given $\z$, minimizing
the arc-cosine in \eqref{fibdist10} is equivalent to maximizing
expression \eqref{fibdist7}. Below, we use the notation \eqref{hyper},
and for any nonzero $\xi\in\bfc$ we set $\hat{\xi}:=\xi/|\xi|$. Facts
used repeatedly in these calculations are that for all $\xi\in \bfc,$
(i) $\xi j$ and $\xi k$ are linear combinations of $j$ and $k$ with
real coefficients, hence are purely imaginary; and (ii) $j\xi =
\bar{\xi}j$ and $k\xi=\bar{\xi}k$. We then compute
\be\label{f1}
f_1(\z,r_U):=
\Re\left(\Bar{r_U}\, (z+wj) \right) =\left\{\begin{array}{ll}
\Re\left(\Bar{r_U}\ z\right) & \mbox{if}\ \z=1, \\
\Re\left(\Bar{r_U}\ w\right) & \mbox{if}\ \z=j, \\
\frac{1}{\sqrt{2}}\Re\left(\Bar{r_U}\ (z+\e w)\right) &
\mbox{if}\ \z=\z_{j,\e}\ , \\
\frac{1}{\sqrt{2}}\Re\left(\Bar{r_U}\ (z-\e iw)\right) &
\mbox{if}\ \z=\z_{k,\e}\ . \\
\end{array}\right.
\ee

\noi Hence for each $\z\in \widehat{Z}_{1,*}$, the value of
$\max_{r_U\in S^1_\bfc} \left|\Re\left(\bar{\z}\,
\Bar{r_U}\,(z+jw)\right)\right|$ is the entry in the last column of
the corresponding line of Table~\ref{classes}; let us denote this as
$|f_2(\z)|$, where $f_2(\z)\in\bfc$. The set
of elements $r_U\in S^1_\bfc$ at which the maximum is attained is
$\{\pm \Hat{f_2(\z)}\}$ if $f_2(\z)\neq 0,$ and all of $S^1_\bfc$ if
$f_2(\z)=0.$

Since $|z|^2+|w|^2=1$, we have $|z\pm w|^2=1\pm 2 \Re(\bar{z}w),
  |z\pm iw|^2=1\mp 2 \Im(\bar{z}w)$. Thus, grouping together the
    elements $\z\in\widehat{Z}_{1,*}$ corresponding to the same
    permutation $\pi_\z$, we have the following:
\bearray
|f_2(\z)| =\max_{r_U\in S^1_\bfc}
\left|\Re\left(\bar{\z}\, \Bar{r_U}\,(z+jw)\right)\right|
&=& \nonumber
\left\{
\begin{array}{ll}
\max\{|z|, |w|\}  &\mbox{if}\ \pi_\z=\piid\,,
\\ \sqrt{\frac{1+2|\Re(\bar{z}w)|}{2}}
&\mbox{if}\ \pi_\z=\pi_{13} \,, \\ \sqrt{\frac{1+2|\Im(\bar{z}w)|}{2}}
&\mbox{if}\ \pi_\z=\pi_{12}\,
\end{array}\right.\\
\label{max123-2to3_b}
\eearray

\noi (assuming $\z\in \widehat{Z}_{1,*}$). Equation
\eqref{fibdist4-2to3} now follows from the definitions
\eqref{varphi}--\eqref{betas}, equations \eqref{max123-2to3_b} and
\eqref{fibdist10}, and the identity $2 \cos^{-1}\sqrt{(1+x)/2} = \cos^{-1}x $ for $0 \le x  \le 1$.

Now let $a=\log d_1, b=\log d_2, c=\log\l_1, d=\log\l_2,
f=\log\l_3$. An easy calculation yields $\|\log(\L D^{-1})\|^2 -
\|\log((\pi_{13}\dotprod \L) D^{-1})\|^2 = -2(a-b)(c-f)$, from which
 \eqref{ellcompar1} follows. The derivations of
\eqref{ellcompar2} and \eqref{ellcompar3} are similar.

 (iii) We use the same strategy as in part (ii), but now with $\z$
ranging only over the set $\widehat{Z}_{1,1}= \{1,j,\z_{j,+}\}$, and
with $r_U,r_V$ both allowed to vary over $S^1_\bfc$. This time we find
\be\label{f3}
f_3(\z,r_U,r_V):=
\Re\left(r_V\,\bar{\z}\,\Bar{r_U}\ (z+wj) \right) =\left\{\begin{array}{ll}
\Re\left(r_V\ \Bar{r_U}\ z\right)& \mbox{if}\ \z=1, \\
\Re\left(r_Vr_U\bar{w}\right)& \mbox{if}\ \z=j, \\
\frac{1}{\sqrt{2}}\ \Re\left( r_V[\Bar{r_U}\, z+
  r_U\bar{w}]\right) &
\mbox{if}\ \z=\z_{j,+}\ . \\
\end{array}\right.
\ee

It is obvious from \eqref{f3} that
\be\label{max2-2_1j}
|f_4(\z)|:=\max_{r_U, r_V\in S^1_\bfc} |\Re\left(r_V\,\bar{\z}\, \Bar{r_U}\, (z+wj)
\right)|=\left\{ \begin{array}{ll} |z| & \mbox{if}\ \z=1, \\
|w| & \mbox{if}\ z=j,\end{array}\right.
\ee

\noi and that the pairs $(r_U,r_V)$ at which the maximum is achieved
are all those for which $r_V\Bar{r_U}=\Bar{\hat{z}}$ if $\z=1$, and
for which $r_Vr_U=\pm \hat{w}$ if $\z=j$. Thus, for these two $\z$'s,
  the set of maximizing pairs $(r_U,r_V)$ is the two circles' worth of
  pairs appearing in the triples $(\z,r_U,r_V)$ in the lines for
  classes ${\rm A}_1'$ and ${\rm A}_2'$ in Table~\ref{classes}, and
  the last entry of each line is the corresponding maximum value
  \eqref{max2-2_1j}.

Now consider $\z=\z_{j,+}=\frac{1+j}{\sqrt{2}}$. Since $r_U,r_V$ are
unit complex numbers, it is clear from \eqref{f3} that for all $r_U,r_V$,
\be\label{upbd}
|f_3(\z_{j,+}, r_U, r_V)|\leq |\Re\left(r_V\,\bar{\z}\, \Bar{r_U}\,
(z+wj) \right)| \leq \frac{|z|+|w|}{\sqrt{2}}.
\ee
First assume that $z\neq 0\neq w$. Then the upper bound on
$|f_3(\z_{j,+}, r_U, r_V)|$ in \eqref{upbd} will be achieved by a pair
$(r_U,r_V)$ if and only if
\be\label{max2.2}
(r_V\, \Bar{r_U}, r_Vr_U) = \pm (\Bar{\hat{z}}, \hat{w})\ .
\ee

\noi But \eqref{max2.2} is easily solved; the solution-set is exacly
the set of four pairs \linebreak $(\pm(\hat{w}\hat{z})^{1/2},\pm
(\hat{w}\Bar{\hat{z}})^{1/2})$ appearing Table~\ref{classes} for Class
B$'$. Thus the upper bound in \eqref{upbd} is actually the maximum
value of $|f_3(\z_{j,+}, \cdot, \cdot)|$.

Now assume that $w=0$ or $z=0$; we define the corresponding set of
pairs in $\E_X\times \E_Y$ (i.e. those pairs $((U\phi(r_U), D),
(V\phi(r_V)\phi(\z_{j,+})^T, \L_{\pi_{\z_{j,+}}}))$ for which
$(r_U,r_V)$ maximizes $|f_3(\z_{j,+},\cdot,\cdot)|$) to be Class C$'$.
If $w=0$ then $|z|=1$, and we need only maximize
$\frac{1}{\sqrt{2}}|\Re\left(r_V\ \Bar{r_U} z\right)|$.  Since
$|z|=1$, the maximum value is $\frac{1}{\sqrt{2}}$, and is achieved at
all pairs $(r_U,r_V)$ for which $r_V\Bar{r_U}=\pm \bar{z}$.
Similarly, if $z=0$ then $|w|=1$, and we need only maximize
$\frac{1}{\sqrt{2}}|\Re\left(r_Vr_U\bar{w}\right)|$. The maximum is
again $\frac{1}{\sqrt{2}}$, now achieved at all pairs $(r_U,r_V)$ for
which $r_Vr_U=\pm w$.

Hence, for $\z=\z_{j,+}$, whether or not $z$ and $w$ are both nonzero,
the right-hand side of \eqref{upbd} is the maximum value of
$|f_3(\z_{j,+}, \cdot, \cdot)|$.  But if $z\neq 0\neq w$ there are
only four maximizing pairs $(r_U,r_V)$, while if $w=0$ or $z=0$ there
are infinitely many.  As noted in Remark \ref{rem:zweq0}, in the
latter case we may replace $V$ with $U$, in which case $(z,w)=(1,0)$
(Convention~\ref{convention5.3}) and the maximizing pairs $(r_U,r_V)$ are exactly
those listed for Class C$'$ in Table~\ref{classes}.

Combining the maximum values computed for $\z=1, j,$ and $\z_{j,+}$,
we have the following: for $\z\in \widehat{Z}_{1,1}$,
\be\label{max123}
|f_4(\z)|=\max_{r_U, r_V\in S^1_\bfc}
\left|\Re\left(r_V\bar{\z}\, \Bar{r_U}\,(z+jw)\right)\right|
=\left\{
\begin{array}{ll}
\max\{|z|, |w|\} &\mbox{if}\ \pi_\z=\pi_{\rm id}, \\
\frac{1}{\sqrt{2}}(|z|+|w|) &\mbox{if}\ \pi_\z=\pi_{13}\ .
\end{array}\right.
\ee
The equality $|z|^2+|w|^2=1$ implies that
$\cos^{-1}\left(\frac{|z|+|w|}{\sqrt{2}}\right)=
 \frac{\pi}{4}
-\cos^{-1}(\max\{|z|,|w|\})$. This fact, combined with equations
\eqref{varphi}, \eqref{max123} and \eqref{fibdist10}, yields
\eqref{fibdist4}.

 (iv) In this case $G_D^0=\So(3),$ $\Gamma_{\J_D}^0=\ts_3^+$, and
the set $Z$ in Proposition \ref{distprop} has only one element $g$,
which we can take to be the identity. The set $\{UR_U: U\in
\Hat{G_D^0}\}$ is simply $\So(3)$, the inner minimum in
\eqref{fibdist8a} is 0, and $\dsr(X,Y)=\|\log (\L D^{-1}\|$.
\qedns

\begin{remark}[Insensitivity to choice of eigen-decompositions]\label{invariant}
{\rm By definition,  $\dsr(X,Y)$ cannot depend on the choice
  of pre-images $(U,D)\in F^{-1}(X) \linebreak = \E_X$, $(V,\L)\in F^{-1}(Y)=\E_Y$, that we have used to
  write down the formulas in Theorem \ref{thm:dsr3}.  However, the
  assumption that $D\in \cald_{\J_1}$ in the parts (ii) and (iii) of
  the theorem limits $(U,D)$ to particular pair of connected
  components of $\E_X$ out of the possible six.  A similar comment
  applies in part (iii) to the choice of $(V,\L)$.  So the individual
  numbers $\ell_{\rm id}, \ell_{(13)}, \ell_{(12)}$ on the right-hand
  sides of \eqref{fibdist4-2to3} and \eqref{fibdist4}, which represent
  distances between the connected component $[(U,D)]$ and the various
  connected components of $\E_Y$, may depend on the choice of $(U,D)$,
  but changing $(U,D)$ to a different pre-image of $X$ (not necessarily with
  $D\in \cald_{\J_1}$) must give us the same {\em set} of
  component-distances, and cannot change any of the the numbers
  $\ell_{\rm id}, \ell_{(13)}, \ell_{(12)}$ at all if the new pre-image
  is in the same connected component as the old. The latter {\em a
    priori} truth is reflected in the formulas given in Theorem
  \ref{thm:dsr3}.  Although the complex numbers $z,w$ in \eqref{hyper}
  depend on the choice of representatives $(U,D)\in \E_X$, $(V,\L)\in
  \E_Y$, when $D$ lies in $\cald_{\J_1}$ the quantities $|z|, |w|,$
  and $\bar{z}w$ depend only on the connected components $[(U,D)],
  [(V,\L)]$.  (Changing $(U,D)$ to $(UR,D)$, with $R\in
  G_{\cald_{\J_1}}^0$, changes $(z,w)$ to $(\xi z, \xi w)$ for some
  $\xi\in S^1_\bfc$; similarly if $\L\in \cald_{\J_1}$, then changing
  $(V,\L)$ to $(VR,\L)$, with $R\in G_{\cald_{\J_1}}^0$, changes
  $(z,w)$ to $(\xi z, \bar{\xi} w)$ for some $\xi\in S^1_\bfc$.) Thus
  when $X\in\smid$, and $Y\in \smid$ or $Y\in \esstop$, in
  \eqref{hyper}--\eqref{betas} we can regard $|z|, |w|,$ and
  $\bar{z}w$ as functions of a pair $([(U,D)],[(V,\L)])$ of connected
  components of fibers. Therefore the same is true of the quantities
  $\ell_{\rm id}, \ell_{(13)}, \ell_{(12)}$ in Theorem \ref{thm:dsr3}.
  (Of course, when $Y\in\esstop$, $[(V,\L)]=\{(V,\L)\}$.)  Furthermore,
  if $D\in \cald_{\J_1}$, then $[(U,D)]$ and $[(U\phi(j),D)]$ are the
  two connected components of $\E_X$ in $\So(3)\times
  \cald_{\J_1}$. Replacing $(U,D)$ by $(U\phi(j), D)$ has the effect
  of replacing $(z,w)$ by $\pm(\bar{w}, -\bar{z})$, which leaves the
  quantities $\varphi,\b,\b'$ in \eqref{varphi}--\eqref{betas}
  unchanged, and hence leaves each of the numbers $\ell_{\rm id},
  \ell_{(13)}, \ell_{(12)}$ in
  \eqref{fibdist4-2to3.piid}--\eqref{fibdist4-2to3.pi12}
  and\eqref{lidiii}--\eqref{l12iii} unchanged. The fact that replacing
  $(U,D)$ by other pre-images of $X$ cannot change the {\em set}
  $\{\ell_{\rm id}, \ell_{(13)}, \ell_{(12)}\}$ is also reflected,
  later in Theorem \ref{allMSR}, by the symmetry of the last column of
  Table~\ref{tableS3} under permutations of $\ell_{\rm id},
  \ell_{(13)}, \ell_{(12)}$.  }
\end{remark}

\setcounter{equation}{0}
 \section{MSSR curves for $\Sym^+(3)$ in the nontrivial cases}\label{sec:mssr3}

Recall that for $X,Y\in\sympp$,  $\M(X,Y)$ denotes the set of
all MSSR curves from $X$ to $Y$. In this section, for $p=3$ we determine the set $\M(X,Y)$ for all $X\in \smid, Y\in \smid\,\union
\,\esstop$ (what we are calling the ``nontrivial cases").

\subsection{Explicit characterization of all MSSR curves in the nontrivial cases}\label{sec:6.1}
%
%
%
For any $X,Y\in \Sym^+(3)$ and $(U,D)\in \E_X$, Proposition
\ref{distprop} assures us that every MSSR curve from $X$ to $Y$
corresponds to some minimal pair whose first element lies in the
connected component $[(U,D)]$.  When $X\in\smid$, by keeping
track of the triples $(\z,r_U,r_V)$ at which the minimum values in
\eqref{fibdist10} are achieved, we can find all the minimal pairs in
$\E_X\times \E_Y$ whose first point lies in $[(U,D)]$ of $\E_X$.  The
following corollary of Proposition \ref{prop:MSRequal} will allow us
to tell when the MSSR curves corresponding to two such minimal pairs
are the same.

\begin{cor} \label{cor:MSRequalp3}
Hypotheses as in Proposition \ref{prop:MSRequal}, but additionally
assume that $p=3$.  For $i=1,2$ let $r_{U,i}, r_{V,i}$, and $\z_i$ be
preimages of $R_{U,i}, R_{V,i}$, and $g_i$ under $\phi$. Then
$\chi_1=\chi_2$ if and only if (i)$'$
\be \label{a2eqa1} r_{V,2}\,\Bar{\z_2}\,\Bar{r_{U,2}} = \pm
r_{V,1}\,\Bar{\z_1}\,\Bar{r_{U,1}} \ee
and (ii)$'$ there exist $\z\in \Hat{\Gamma}, r\in
\Hat{G_{D,\L_1}^0}$ such that
\begin{eqnarray}
\label{d2pdd1-final}
D &=& \pi_\z\dotprod D, \\
\label{l2pdl1-final}
\L_2 &=& \pi_\z\dotprod \L_1, \\
\mbox{and} \ \ \ \Bar{r_{U,1}}\, r_{U,2}  &=& r\bar{\z}.
\label{rppp-final}
\end{eqnarray}

\end{cor}

\pf  This follows immediately from Proposition \ref{prop:MSRequal}.\qedns

The classification we will give of MSSR curves involves six classes of
scaling-rotation curves when $Y\in\esstop$, and four classes when $Y\in
\smid$. Not all of these classes occur for a given $X$ and $Y$, and
when they do occur they are not necessarily minimal. The (potentially)
minimal pairs giving rise to the various classes of scaling-rotation
curves can be described in terms of the data $z,w$ and the triple
$(\z,r_U,r_V)$. Our names for these classes of pairs and curves, and
the data $(\z,r_U,r_V)$ corresponding to each class, are listed in
Table~\ref{classes}.  For the $\z$ appearing in each line of the
table, the accompanying values of $(r_U,r_V)$ are all those that
minimize the arc-cosine term in the corresponding line of
\eqref{fibdist10}, provided that any unit complex number
$\hat{\xi}=\xi/|\xi|$ appearing in that line's indicated formula for
$(r_U,r_V)$ is defined (i.e. provided $\xi\neq 0$); see the proof of
Theorem \ref{allMSR} in later in this section.  The corresponding
pairs in $\E_X\times \E_Y$ determine a class $\M_l([(U,D)],[(V,\L)])$
of scaling-rotation curves, where $l$ is the corresponding class-name
appearing in Table~\ref{classes}.  As the notation suggests, each
class of curves depends only on the connected components $[(U,D)],
[(V,\L)]$ in $\E_X, \E_Y$, although the data $(r_U,r_V)$ for a given
$\z$ will depend fully on the matrices $U,V$.  For the
scaling-rotation curves in $\M_l([(U,D)],[(V,\L)])$ to be {\em
  minimal} there are restrictions on the component-pair
$([(U,D)],[(V,\L)])$, reflected by restrictions on $z$ and $w$ that
depend only on these connected components; e.g. for Class B$_1$ to be
minimal we need $\Re(\bar{z}w)\geq 0$, and for Class ${\rm A}_1'$ to
be minimal we need $|z|\geq |w|$. The full set of restrictions can be
read off from Tables~\ref{tableS3} and \ref{tableS2}, which are part
of Theorem \ref{allMSR} below.


\begin{remark} {\rm In our application of Corollary \ref{cor:MSRequalp3}
    to the proof of Theorem \ref{allMSR} below, we will have $D\in
    \cald_{\J_1}$, and hence the quaternions $r_{U,i}$, $r_{V,i}$ in
    \eqref{a2eqa1} and \eqref{rppp-final} will lie in $S^1_\bfc$. Note
    also that the only permutations $\pi$ for which $\pi\dotprod D = D$
    are the identity and the transposition $\pi_{23}$. Thus the only
    $\z$'s that can satisfy \eqref{d2pdd1-final} are those that lie in
    the group $\Hat{\Gamma_1} =\left\{\pm 1, \pm i, \frac{\pm 1 \pm
      i}{\sqrt{2}}\right\} = \{e^{2\pi im/8}: m\in {\mathbb Z}\}
    \subset S^1_\bfc$. However, in general the $\z_i$ in
    \eqref{a2eqa1} need not lie in $\bfc$.}
\end{remark}

\begin{table}
\caption{\label{classes} {\footnotesize Names and data for classes of
    pairs in $\E_X\times \E_Y$ that, for {\em some} $X$ in $\smid$ and
    {\em some} $Y$ in $\esstop$ or $\smid$, determine at least one MSSR
    curve from $X$ to $Y$. For any nonzero $\xi\in \bfc$, $\hat{\xi}$
    is the unit complex number $\xi/|\xi|$, and $\xi^{1/2}$ is an
    arbitrary choice of one of the two square roots of $\xi$. Wherever
    a number of the form $\hat{\xi}$ appears in this table, the
    corresponding class is defined only for $\xi\neq 0$. In case ${\rm C'}$
    we have used Convention~\ref{convention5.3} to simplify this line of the table;
    for general definition of Class ${\rm C'}$ see the proof of Theorem
    \ref{thm:dsr3} in Section \ref{sec:distform3}.  The last column of
    the table is included for the proof of Theorem \ref{thm:dsr3}.}}
{\footnotesize
\begin{tabular}{cclc}
 & Class & $\{(\z,r_U,r_V)\}$ & $\left|\Re\left(r_V\, \bar{\z}\,
\Bar{r_U}\, (z+wj)\right)\right|$
\\ \hline
\raisebox{-1ex}{For $Y\in\esstop$:}
& \raisebox{-1ex}{${\rm A}_1$}
& \raisebox{-1ex}{$\{(1,\pm \hat{z},1)\}$} & \raisebox{-1ex}{$|z|$}
\\
& \raisebox{-1ex}{${\rm A}_2$}
& \raisebox{-1ex}{$\{(j,\pm \hat{w},1)\}$}
& \raisebox{-1ex}{$|w|$}\\
& \raisebox{-1ex}{${\rm B}_1$}
& \raisebox{-1ex}{$\{(\z_{j,+},\pm \Hat{(z+w)},1)\}$}
& \raisebox{-1ex}{$\frac{1}{\sqrt{2}}|z+w|$}\\
& \raisebox{-1ex}{${\rm B}_2$}
& \raisebox{-1ex}{$\{(\z_{j,-},\pm \Hat{(z-w)},1)\}$}
& \raisebox{-1ex}{$\frac{1}{\sqrt{2}}|z-w|$}\\
& \raisebox{-1ex}{${\rm C}_1$}
& \raisebox{-1ex}{$\{(\z_{k,+},\pm \Hat{(z-iw)},1)\}$}
& \raisebox{-1ex}{$\frac{1}{\sqrt{2}}|z-iw|$}\\
& \raisebox{-1ex}{${\rm C}_2$}
& \raisebox{-.5ex}{$\{(\z_{k,-},\pm \Hat{(z+iw)},1)\}$}
& \raisebox{-1ex}{$\frac{1}{\sqrt{2}}|z+iw|$}\\
\hline
\raisebox{-1ex}{For $Y\in\smid$:}
& \raisebox{-1ex}{${\rm A}_1'$}
& \raisebox{-1ex}{$\{(1,r,\pm  r\Bar{\hat{z}})\ : r\in S^1_\bfc\}$}
& \raisebox{-1ex}{$|z|$}\\
& \raisebox{-1ex}{${\rm A}_2'$}
& \raisebox{-1ex}{$\{(j,r,\pm  r\hat{w})\ : r\in S^1_\bfc\}$}
& \raisebox{-1ex}{$|w|$}\\
& \raisebox{-1ex}{${\rm B}'$}
& \raisebox{-1ex}{$\{(\z_{j,+},\pm(\hat{w}\hat{z})^{1/2},\pm
  (\hat{w}\Bar{\hat{z}})^{1/2})\}$}
& \raisebox{-1ex}{$\frac{1}{\sqrt{2}}(|z|+|w|)$}\\
&&\raisebox{-.5ex}{\mbox{\hspace{3ex}} \mbox{(all
  sign-combinations allowed)}} & \\
& \raisebox{-1ex}{${\rm C}'$}
& \raisebox{-1ex}{$\{(\z_{j,+}, r, \pm r)
   : r\in S^1_\bfc\}$
  \ if $(z,w)=(1,0)$; }
& \raisebox{-1ex}{$\frac{1}{\sqrt{2}}$}\\
 && \raisebox{-.5ex}{\mbox{\hspace{3ex}} class defined if $w=0$ or $z=0$}\\
 && \raisebox{-.5ex}{\mbox{\hspace{3ex}} but left undefined otherwise.}\\
\hline

\end{tabular}
}
\end{table}

\begin{thm}\label{allMSR}
Assume that $X\in\smid$, $Y\in\esstop\ \union\ \smid$, $(U,D) \in \E_X$,
$(V,\L)\in\E_Y$, and that the first two diagonal entries of each of
the matrices $D,\L$ are distinct. (Thus $D\in \cald_{\J_1}$, and if $Y\in
\smid$ then $\L\in\cald_{\J_1}$.) Let $l$ stand for the class-names
in Table~\ref{classes}, and abbreviate $\M_l([(U,D)], [V,\L])$
as $\M_l$.

 (i) Except for $\M_{{\rm C}'}$, every class $\M_l$, when defined,
consists of a single curve $\chi_l=\chi_l([(U,D)], [V,\L])$. The class
$\M_{{\rm C}'}$, which we define only when $w=0$ or $z=0$, is an
infinite family of scaling-rotation curves, in natural one-to-one
correspondence with a circle. The
class $\M_{{\rm C}'}$ does not depend on the choice of components
$[(U,D)], [(V,\L)] \subset \So(3)\times \cald_{\J_1}$, so can
unambiguously be written as $\M_{{\rm C}'}(X,Y)$.

 (ii) For any data-triple $(\z,r_U,r_V)$ as in Table~\ref{classes},
let $R_U=\phi(r_U), R_V=\phi(r_V)$. For both $Y\in\esstop$ and $Y\in
\smid$, the pair
\be\label{minpr}
((UR_U, D), (VR_V\phi(\z)^{-1}, \L_{\pi_\z}))\in \E_X\times \E_Y
\ee

\noi is a minimal pair in each case listed in Tables~\ref{tableS3} and
\ref{tableS2}, with $(UR_U,D)$ lying in the connected component
$[(U,D)]$ of $\E_X$.  Conversely, every minimal pair in $\E_X\times
\E_Y$ whose first point lies in $[(U,D)]$ is given by the data in
Table~\ref{classes} and either Table~\ref{tableS3} or Table
\ref{tableS2}.

 (iii) For $Y\in \esstop$, depending on the value of $Y$ the set
$\M(X,Y)$ can consist of one, two, three, or four curves, as detailed
in Table~\ref{tableS3}.  In Tables~\ref{tableS3} and \ref{tableS2},
note that ``$|\M(X,Y)|=1$'' means precisely that there is a {\em
  unique} MSSR curve from $X$ to $Y$.

\begin{table}

\caption{\label{tableS3} {
\footnotesize The set $\M(X,Y)$ of minimal smooth scaling-rotation
curves from $X$ to $Y$ when $X$ has exactly two distinct eigenvalues
and $Y$ has three.
Data-combinations that are mutually exclusive are not shown (e.g. if
$\ell_{\rm id} = \ell_{(13)} <\ell_{(12)}$, it is impossible to have
$|z|-|w|= 0=\Re(\bar{z}w)$).  In the subcase of $\ell_{{\rm id}}
=\ell_{(13)} = \ell_{(12)}$ in which $|z|=|w|$, the hypothesis
$\Re(\bar{z}w)\neq 0\neq \Im(\bar{z}w)$ is redundant; it is already
implied by the case/subcase hypotheses. (This follows from Theorem
\ref{thm:dsr3}; see the proof of Theorem \ref{allMSR}.)}}}

{\footnotesize
\mbox{\hspace{-.5in}}
\begin{tabular}{cccc}
Case & Subcase & $\M(X,Y)$ & $|\M(X,Y)|$
\\ \hline
   \raisebox{-1ex}{$\ell_{\rm id}<\min\{\ell_{(13)},\ell_{(12)}\}$}
 &    \raisebox{-1ex}{$|z|\neq |w|$} &
\raisebox{-1ex}{$\{\chi_{A_1}\}$ if $|z|>|w|$; $\{\chi_{A_2}\}$ if $|z|<|w|$}
 &    \raisebox{-1ex}{1} \\
& \raisebox{-1ex}{$|z|=|w|$}
& \raisebox{-.5ex}{$\{\chi_{A_1}, \chi_{A_2}\}$}
& \raisebox{-1ex}{2}\\ \hline

\raisebox{-1ex}{$\ell_{(13)}<\min\{\ell_{{\rm id}},\ell_{(12)}\}$}
& \raisebox{-1ex}{$\Re(\bar{z}w)\neq 0$}
& \raisebox{-1ex}{$\{\chi_{B_1}\}$ if $\Re(\bar{z}w)>0$;
  $\{\chi_{B_2}\}$ if $\Re(\bar{z}w)<0$}
 & \raisebox{-1ex}{1}\\
&\raisebox{-1ex}{$\Re(\bar{z}w)=0$}
& \raisebox{-.5ex}{$\{\chi_{B_1}, \chi_{B_2}\}$}
& \raisebox{-1ex}{2}\\
\hline

\raisebox{-1ex}{$\ell_{(12)}<\min\{\ell_{{\rm id}},\ell_{(13)}\}$}
& \raisebox{-1ex}{$\Im(\bar{z}w)\neq 0$}
& \raisebox{-1ex}{$\{\chi_{C_1}\}$ if $\Im(\bar{z}w)>0$;
  $\{\chi_{C_2}\}$ if $\Im(\bar{z}w)<0$}
 & \raisebox{-1ex}{1}\\
&\raisebox{-1ex}{$\Im(\bar{z}w)=0$}
& \raisebox{-.5ex}{$\{\chi_{C_1}, \chi_{C_2}\}$}
& \raisebox{-.5ex}{2}\\
\hline

\raisebox{-1ex}{$\ell_{\rm id} = \ell_{(13)} < \ell_{(12)}$}
& \raisebox{-1ex}{$|z|-|w|\neq 0\neq \Re(\bar{z}w)$}
&
\begin{tabular}{c}
\raisebox{-1ex}{$\{\chi_{A_m}, \chi_{B_n}\}$} \\
\raisebox{-1ex}{$m=1$ (resp. 2) if $|z|-|w|>0$ (resp. $<0$),}\\
\raisebox{-1ex}{$n=1$ (resp. 2) if $\Re(\bar{z}w)>0$ (resp. $<0$)}
\end{tabular}

& \raisebox{-1ex}{2}\\

& \raisebox{-1ex}{$|z|-|w|=0\neq \Re(\bar{z}w)$}
&
\begin{tabular}{c}
\raisebox{-1ex}{$\{\chi_{A_1}, \chi_{A_2}, \chi_{B_n}\}$}\\
\raisebox{-1ex}{$n=1$ (resp. 2) if $\Re(\bar{z}w)>0$ (resp. $<0$)}
\end{tabular}
& \raisebox{-1ex}{3}\\
& \raisebox{-1ex}{$|z|-|w|\neq 0 =\Re(\bar{z}w)$}
&
\begin{tabular}{c}
\raisebox{-1ex}{$\{\chi_{A_m}, \chi_{B_1}, \chi_{B_2}\}$}\\
\raisebox{-1ex}{$m=1$ (resp. 2) if $|z|-|w|>0$ (resp. $<0$)}
\end{tabular}
& \raisebox{-1ex}{3}\\
\hline

\raisebox{-1ex}{$\ell_{\rm id} = \ell_{(12)} < \ell_{(13)}$}
& \raisebox{-1ex}{$|z|-|w|\neq 0\neq \Im(\bar{z}w)$}
&
\begin{tabular}{c}
\raisebox{-1ex}{$\{\chi_{A_m}, \chi_{C_n}\}$}\\
\raisebox{-1ex}{$m=1$ (resp. 2) if $|z|-|w|>0$ (resp. $<0$),}\\
\raisebox{-1ex}{$n=1$ (resp. 2) if $\Im(\bar{z}w)>0$ (resp. $<0$)}
\end{tabular}
& \raisebox{-1ex}{2}\\

& \raisebox{-1ex}{$|z|-|w|=0\neq\Im(\bar{z}w)$}
& \begin{tabular}{c}
\raisebox{-1ex}{$\{\chi_{A_1}, \chi_{A_2},\chi_{C_n}\}$}\\
\raisebox{-1ex}{$n=1$ (resp. 2) if $\Im(\bar{z}w)>0$ (resp. $<0$)}
\end{tabular}
& \raisebox{-1ex}{3}\\

& \raisebox{-1ex}{$|z|-|w|\neq 0 =\Im(\bar{z}w)$}
& \begin{tabular}{c}
\raisebox{-1ex}{$\{\chi_{A_m}, \chi_{C_1},\chi_{C_2}\}$}\\
\raisebox{-1ex}{$m=1$ (resp. 2) if $|z|-|w|>0$ (resp. $<0$)}
\end{tabular}
& \raisebox{-1ex}{3}\\
\hline

\raisebox{-1ex}{$\ell_{(13)} =\ell_{(12)} < \ell_{\rm id}$}
& \raisebox{-1ex}{$\Re(\bar{z}w)\neq 0\neq \Im(\bar{z}w)$}
& \begin{tabular}{c}
\raisebox{-1ex}{$\{\chi_{B_m}, \chi_{C_n}\}$}\\
\raisebox{-1ex}{$m=1$ (resp. 2) if $\Re(\bar{z}w)>0$ (resp. $<0$),}\\
\raisebox{-1ex}{$n=1$ (resp. 2) if $\Im(\bar{z}w)>0$ (resp. $<0$)}
\end{tabular}
& \raisebox{-1ex}{2}\\

& \raisebox{-1ex}{$\Re(\bar{z}w)= 0\neq \Im(\bar{z}w)$}
& \begin{tabular}{c}
\raisebox{-1ex}{$\{\chi_{B_1}, \chi_{B_2}, \chi_{C_n}\}$}\\
\raisebox{-1ex}{$n=1$ (resp. 2) if $\Im(\bar{z}w)>0$ (resp. $<0$)}
\end{tabular}
& \raisebox{-1ex}{3}\\

& \raisebox{-1ex}{$\Re(\bar{z}w)\neq 0= \Im(\bar{z}w)$}
& \begin{tabular}{c}
\raisebox{-1ex}{$\{\chi_{B_n}, \chi_{C_1}, \chi_{C_2}\}$}\\
\raisebox{-1ex}{$n=1$ (resp. 2) if $\Re(\bar{z}w)>0$ (resp. $<0$)}
\end{tabular}
& \raisebox{-1ex}{3}\\
\hline

\raisebox{-1ex}{$\ell_{\rm id}=\ell_{(13)} =\ell_{(12)}$}
&
\begin{tabular}{c}
\raisebox{-1ex}{$|z|\neq |w|$ \ and}\\
\raisebox{-1ex}{$\Re(\bar{z}w)\neq 0 \neq \Im(\bar{z}w)$}
\end{tabular}
&
\begin{tabular}{c}
\raisebox{-1ex}{$\{\chi_{A_1}, \chi_{B_m}, \chi_{C_n}$\}}\\
\raisebox{-1ex}{$l=1$ (resp. 2) if $|z|-|w|>0$ (resp. $<0$),}\\
\raisebox{-1ex}{$m=1$ (resp. 2) if $\Re(\bar{z}w)>0$ (resp. $<0$),}\\
\raisebox{-1ex}{$n=1$ (resp. 2) if $\Im(\bar{z}w)>0$ (resp. $<0$)}
\end{tabular}
& \raisebox{-1ex}{3}\\

& \begin{tabular}{c}
\raisebox{-1ex}{$|z|=|w|$ \ and}\\
\raisebox{-1ex}{$\Re(\bar{z}w)\neq 0  \neq \Im(\bar{z}w)$}
\end{tabular}
&
\begin{tabular}{c}
\raisebox{-1ex}{$\{\chi_{A_1}, \chi_{A_2}, \chi_{B_m},\chi_{B_n}$\},}\\
\raisebox{-1ex}{$m=1$ (resp. 2) if $\Re(\bar{z}w)>0$ (resp. $<0$),}\\
\raisebox{-1ex}{$n=1$ (resp. 2) if $\Im(\bar{z}w)>0$ (resp. $<0$)}
\end{tabular}
& \raisebox{-1ex}{4}\\

& \begin{tabular}{c}
\raisebox{-1ex}{$|z|\neq |w|$ \ and}\\
\raisebox{-1ex}{$\Re(\bar{z}w)= 0 \neq \Im(\bar{z}w)$}
\end{tabular}
&
\begin{tabular}{c}
\raisebox{-1ex}{$\{\chi_{A_m}, \chi_{B_1},\chi_{B_2}, \chi_{C_n}\}$} \\
\raisebox{-1ex}{$m=1$ (resp. 2) if $|z|-|w|>0$ (resp. $<0$),}\\
\raisebox{-1ex}{$n=1$ (resp. 2) if $\Im(\bar{z}w)>0$ (resp. $<0$)}
\end{tabular}
& \raisebox{-1ex}{4}\\

& \begin{tabular}{c}
\raisebox{-1ex}{$|z|\neq |w|$ \ and}\\
\raisebox{-1ex}{$\Re(\bar{z}w)\neq 0= \Im(\bar{z}w)$}
\end{tabular}

&\begin{tabular}{c}
\raisebox{-1ex}{$\{\chi_{A_m}, \chi_{B_n},\chi_{C_1}, \chi_{C_2}\}$} \\
\raisebox{-1ex}{$m=1$ (resp. 2) if $|z|-|w|>0$ (resp. $<0$),}\\
\raisebox{-1ex}{$n=1$ (resp. 2) if  $\Re(\bar{z}w)>0$ (resp. $<0$)}
\end{tabular}
& \raisebox{-1ex}{4}\\
\hline

\end{tabular}

\end{table}

 (iv) For $Y\in \smid$, depending on the value of $Y$ the set
$\M(X,Y)$ can consist of one, two, three, or infinitely many curves,
as detailed in Table~\ref{tableS2}. When $|z|\geq |w|$ (respectively,
$|z|\leq |w|$), all minimal pairs in Class ${\rm A}_1'$ (resp. ${\rm
  A}_2'$) determine the same MSSR curve, so to write down this curve
it suffices to take $r=1$ in the data-triple for this class in
Table~\ref{classes}.

\begin{table}
\caption{\label{tableS2} {\footnotesize The set $\M(X,Y)$ of minimal
    smooth scaling-rotation curves from $X$ to $Y$ when each of $X$
    and $Y$ has exactly two distinct eigenvalues.  See text for
    notation.}}
{\footnotesize
\begin{tabular}{cccc}
Case & Subcase & $\M(X,Y)$ & $|\M(X,Y)|$
\\ \hline
   \raisebox{-1ex}{$\ell_{\rm id}<\ell_{(12)}$}
 & \raisebox{-1ex}{$|z|\neq|w|$}
&    \raisebox{-1ex}{$\{\chi_{{\rm A}_1'}\}$ if $|z|>|w|$;
$\{\chi_{{\rm A}_2'}\}$ if $|z|<|w|$}
 &    \raisebox{-1ex}{1} \\
& \raisebox{-1ex}{$|z|=|w|$}
& \raisebox{-1ex}{$\{\chi_{{\rm A}_1'}, \chi_{{\rm A}_2'}\}$}
& \raisebox{-1ex}{2}\\ \hline

\raisebox{-1ex}{$\ell_{\rm id}>\ell_{(12)}$}
& \raisebox{-1ex}{$z\neq 0\neq w$} & \raisebox{-1ex}{$\{\chi_{{\rm B}'}\}$}
 & \raisebox{-1ex}{1}\\
&\raisebox{-1ex}{$w=0$ or $z=0$} & \raisebox{-1ex}{$\M_{{\rm C}'}$}
& \raisebox{-1ex}{$\infty$}\\
\hline

\raisebox{-1ex}{$\ell_{\rm id}=\ell_{(12)}$}
& \raisebox{-1ex}{$0\neq |z|\neq |w|\neq 0$}
& \raisebox{-1ex}{$\{\chi_{{\rm A}_1'}, \chi_{{\rm B}'}\}$ if
  $|z|>|w|$;
$\{\chi_{{\rm A}_2'}, \chi_{{\rm B}'}\}$ if $|z|<|w|$}
& \raisebox{-1ex}{2}\\

& \raisebox{-1ex}{$|z|=|w|$}
& \raisebox{-1ex}{$\{\chi_{{\rm A}_1'}, \chi_{{\rm A}_2'},
\chi_{{\rm B}'}\}$}
& \raisebox{-1ex}{3} \\

& \raisebox{-1ex}{$w=0$ or $z=0$}
& \raisebox{-1ex}{$\{\chi_{{\rm A}_1'}\}\ \union\ \M_{{\rm C}'}$ if
  $w=0$;
$\{\chi_{{\rm A}_2'}\}\ \union\ \M_{{\rm C}'}$ if $z=0$}
& \raisebox{-1ex}{$\infty$}\\

\hline

\end{tabular}  }
\end{table}

\end{thm}

\begin{remark}[Symmetries in Theorem \ref{thm:dsr3} and Tables
    \ref{tableS3} and \ref{tableS2}]\label{uphij} {\rm As noted
    in Remark \ref{invariant}, when $(U,D)$ is a pre-image (eigen-decomposition) of $X$
    with $D\in \cald_{\J_1}$, replacing $(U,D)$ by the pre-image
    $(U\phi(j), D)$ in the other connected component of $\E_X$ in
    $\So(3)\times\cald_{\J_1}$ has the effect of replacing $(z,w)$ by
    $(z_{\rm new}, w_{\rm new})=\pm(\bar{w}, -\bar{z})$. Observe that
    $\Bar{z_{\rm new}}\ w_{\rm new} = -\bar{z}w$, and that if
    $|z|<|w|$ then $|z_{\rm new}|> |w_{\rm new}|$. Thus when $X\in
    \smid$ we can always choose our pair of pre-images $(U,D), (V,\L)$
    (with $D\in \cald_{\J_1}$) to satisfy $|z|\geq |w|$, or
    $\Re(\bar{z}w)\geq 0$, or $\Im(\bar{z}w)\geq 0$ (though not
    necessarily more than one of these inequalities at the same
    time). This explains the ``symmetry'' in Tables~\ref{tableS3} and
    \ref{tableS2} under interchange of $|z|$ and $|w|$ and under
    sign-changes of $\Re(\bar{z}w)$ and $\Im(\bar{z}w)$: we have
    $\chi_{{\rm A}_1}([(U\phi(j),D)], [(V,\L)]) = \chi_{{\rm
        A}_2}([(U,D], [(V,\L)])$, and a similar relation for the
    class-pairs $({\rm B}_1, {\rm B}_2), ({\rm C}_1, {\rm C}_2)$, and
    $({\rm A}_1', {\rm A}_2')$.

}
\end{remark}

\begin{remark}[Condition for curves in family $\M_{{\rm C}'}$ to be
minimal]\label{rem:proob} {\rm The \linebreak
conditions under which the infinite
    family $\M_{{\rm C}'}$ arises in Table~\ref{tableS2}---namely,
    $\ell_{\rm id} \geq \ell_{(12)}$ and either $w=0$ or $z=0$---can
    be described more explicitly and geometrically in terms of the
    ellipsoids of revolution corresponding to $X$ and $Y$. Recall from
    Remark \ref{rem:zweq0} that ``$w=0 \ \mbox{or}\ z=0$'' is
    equivalent to the condition that these ellipsoids have the same
    axis of symmetry, and to the condition $\varphi=0$. But
    \eqref{l1l2-alt1b} shows that when $\varphi=0$, the condition
    $\ell_{\rm id} \geq \ell_{(13)}$ is equivalent to
\be\label{proob}
-\log\left(\frac{d_1}{d_2}\right)
\log\left(\frac{\l_1}{\l_2}\right)\geq
k\,\frac{\pi^2}{8}\ .
\ee

\noi In particular, $\log\left(\frac{d_1}{d_2}\right)$ and
$\log\left(\frac{\l_1}{\l_2}\right)$ must have opposite signs for
\eqref{proob} to hold, so one of the ellipsoids must be prolate and
the other oblate.  Conversely, given two ellipsoids of revolution with
the same axis of symmetry, one prolate and the other oblate, if their
``prolateness-oblateness product'' is sufficiently large---i.e. if
\eqref{proob} holds---then the set of MSSR curves from $X$ to $Y$ will
include the 1-parameter family $\M_{{\rm C}'}$. The proof below
of Theorem \ref{allMSR}(i) shows that for such $X$ and
$Y$, a choice of orientation of the common axis
of symmetry naturally determines a continuous one-to-one
correspondence between the ``equator'' of $X$ (or $Y$) and the family
$\M_{{\rm C}'}$.   For a graphical example
illustrating several members of the family $\M_{{\rm C}'}$ as evolutions of
the $X$-ellipsoid to the $Y$-ellipsoid, see Fig.~\ref{fig:Cp} in Section~\ref{sec:appendix-MSSR3}.}
\end{remark}

\bsn{\bf Proof of Theorem \ref{allMSR}.}
(i) By definition, each curve-class $\M_l$ is a set of (not
necessarily minimal) smooth scaling-rotation curves corresponding to
pairs $((U\phi(r_U),D), (V\phi(r_V)\phi(\z)^{-1},
\pi_\z\dotprod\L))\in\E_X\times\E_Y$ (where we set $r_V=1$ if $Y\in
\esstop$) for which $r_U$ maximizes the function $|f_1(\z, \cdot)|$
given by \eqref{f1} if $Y\in \esstop$, or for which $(r_U,r_V)$
maximizes the function $|f_3(\z, \cdot)|$ given by \eqref{f3} if $Y\in
\smid$.  In the proof of Theorem \ref{thm:dsr3} we established that
Table~\ref{classes} lists all the corresponding triples
$(\z,r_U,r_V)$, with the exception that for class C$'$ we followed
Convention~\ref{convention5.3} and listed the corresponding triples only for the case
$(z,w)=(1,0)$ (see Remark \ref{rem:zweq0}). For $i\in \{1,2\}$ let
$(\z,r_{U,i},r_{V,i})$ be two such triples listed in
Table~\ref{classes} corresponding to the same class $\M_l$, and let
$\chi_i$ be the MSSR curves they determine.

First assume that $l\neq {\rm C}'$. Then $r_{U_2}=\pm r_{U,1}$ and
$r_{V,1}=\pm r_{V,2}$, so $\phi(r_{U_2})=\phi(r_{U,1})$ and
$\phi(r_{V,1})=\phi(r_{V,2})$.  Hence the minimal pair in $\Mthree$
determined by $(\z,r_{U,i},r_{V,i})$ is the same for both values of
$i$, so $\chi_2=\chi_1$.  Thus $\M_l$ consists of a single
curve.

Now assume that the $(\z,r_{U,i},r_{V,i})$ are associated with class
C$'$ and that $(z,w)=(1,0)$. Then $(\z,r_{U,i},r_{V,i})=(\z_{j,+},
r_i,\e_i r_i)$ for some $r_i\in S^1_\bfc$ and $\e_i\in\{\pm 1\}$. A
straightforward computation yields $r_{V,i}\,\Bar{\z_{j,+}}\,
\Bar{r_{V,i}} = \e_i\frac{1}{\sqrt{2}}\ (1- (r_i)^2j).$ Thus by
Corollary \ref{cor:MSRequalp3}, a necessary condition to have
$\chi_2=\chi_1$ is
\be\label{scb6}
1-(r_2)^2j =\pm (1-r_1^2j),
\ee

\noi But $(r_i)^2j\in\spn\{j,k\}$, which holds only if $r_2=\pm
r_1$. Thus if $r_2\neq\pm r_1$, then $\chi_2\neq \chi_1$.  Conversely,
suppose that $r_2=\e r_1$, where $\e=\pm 1$. Then \eqref{scb6} holds,
$\Bar{r_{U,1}}r_{U,2}=\e$, and
\eqref{d2pdd1-final}--\eqref{rppp-final} are satisfied with $\z=1$ and
$r=\e$. Corollary \ref{cor:MSRequalp3} then implies that
$\chi_2=\chi_1$.

Thus for triples $(\z_{j,+},r_{U,i},r_{V,i})$ associated with class
C$'$, a necessary and sufficient condition for $\chi_1,\chi_2$ to
coincide (following Convention~\ref{convention5.3}) is $r_{U,2}=\pm r_{U,1}$, which is
equivalent to $R_{U,1}=R_{U,2}$ in $\So(3)$. Since $R_{U,i}\in G_D^0$,
the preceding sets up a one-to-one correspondence between $\M_{{\rm
    C}'}$ and the circle $G_D^0$:

\be\label{corr1}
\left.
\begin{array}{rcl}
(R\in G_D^0) & \longleftrightarrow& F_1(R):=\chi_{{\rm C}'}^R \\
\phantom{\int_0^1}\mbox{where}\ \ \chi_{{\rm C}'}^R
&=&
\mbox{the SR curve $[0,1]\to\Sym^+(3)$ determined by}\\
&& \mbox{the minimal pair}\ ((UR,D),(UR\phi(\z_{j,+})^{-1},
\pi_{13}\dotprod \L)).
\end{array} \right\}
\ee

\noi Note the the map $G_D^0\times [0,1]\to {\rm Sym}^+(3),
(R,t)\mapsto \chi_{{\rm C}'}^R(t),$ is continuous, hence uniformly
continuous since $G_D^0\times [0,1]$ is compact.  Thus the injective
map $F_1: G_D^0\to C([0,1],\Sym^+(3))$ (the space of continuous maps
$[0,1]\to \Sym^+(3)$) is continuous, hence a homeomorphism onto its
image (since $G_D^0$ is compact and $C([0,1],\Sym^+(3))$ is
Hausdorff), which is $\M_{{\rm C}'}$.  Therefore, in this natural
topology, $\M_{{\rm C}'}$ is homeomorphic to a circle.

While the above map $F_1$ explicitly parametrizes $\M_{{\rm C}'}$ by
the circle $G_D^0$, this parametrization is not canonical---it depends
on several non-unique choices, such as a particular matrix $U\in\Sop$
among all those that satisfy $UDU^T=X$, and our choice of
representative $\z_{j,+}$ of the double-coset
$(\Hat{\Gamma_D^0},\Hat{\Gamma_\L^0})$ double-coset (in
$\Hat{\Gamma}$) in which $\z_{j,+}$ lies.  There is a more directly
geometric parametrization of $\M_{{\rm C}'}=\M_{{\rm C}'}(X,Y)$, which
we exhibit next, by a circle in $\bfr^3$ determined by the ellipsoids
$\Sigma_X, \Sigma_Y$ to which $X,Y$ correspond.

Recall that under the Class C$'$ hypotheses ($w=0$ or $z=0$), $X$ and
$Y$ have the same, unique, axis of circular symmetry $L$ (see Remark
\ref{rem:proob}), and hence also have a common ``equatorial plane''
$L^\perp$. For $t\in [0,1]$ and $R\in G_D^0$, let $\Sigma_t^R$ be the
ellipsoid in $\bfr^3$ corresponding to $\chi_{{\rm C}'}^R(t)$; note
that $\Sigma_0^R=\Sigma_X$ and $\Sigma_1^R=\Sigma_Y$ for all $R$.

Let $Q=\log(\phi(\z_{j,+})^{-1})= \log(\phi(\z_{j,-}))$; explicitly,
\ben
Q=\frac{\pi}{2}\left[\begin{array}{rrr} 0 & 0 &
    -1\\ 0 & 0 & 0\\ 1& 0 & 0\end{array}\right]
\een

\noi (the matrix $\phi(\z_{j,-}))$ is given in Table~\ref{zetas}). For
$R\in G_D^0$ let $A(R)=UR\,Q(UR)^{-1}$.  Then
$A=\log\left(UR\phi(\z_{j,+})^{-1}(UR)^{-1}\right)$, and from
equation(s) \eqref{geod}--\eqref{defchi} we have
\bearray\nonumber
\chi_{{\rm C}'}^R(t)&=& e^{tA(R)} UR
D^{1-t}(\pi_{13}\dotprod \L)^t  (e^{tA(R)} UR)^{-1} \\
&=&U_R(t)D^{1-t}(\pi_{13}\dotprod \L)^t U_R(t)^{-1}
\label{chitht-2}
\eearray

\noi where $U_R(t)= e^{tA(R)}UR=URe^{tQ}$.  Let ${\bf e}_1=(1,0,0)^T$,
let $\e\in\{\pm 1\}$, and define $\g_R(t)= U_R(t)\e{\bf e}_1$. For all
$R\in G_D^0$, the vector $v_0=\g_R(0)$ is one of the two unit vectors
lying on the axis $L$, and $\g_R(t)$ is a unit vector lying on one of
the principal axes of $\Sigma^R_t$; equivalently, a unit eigenvector
of $\chi_{{\rm C}'}^R(t)$.  For each $R$, the inner product of $v_0$
with $\g_R(1)$ is
\ben
(UR{\bf e}_1)\dotprod (UR \phi(\z_{j,-}){\bf e}_1) =
{\bf e}_1\dotprod (\phi(\z_{j,-}){\bf e}_1) = (1,0,0)\dotprod
(0,0,1)=0.
\een

\noi Hence $\g_R(1)$ lies in the unit circle $C$ in the equatorial
plane $L^\perp$.  It is easily checked that the continuous map
$F_2:=F_{2,v_0}: G_D^0\to C$ given by $F_2(R)=\g_R(1)$ is a
bijection, hence a homeomorphism. Thus the map
$F_{v_0}:=F_1\circ (F_{2,v_0})^{-1}: C\to \M_{{\rm C}'}$
parametrizes $\M_{{\rm C}'}(X,Y)$ by the circle $C$.

One can easily check that there is at most one $t\in (0,1)$ for which
the eigenvalues of $\chi_{{\rm C}'}^R(t)$ are not all distinct.  Thus
$\g_R$ is the unique continuous map $[0,1]\to \bfr^3$ such that
$\g_R(0)=v_0$, $\g_R(1)= F_2(R)$, and $\g_R(t)$ is a unit eigenvector
of $\chi_{{\rm C}'}^R(t)$ for all $t\in [0,1]$.  Hence for each
$\chi\in \M_{{\rm C}'}(X,Y)$, there is a unique continuous map
$\tilde{\g}_\chi = \tilde{\g}_{\chi,v_0}: [0,1]\to\bfr^3$ such that
$\tilde{\g}_\chi(0)=v_0$ and $\tilde{\g}_\chi(t)$ is a unit
eigenvector of $\chi(t)$ for all $t\in [0,1]$.

This characterization shows that the parametrization $F_{v_0}:C\to
\M_{{\rm C}'}(X,Y)$ is canonical up to the choice $v_0$ of one of the
two unit vectors $L$. Given $v_0$ and a vector $w\in C$, there is a
unique $\chi=\chi_{w,v_0}\in \M_{{\rm C}'}(X,Y)$ such that the curve
$\tilde{\g}_\chi$ defined above has $\tilde{\g}_\chi(0)=v_0$ and
$\tilde{\g}_\chi(1)=w$. Moreover, $\tilde{\g}_{\chi,-v_0}(1)=
-\tilde{\g}_{\chi,v_0}(1)$, so the two parametrizations are simply
related to each other $(F_{-v_0})^{-1}=-(F_{v_0})^{-1}$.

(ii) Our proof of Theorem \ref{thm:dsr3} established that all the MSSR
curves from $X$ to $Y$ are accounted for by the curves coming from
minimal pairs in the classes listed in Table~\ref{classes}. The
first element $(UR_U,D)$ of each such pair lies in $[(U,D)]$, since
$R_U\in G_D^0$.

It remains only to establish that all MSSR curves are accounted for by
one of the (sub)cases listed in Table~\ref{tableS3} or
Table~\ref{tableS2}, and that necessary and sufficient conditions for
the curve(s) in a given class $\M_l$ to be minimal are the conditions
that can be read off from Table~\ref{tableS3} if $Y\in \esstop$, or
Table~\ref{tableS2} if $Y\in \smid$. (For example, if $Y\in \esstop$, to
read off from Table~\ref{tableS3} the conditions for the (unique)
curve $\chi_{A_1}$ in $\M_{{\rm A}_1}$ to be minimal, we simply take
the union of all the cases for which $\chi_{A_1}$ is an element of
$\M(X,Y)$, as indicated by the third column of the table.  These
conditions reduce to: $|z|\geq |w|$ and $\ell_{\rm
  id}\leq\min\{\ell_{(13)},\ell_{(12)}\}$.)

First assume that $Y\in \esstop$.  Equations
\eqref{ellcompar1}--\eqref{ellcompar3} show that no nonvacuous
subcases have been omitted in Table~\ref{tableS3}.  (For example, if
$|z|-|w|=0=\Re(\bar{z}w)$ then $\varphi=\frac{\pi}{4}=\b$, so
\eqref{ellcompar1} shows that $\ell_{\rm id}^2-\ell_{(13)}^2\neq 0$,
since, by hypothesis, $d_1\neq d_2$ and $\l_1\neq \l_3$; thus for the
two cases in Table~\ref{tableS3} in which $\ell_{\rm id}=\ell_{(13)}$,
there are no ``$|z|-|w|=0=\Re(\bar{z}w)$'' subcases. In the last case
in the table, equations \eqref{ellcompar1}--\eqref{ellcompar3} show
that no two of the three angles $\varphi, \b, \b'$ can be equal, and
hence that if $|z|=|w|$ [equivalently, $\varphi=\frac{\pi}{4}$], then
automatically $\Re(\bar{z}w)\neq 0\neq \Im(\bar{z}w)$; else we would
have $\b=\frac{\pi}{4}$ or $\b'=\frac{\pi}{4}$. Thus the hypothesis
$\Re(\bar{z}w)\neq 0\neq \Im(\bar{z}w)$ in the $|z|=|w|$ subcase of
$\ell_{{\rm id}} =\ell_{(13)} = \ell_{(12)}$ in which $|z|=|w|$ is
redundant, as asserted in the table's caption.)  Hence every MSSR
curve from $X$ to $Y$ occurs in one of the subcases listed in column 2
of Table~\ref{tableS3}.

For $l\in\{{\rm A}_1, {\rm A}_2, {\rm B}_1, {\rm B}_2, {\rm C}_1, {\rm
  C}_2\}$, let $\z_l$ be the element of $\widehat{Z}_{1,*}$ that appears
in the triples to the right of class-name $l$ in Table~\ref{classes},
and define $\ell_l\geq 0$ by
\be\label{elll}
\ell_l^2=4k\left(\cos^{-1}|f_2(\z_l)|\right)^2
+\|\log(\L_{\pi_{\z_l}}D^{-1})\|^2,
\ee

\noi where $f_2$ is as in the proof of Theorem \ref{thm:dsr3};
cf. \eqref{fibdist10}.  Then $\ell_{\rm id}=\min\{\ell_{{\rm A}_1},
\ell_{{\rm A}_2}\}$, $\ell_{(13)}=\min\{\ell_{{\rm B}_1}, \ell_{{\rm
    B}_2}\}$, and $\ell_{(12)}=\min\{\ell_{{\rm C}_1}, \ell_{{\rm
    C}_2}\}$. A set of necessary and sufficient conditions to have
$\dsr(X,Y)= \ell_{A_i}$ (respectively $\ell_{B_i}, \ell_{C_i}$) is:
(a) $\ell_{{\rm A}_i}\leq \ell_{{\rm A}_{i'}}$ (resp. $\ell_{{\rm
    B}_i}\leq \ell_{{\rm B}_{i'}}$, $\ell_{{\rm C}_i}\leq \ell_{{\rm
    C}_{i'}}$), where $\{i,i'\}=\{1,2\}$, and (b) $ \min\{\ell_{\rm
  id},\ell_{(13)}, \ell_{(12)}\}=\ell_{\rm id}$ (resp. $\ell_{(13)},
\ell_{(12)}$).

For $l={\rm A}_1$ and $l={\rm A}_2$ the contributions to $\ell_l$ from
the term involving $D$ on the right-hand side of \eqref{elll} are
identical, so $\ell_{{\rm A}_1}\leq \ell_{{\rm A}_2}$ if and only if
$|f_2(1)|\geq |f_2(j)|$. Thus, $\ell_{{\rm A}_1}\leq \ell_{{\rm A}_2}
\iff |z|\geq |w|$. Similarly, $\ell_{{\rm B}_1}\leq \ell_{{\rm B}_2}
\iff |z+w|\geq |z-w| \iff \Re(\bar{z}w) \geq 0$, and $\ell_{{\rm
    C}_1}\leq \ell_{{\rm C}_2} \iff |z-iw|\geq |z+iw| \iff
\Im(\bar{z}w) \geq 0$.  Note that $\max\{|z|,|w|\}$,
$\max\{|z+w|,|z-w|\}$, and $\max\{|z+iw|,|z-iw|\}$ are all strictly
positive. Hence, if $\dsr(X,Y)= \ell_l$, then $f_2(\z_l)\neq 0$,
$\Hat{f_2(\z_l)}$ is defined, and the curve $\chi_l$ associated with
the data listed in Table~\ref{classes} is defined.

Thus, for $Y\in \esstop$, all MSSR curves from $X$ to $Y$ are accounted for
in Table~\ref{tableS3}, and in each case listed in the table, a curve
$\chi_l$ is minimal (where $l\in\{{\rm A}_1, {\rm A}_2, {\rm B}_1,
{\rm B}_2, {\rm C}_1, {\rm C}_2\}$) if and only if the conditions
indicated in the table are satisfied.

Now assume that $Y\in \smid$. Analogously to the case $Y\in \esstop$,
define $\ell_l\geq 0$ by
\be\label{elll2}
\ell_l^2=4k\left(\cos^{-1}|f_4(\z_l)|\right)^2
+\|\log(\L_{\pi_{\z_l}}D^{-1})\|^2,
\ee

\noi where $l\in\{{\rm A}_1, {\rm A}_2, {\rm B}_1, {\rm B}_2, {\rm
  C}_1, {\rm C}_2\}$, $\z_l$ is the element of $\widehat{Z}_{1,*}$ that
appears in the triples to the right of class-name $l$ in
Table~\ref{classes}, and $f_4$ is as in the proof of Theorem
\ref{thm:dsr3} (again cf. \eqref{fibdist10}). Then $\ell_{\rm
  id}=\min\{\ell_{{\rm A}_1'}, \ell_{{\rm A}_2'}\}.$ A set of
necessary and sufficient conditions to have $\dsr(X,Y)= \ell_{A_i'}$
is: (a)$'$ $\ell_{{\rm A}'_i}\leq \ell_{{\rm A}'_{i'}}$ , where
$\{i,i'\}=\{1,2\}$, and (b)$'$ $\ell_{\rm id}\leq
\ell_{(13)}\}$. Noting that that exactly one of the classes B$'$ and
C$'$ is defined for a given $Y$, a necessary and sufficient conditions
to have $\dsr(X,Y)= \ell_{{\rm B}'}$ (respectively $\ell_{{\rm C}'}$)
is $\ell_{(13)}\leq \ell_{\rm id}$.

The same reasoning used in the case $Y\in \esstop$ shows now that
$\ell_{{\rm A}_1'}\leq \ell_{{\rm A}_2'}$ if and only if $|z|\geq
|w|$, and that if $\dsr(X,Y)= \ell_l$, then the curve-class $\M_l$
associated with the data listed in Table~\ref{classes} is defined.  It
follows that for $Y\in \smid$, all MSSR curves from $X$ to $Y$ are
accounted for in Table~\ref{tableS2}, and in each case listed in the
table, the curve(s) $\chi$ in the class $\M_l$ (where $l\in\{{\rm
  A}_1', {\rm A}_2', {\rm B}', {\rm C}'\}$) is/are minimal if and only
if the conditions indicated in the table are satisfied (modulo
Convention~\ref{convention5.3} in the case of class C$'$).

 (iii) Since we have now established that $\M(X,Y)$ consists of
precisely those curves listed in Table~\ref{tableS3} for the subcase
corresponding to the given data $((U,D), \linebreak (V,\L))$, it
suffices to show that if $\chi_1,\chi_2$ are MSSR curves in distinct
classes $l_1, l_2 \in \{{\rm A}_1, {\rm A}_2, {\rm B}_1, {\rm B}_2,$
${\rm C}_1, {\rm C}_2\}$, then $\chi_1\neq \chi_2$.

Given such $l_1, l_2$, for $i\in\{1,2\}$ let $(\z_i, r_{U,i},r_{V,i})$
be a triple from Table~\ref{classes} correspond to class $l_i$.  Since
$r_{V,1}=1=r_{V,2}$ for all such triples in all classes corresponding
to $Y\in \esstop$, we may rewrite \eqref{a2eqa1} as
$\Bar{\z_1}\ \z_2=\pm \Bar{r_{U,1}}\ r_{U,2}$.  But
$\Bar{r_{U,1}}\ r_{U,2} \in \bfc$, so by Corollary
\ref{cor:MSRequalp3} a necessary condition to have $\chi_1=\chi_2$ is
\be\label{neccond}
\Bar{\z_1}\ \z_2 \in \bfc.
\ee

\noi We compute the following:
\bearray
\label{12eq1j}
(\z_1,\z_2) = (1, j) &\implies& \Bar{\z_1}\ \z_2 =j, \\
\nonumber
(\z_1,\z_2) = (\z_{j,+}, \z_{j,-}) &\implies& \Bar{\z_1}\ \z_2 =-j, \\
\nonumber
(\z_1,\z_2) = (\z_{k,+}, \z_{k,-}) &\implies& \Bar{\z_1}\ \z_2 =-k, \\
\nonumber
(\z_1,\z_2) \in \{1,j\}\times \{\z_{j,\pm}\} &\implies&
\Bar{\z_1}\ \z_2 \in \left\{\frac{\pm 1 \pm j}{\sqrt{2}}\right\}, \\
\nonumber
(\z_1,\z_2) \in \{1,j\}\times \{\z_{k,\pm}\} &\implies&
\Bar{\z_1}\ \z_2 \in \left\{\frac{1 \pm k}{\sqrt{2}} \,,
\frac{\pm i -j}{\sqrt{2}}\right\}, \\
\nonumber
(\z_1,\z_2) \in \{\z_{j,\pm}\}\times \{\z_{k,\pm}\} &\implies&
\Bar{\z_1}\ \z_2 \in \left\{\frac{1 \pm i \pm j \pm k}{2}\right\}.
\eearray

\noi Since $\Bar{\z}_1\z_2\notin \bfc$ in every case, it follows that
$\chi_1\neq \chi_2$.

 (iv) Analogously to part (iii), since we have now established that
$\M(X,Y)$ consists of precisely those curves listed in
Table~\ref{tableS2} for the subcase corresponding to the given data
$((U,D),(V,\L))$, and that $\M_{{\rm C'}}$ (when defined) contains
infinitely many curves, it suffices to show that if $\chi_1,\chi_2$
are MSSR curves in distinct classes $l_1, l_2 \in \{{\rm A}_1', {\rm
  A}_2', {\rm B}', {\rm C}'\},$ then $\chi_1\neq \chi_2$. Since
exactly one of the classes $\M_{{\rm B'}}\,,\M_{{\rm C'}}$ is
nonempty, we do not need to consider the case $\{l_1,l_2\}=\{{\rm B}',
{\rm C}'\}$. Because of Convention~\ref{convention5.3}, we
also do not need to consider the case $\{l_1,l_2\}=\{{\rm A_2}', {\rm
  C}'\}$. Thus we need only consider the case-pairs $(l_1,l_2)\in\{
({\rm A}_1', {\rm A}_2'), ({\rm A}_1', {\rm B}'), ({\rm A}_1', {\rm
  C}'),$ $({\rm A}_2', {\rm B}')\}$.

Given such $(l_1, l_2)$, for $i\in\{1,2\}$ let
$(\z_i,r_{U,i},r_{V,i})$ again be a data-triple from
Table~\ref{classes} correspond to class $l_i$.  By Corollary
\ref{cor:MSRequalp3}, to show $\chi_1\neq\chi_2$ it suffices to show
that \eqref{a2eqa1} is not satisfied. If $l_1={\rm A}_1'$ then
$\z_1=1$, and \eqref{a2eqa1} cannot be satisfied unless $\z_2\in
\bfc$, which does not hold since $l_2\neq {\rm A}_1'$.  For the case
$(l_1,l_2)=({\rm A}_2',{\rm B}')$, if \eqref{a2eqa1} were satisfied we
would have $\Bar{\z_{j,+}}=\xi_1 j\xi_2$ for some $\xi_1,\xi_2\in
\bfc$, an impossibility since $\xi_1 j\xi_2 =\xi_1\Bar{\xi_2}j \in
\bfc j =\spn\{j,k\}.$ \qedns
%

   \subsection{Algorithm for computing MSSR curves
for $p=3$ in the nontrivial cases}\label{algo}

Let $X,Y\in\Sym^+(3)$ be as in Theorem \ref{allMSR}.  Starting with eigen-decompositions
$(U,D)$ of $X$, $(V,\L)$ of $Y$, an algorithm to compute all the MSSR curve(s)
from $X$ to $Y$ is as follows.  This algorithm applies {\em only} when $p=3$, and only to the
nontrivial cases.  

\begin{itemize}

\item[Step 1.] If $U^{-1}V$ is not an involution, proceed to Step 2.  If
$U^{-1}V$ is an involution, find an even sign-change matrix $I_\mbs$  for which $d_{SO}(U^{-1}VI_\mbs,I)< d_{SO}(U^{-1}V,I)$.
The pair $(VI_\sigma, \L)$ is still a
pre-image of $Y$ since the action of sign-change matrices on diagonal
matrices is trivial.  Replace $V$ by $VI_\sigma$, renamed to
$V$. Proceed to Step 2.

\item[Step 2.] Find $\th\in [0,\pi)$ and a unit vector $\tilde{a}\in
  \bfr^3$ such that $U^{-1}V=R_{\th,\tilde{a}}$. There is a unique such
  $\th$ and, if $\th\neq 0$, a unique such $\tilde{a}$; writing
  $R=U^TV$ these can be computed using
\begin{align}
\label{rodinv1}
\th &=\cos^{-1}\frac{\tr(R)-1}{2}\ ,\\
\tilde{a}
&= \left\{ \begin{array}{ll}
\frac{1}{2\sin\th}(R_{32}-R_{23}, R_{13}-R_{31}, R_{21}-R_{12})^T&
\mbox{if}\ \theta\neq 0,\\
{\bf 0} & \mbox{if}\  \th = 0.
\end{array}
\right.
\label{rodinv2}
\end{align}

\noi (Equations \eqref{rodinv1}--\eqref{rodinv2} are consequences of
the well-known ``Rodrigues formula''.)

\item[Step 3.] Define $z,w\in \bfc$ by $z+wj=s(U^{-1}V)$, where
$s:\So(3)_{<\pi}\to S^3$ is the map given by \eqref{halfangle}.
The ``minimal classes", i.e. the curve-classes containing an element of $\M(X,Y)$,
as well as the cardinality of $\M(X,Y)$, can then
be read off from Table~\ref{tableS3} if $Y\in\esstop$, or Table~\ref{tableS2} if $Y\in\smid$.  (If $Y\in\esstop$, first compute the numbers $\ell_{\rm id}, \ell_{(13)},$ and $\ell_{(12)}$, defined
in \eqref{fibdist4-2to3.piid}--\eqref{fibdist4-2to3.pi12}, to use Table~\ref{tableS3}.)  The
appropriate line of Table~\ref{classes} then gives the pairs
$(r_U,r_V)\in S^1_\bfc\times S^1_\bfc$ for each curve-class. The remaining
steps of this algorithm are applied to each minimal class.

\indenum Note that for any class other than
C$'$, all the $(r_U,r_V)$ pairs in Table~\ref{classes} determine the same scaling-rotation
curve, so just choose one pair from this line of the table. If the
data are in class C$'$, there will be one MSSR curve for each $r\in
S^1_\bfc$, but the $\pm$ sign in the table can be ignored (treated as
+), since the sign does not affect the image under $\phi$.

\item[Step 4.] For the chosen $(r_U,r_V)$ in each minimal class (there will only be one in each class
except for Class C$'$), compute the rotations $R_U=\phi(r_U),
R_V=\phi(r_V)$ from the unit complex numbers $r_U, r_V$ using the
general formula

\be
\phi(e^{ti}) = \left[ \begin{array}{lcr} 1 & 0 & 0\\
0& \cos 2t & -\sin 2t \\ 0 & \sin 2t
    & \cos 2t\end{array}\right] \ \ \mbox{for}\ \ t\in\bfr.
\ee

\noi (Note that if we identify the $x_2x_3$ plane with $\bfc$ via
$(x_2,x_3)\leftrightarrow x_2+x_3i$, then for $\xi\in S^1_\bfc$ the
lower right $2\times 2$ submatrix of $\phi(\xi)$ corresponds simply to
multiplication by $\xi^2$.  In Case B$'$ this conveniently ``undoes''
the square roots in Table~\ref{classes}; for example, if
$r_U=\pm(\hat{z}\hat{w})^{1/2}$, then $\phi(r_U)$ is the rotation
about the $x_1$ axis that corresponds to multiplying $x_2+x_3i$ by
$\hat{z}\hat{w}$.)

\item[Step 5.] Read off the value of $\phi(\z)$ from Table~\ref{zetas}. Then
plug this and the values of $(R_U,R_V)$ computed in Step 4 into
\eqref{minpr}, yielding (for each of these pairs) the endpoints of a
geodesic from $\E_X$ to $\E_Y$ whose projection to ${\rm Sym}^+(3)$ is
an MSSR curve.

\item[Step 6.] For each of the endpoint-pairs computed in Step 5, writing
the endpoints as $(U_1,D)\in\E_X, (V_1,\L_1)\in\E_Y$, set
$A=\log(U_1^{-1}V_1), L=\log(D^{-1}\L_1)$.  Then use formulas
\eqref{geod} and \eqref{defchi} (with $U_1$ playing the role of $U$ in
these formulas) to compute the formula for the corresponding MSSR
curve $\chi:[0,1]\to \sympp$.

\end{itemize}

\begin{remark} {\rm  For the case in which each of $X$ and $Y$ has exactly two distinct eigenvalues,
this algorithm for computing closed-form expressions for MSSR curves
in the $p=3$ nontrivial cases replaces the numerical algorithm
in \cite{JSG2015scarot} described therein after Theorem 4.3. }
\end{remark}

\setcounter{equation}{0}

\section{Unique and non-unique cases of MSSR curves in $\symp(3)$}\label{sec:appendix-MSSR3}
In the section, we give examples, graphical illustrations,
and further discussion of unique and non-unique cases of MSSR curves in $\symp(3)$. As shown in
Theorem \ref{thm:dsr3}, we can divide our analysis into four possibilities for the strata in which $X$ and $Y$, the endpoints of the MSSR curve, lie.

\begin{enumerate}
\item[(i)]  $X,Y \in \Stop$.
\item[(ii)]  $X \in \Sc_{\rm mid}$, $Y \in \Stop$.
\item[(iii)]  $X,Y \in \Sc_{\rm mid}$.
\item[(iv)] $X \in \Sbot$.
\end{enumerate}

The case (i) in which both $X$ and $Y$ have three distinct eigenvalues is discussed in Section~\ref{sec:caseI}.
For cases (ii) and (iii), graphical illustrations  and further discussion of all classes of MSSR curves from $X$ to $Y$ are provided in Section \ref{sec:mid-top} for case (ii), and in Section \ref{sec:case_mid} for case (iii).
It is easy to see that in case (iv) there is a unique MSSR curve from $X$ to $Y$ for any $Y \in \symp(3)$.

\subsection{The case in which both $X$ and $Y$ have three distinct eigenvalues}\label{sec:caseI}
Let $X,Y \in \Stop := \Sc_{[\J_{\rm top}]}$, and let $(U,D) \in \Ec_X$, $(V,
\Lambda) \in \Ec_Y$. 
By Proposition \ref{distprop},
$((U,D), (VP_{g}^{-1}, \pi_{g} \dotprod \Lambda))$ is a minimal pair if $g$ is in the set
\begin{equation}\label{eq:dist_ing}
\Nsf_{((U,D), (V,\Lambda))} = \argmin_{g \in \tilde{S}_3^+}
                     \left\{  k d_{SO(3)}^2 (U, VP_g^{-1}) + d_\Dc^2 (D, \pi_g\dotprod\Lambda)
                      \right\}.
\end{equation}
Depending on $(U,D), (V,\Lambda)$, any $g \in \tilde{S}_3^+$ can provide a minimal pair.
Let $n_{(X,Y)} :=  |\Nsf_{((U,D), (V,\Lambda))}|$, which is insensitive to particular choices of $((U,D), (V,\Lambda))$. If  $ n_{(X,Y)}  = 1$, then the MSSR curve from $X$ to $Y$ is unique; more generally, there are exactly $n_{(X,Y)}$ MSSR curves.


 Given a particular pair $((U,D), (V,\Lambda))$, the 24 elements of $\tilde{S}_3^+$ label the corresponding scaling-rotation curves, which are candidates for the MSSR curves.
A strategy to characterize all unique and non-unique cases of MSSR curves is to divide the minimization problem into smaller subproblems. To this end,
we classify these subproblems  according to the six possibilities for $\mbox{proj}_2(g) = \pi_g \in S_3$.
Recall that we write $\pi_{\rm id}$ for the identity permutation, and, for
distinct $a,b \in\{1,2,3\}$, we write $\pi_{ab}$ for the transposition $(ab)$, the permutation
that just interchanges $a$ and $b$. The six elements of $S_3$ are
\begin{align}
  \pi_{\rm id}, &  \nonumber\\
  \pi_{12}, & \quad  \pi_{23},\nonumber\\
  \pi_{123}:=\pi_{23}\pi_{12}, & \quad \pi_{132}:=\pi_{12}\pi_{23}, \nonumber\\
  \pi_{13} &(=  \pi_{12} \pi_{23}\pi_{12} =  \pi_{23}\pi_{12}\pi_{23}). \label{eq:allpermutations}
\end{align}
For each permutation $\pi_\star \in S_3$, we can find $X,Y$ such that there are unique or non-unique $g$'s projecting to $\pi_\star$ that give the smallest distance.
The subproblem for each $\pi_\star$ has many subcases. For example, among the four $g$'s with $\pi_g = \pi_{\rm id}$, there are 8 subcases of (possibly) unique or non-unique MSSR curves, determined by the value of $U^{-1}V$. Instead of analyzing all subcases, we focus on the classification by the values of $\pi_g$, which provides interesting information concerning the  corresponding MSSR curves.

\begin{figure}[tb!]
  \centering
  \vspace{-1.0in}
\includegraphics[trim={3cm 0cm 2cm 1cm},width=0.9\textwidth]{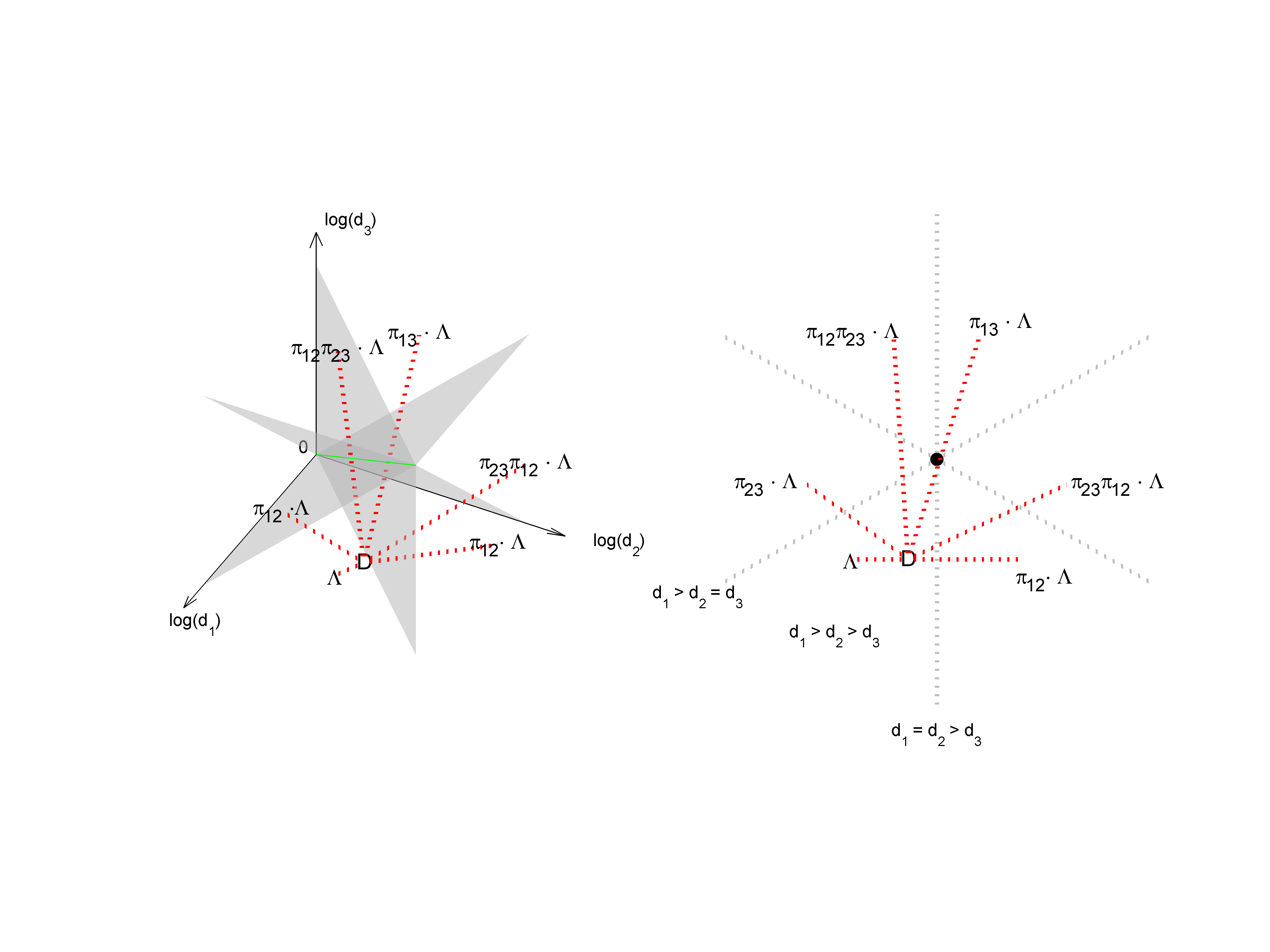}
  \vspace{-1.0in}
  \caption{ The space $\diag^+(3)$ with the line and planes representing the its stratification (left) and a cross section (right) of the stratified $\diag^+(3)$; see Fig. \ref{fig:Symp3stratification}.
  For $D, \Lambda$ in the same connected component of $\Dc_{\J_{\rm top}} \subset \diag^+(3)$, the distance $d_\Dc(D,\Lambda)$ is represented by the ``length'' of the red dotted line. Also shown are $\pi_g \cdot \Lambda$ and $d_\Dc(D, \pi_g \cdot \Lambda)$ for all possible $\pi_g$. Here, we chose $D = \mbox{diag}(8,6,3)$, $\Lambda = \mbox{diag}(15, 8, 6)$.
  \label{fig:distinct_case}}
\end{figure}

For this purpose, we assume that $D$ and $\Lambda$ are in the same connected component of $\Dc_{\J_{\rm top}}$ (e.g., $d_1>d_2>d_3$ and $\lambda_1>\lambda_2>\lambda_3$), so that the minimum of $d_\Dc (D, \pi_{g}\dotprod\Lambda)$ is achieved by choosing $g$ to satisfy $\pi_g = \pi_{\rm id}$. Moreover, if $D$ and $\Lambda$ satisfy $d_1 > d_2 > d_3$, $\lambda_1 > \lambda_2 > \lambda_3$, then
$$
d_\Dc (D, \pi_{\rm id}\dotprod\Lambda) \le \left\{ \begin{matrix}
                                         d_\Dc (D, \pi_{12}\dotprod\Lambda) \\
                                           \mbox{ or } \\
                                         d_\Dc (D, \pi_{23}\dotprod\Lambda) \\
                                        \end{matrix} \right\}
                            \le \left\{ \begin{matrix}
                                         d_\Dc (D, \pi_{123}\dotprod\Lambda) \\
                                           \mbox{ or } \\
                                         d_\Dc (D, \pi_{132}\dotprod\Lambda) \\
                                        \end{matrix} \right\}
                            \le d_\Dc (D, \pi_{13}\dotprod\Lambda).
$$
In order for a $g$ such that $\pi_g \neq \pi_{\rm id}$ to give a minimal pair, the corresponding ``rotation distance'' (the first term of (\ref{eq:dist_ing})) needs to be sufficiently smaller than the four rotation distances associated with the identity permutation. A typical example of the distances $d_\Dc (D, \pi_{g}\dotprod\Lambda)$ is illustrated in Fig.~\ref{fig:distinct_case}.



In this example, because $\pi_g = \pi_{\rm id}$ gives the smallest $d_\Dc (D, \pi_{g}\dotprod\Lambda)$, for sufficiently small $k$ the minimizer $g$ of (\ref{eq:dist_ing}) satisfies $\pi_g = \pi_{\rm id}$, and the ellipsoids corresponding to the MSSR curves from $t = 0$ to $1$ are always tri-axial. For fixed $k$, other choices of $g$ can provide a minimal pair, depending on the values of $(U,D)$, $(V, \Lambda)$. Below, we list
the shape-classification changes of the MSSR curve $\chi_{g}$ corresponding to the particular $g \in \Nsf_{((U,D), (V,\Lambda))} \subset \tilde{S}^+_3$. A $3 \times 3$ SPD matrix with eigenvalues $a \ge b \ge c$ has the shape of a sphere if $a = b = c$, oblate spheroid (or oblate, for short) if $a = b > c$, prolate spheroid (or prolate) if $a> b = c$, and tri-axial ellipsoid (or tri-axial) if $a > b > c$.
We assume below that $D$ and $\Lambda$ have been chosen to lie in the same connected component of $\Dc_{\J_{\rm top}}$.

\begin{enumerate}
\item For all $g$ such that $\pi_g = \pi_{\rm id}$, for all $t \in [0,1]$, the ellipsoid corresponding to $\chi_g(t) $ is always tri-axial.
\item For all $g$ such that $\pi_g = \pi_{\rm 12}$, the shape-classification changes of the MSSR curves $\chi_g(t)$ from $t= 0$ to $t = 1$ are (tri-axial $\to$ oblate $\to$ tri-axial).
\item For all $g$ such that $\pi_g = \pi_{\rm 23}$, the shape-classification changes  are (tri-axial $\to$ prolate $\to$ tri-axial).
\item For all $g$ such that $\pi_g = \pi_{\rm 123}$, the shape-classification changes  are (tri-axial $\to$ oblate $\to$ tri-axial $\to$ prolate $\to$ tri-axial).
\item For all $g$ such that $\pi_g = \pi_{\rm 132}$, the shape-classification changes  are (tri-axial $\to$ prolate $\to$ tri-axial $\to$ oblate $\to$ tri-axial).
\item For all $g$ such that $\pi_g = \pi_{\rm 13}$, the shape-classification changes  are either \\
    (tri-ax. $\to$ oblate $\to$ tri-ax. $\to$ prolate $\to$ tri-ax. $\to$ oblate $\to$ tri-ax.),\\
    (tri-ax. $\to$ prolate  $\to$ tri-ax. $\to$ oblate$\to$ tri-ax. $\to$ prolate  $\to$ tri-ax.), \\
    or (tri-axial $\to$ sphere $\to$ tri-axial).
\end{enumerate}

\begin{figure}[tb!]
  \centering
\includegraphics[trim={4cm 1cm 3cm 1cm},width=0.8\textwidth]{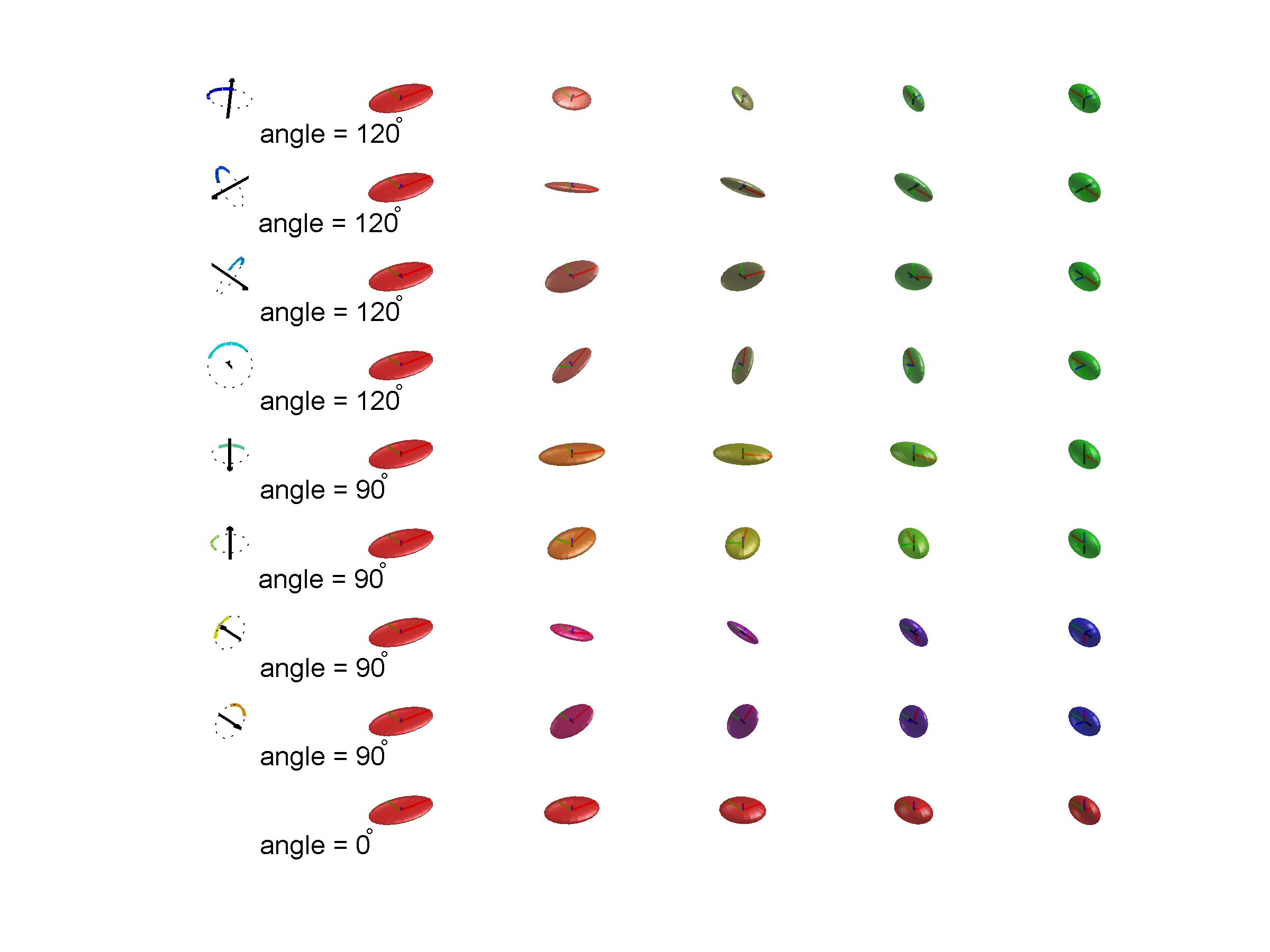}
  \caption{An example for one case of non-unique MSSR curves with $X,Y \in  \Stop$. There are nine MSSR curves.
 The rotation parameter $A$ of each MSSR curve is depicted as the axis-angle figure in the left-most panels. 
  The top four MSSR curves correspond to the minimal pairs provided by $g$ with $\pi_g = \pi_{\rm id}$. The next two MSSR curves correspond to two choices of $g$ with $\pi_g = \pi_{12}$. The next two MSSR curves correspond to two choices of $g$ with $\pi_g = \pi_{23}$. The last MSSR curve corresponds to $g$ with $\pi_g = \pi_{123}$.\label{fig:distinct_worst case}}
\end{figure}

We close this section by providing the worst-case example of non-uniqueness we have found (for $X$ and $Y$ both having three distinct eigenvalues), in which there are 9 MSSR curves. We choose $X = \mbox{diag}(e^a,e^b,e^c)$, where $c = 1$, $b = c + \frac{7\sqrt{7}}{6\sqrt{8}}\pi$, $a = b + \sqrt{\frac{7}{72}}\pi$, and $Y = R \mbox{diag}(e^x,e^y,e^z) R^T$ where $ \phi^{-1}(R) = \frac{\pm 1}{2}(1+i+j+k)$ and
$z = 0$, $y= \frac{1}{3\sqrt{14}}\pi$, $x  = y + \sqrt{\frac{7}{72}}\pi$.
The nine MSSR curves are illustrated in Fig.~\ref{fig:distinct_worst case}. (See Section~\ref{sec:5.1.1.relhso3} for the definition of $\phi^{-1}(R)$.)

\subsection{The case in which $X \in \Smid$ and $Y \in \Stop$}\label{sec:mid-top}
For this special case, $X$ has just two distinct eigenvalues, while $Y$ has three. We parameterize $X$ with $a,b \in \Real$ (arbitrary, not size-ordered), and any $U \in \SO(3)$, and parameterize $Y$ with $c>d>f>0$ and a unit quaternion $q = z + wj \in S_{\mathbf{H}}^3$, where $z, w \in \Complex$ and, $|z|^2 + |w|^2 = 1$, as follows:
\begin{align} \label{eq:case3par}
X = U
     \begin{pmatrix}
       e^a & 0 & 0 \\
       0 & e^b & 0 \\
       0 & 0 & e^b \\
     \end{pmatrix} U^T, \quad Y = U \phi(q)
                             \begin{pmatrix}
       e^c & 0 & 0 \\
       0 & e^d & 0 \\
       0 & 0 & e^f \\
     \end{pmatrix} (U \phi(q) )^T
\end{align}
where $\phi: S_{\mathbf{H}}^3 \to \SO(3)$ is the natural two-to-one Lie-group homomorphism defined in Section \ref{sec:5.1.1.relhso3}. Without loss of generality, we assume $\mbox{Re}(q) = \mbox{Re}(z) >0$ so that $\phi(q)$ is not an involution.

There are six different cases of MSSR curves as summarized in Table \ref{classes}, named ${\rm A}_1, {\rm A}_2, {\rm B}_1,{\rm B}_2, {\rm C}_1, {\rm C}_2$. To give representative examples of these cases, we partially rewrite the conditions in Table \ref{tableS3}.

Recall that $\varphi = \cos^{-1} ( \max\{|z|,|w|\}) \in [0, \frac{\pi}{4}]$.
For each $\varphi \in (0,\pi/4)$ there are exactly two cases: $\cos(\varphi) = |z| > |w|$ and $|z| < |w| = \cos(\varphi)$. If $\varphi = 0$, then $|z| = 1, |w| = 0$. If $\varphi = \frac{\pi}{4}$, then $|z| = |w| = 1/\sqrt{2}$. We define
\begin{equation} \label{eq:alpha}
\alpha = \left\{
           \begin{array}{ll}
             \cos^{-1} \left(    \frac{|{\rm Re}(\bar{z} w)|}{|\bar{z} w|}\right), & \hbox{if $\varphi > 0$;} \\
             \pi/2, & \hbox{if $\varphi = 0$,}
           \end{array}
         \right.
 \in [0, \frac{\pi}{2}],
\end{equation}
so that the parameters $\beta$ and $\beta'$ appearing in  Theorem~\ref{thm:dsr3} and Table \ref{tableS3} will be represented by
\begin{align}
   \beta  &= \frac{1}{2} \cos^{-1} ( \sin(2\varphi) \cos\alpha ) \in [0, \frac{\pi}{4}],  \label{eq:beta and betap1}\\
   \beta' &= \frac{1}{2} \cos^{-1} ( \sin(2\varphi) \sin\alpha ) \in [0, \frac{\pi}{4}].  \label{eq:beta and betap2}
\end{align}
Note that
\begin{align}
\varphi = 0 \ \Longleftrightarrow \  |z| = 1 \ \Longleftrightarrow \ \beta = \beta' = \frac{\pi}{4}
\end{align}
and that
\begin{align}
\mbox{Re}(\bar{z} w) = 0 \ \Longleftrightarrow \  \cos\alpha = 0  \ \Longleftrightarrow \  \alpha = \frac{\pi}{2}, \\
\mbox{Im}(\bar{z} w) = 0 \ \Longleftrightarrow \  \cos\alpha = 1  \ \Longleftrightarrow \  \alpha = 0.
\end{align}
For any $\alpha < \pi/2$,  $\mbox{Re}(\bar{z} w)$ can have either sign. For $\alpha > 0$,  $\mbox{Im}(\bar{z} w)$ can have either sign.

\begin{figure}[tb!]
  \centering
  \vspace{0.3cm}
\includegraphics[trim={3cm 1cm 3cm 1cm},width= .6\textwidth]{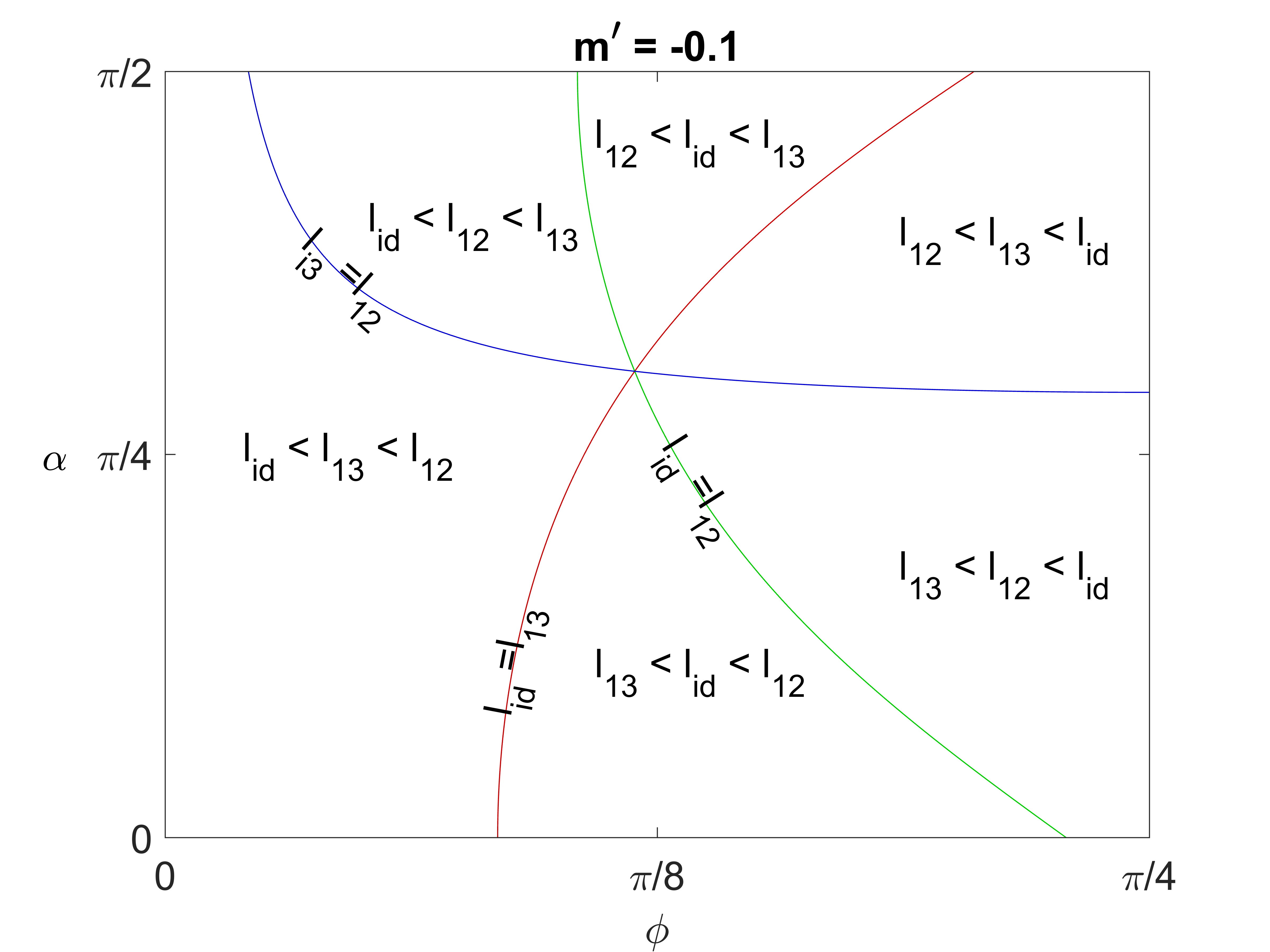}
  \caption{The ordering of $\ell_{\rm id}, \ell_{13}, \ell_{12}$ corresponding to the value  of $(\varphi,\alpha)\in [0,\pi/4] \times [0, \pi/2]$. Note that $m' = -0.1$ is fixed. \label{fig:case3graphic1}
}
\end{figure}

\begin{figure}[tb!]
  \centering
  \vspace{0.3cm}
\includegraphics[trim={3cm 1cm 3cm 1cm},width= .6\textwidth]{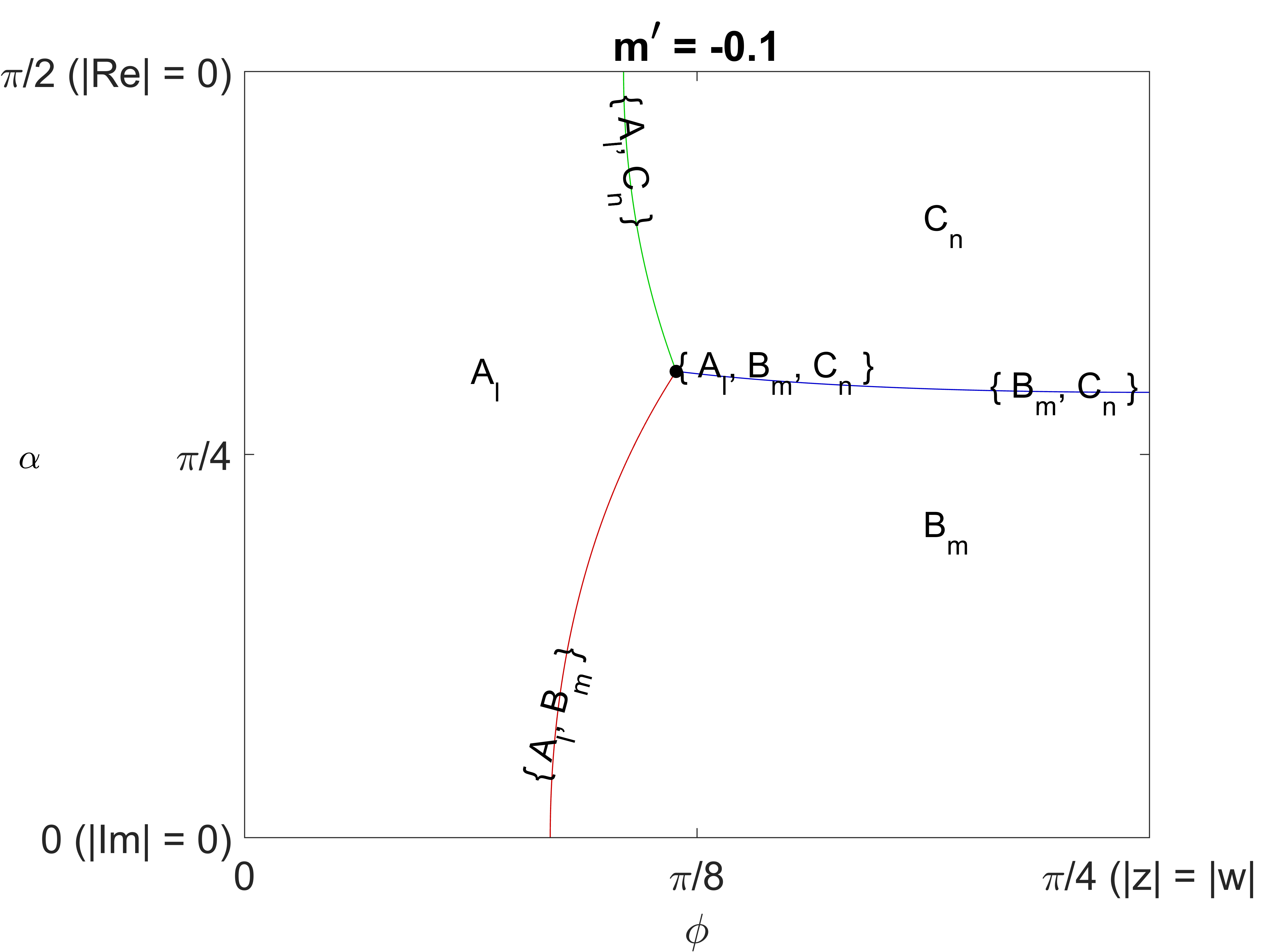}
  \caption{Cases of unique and non-unique MSSR curves in $\symp(3)$ when $X$ has just two distinct  eigenvalues, and $Y$ has three distinct eigenvalues. This figure uses the same parametrization of Fig.~\ref{fig:case3graphic1}, with $m' = -0.1$. The case ${\rm A}_l$ stands for either ${\rm A}_1$ (if $|z| > |w|$), ${\rm A}_2$ (if $|z| < |w|$) or $\{{\rm A}_1,{\rm A}_2\}$ (if $|z| = |w|$). The latter case ($|z| = |w|$) can only occur if $\varphi = \pi/4$. The case ${\rm B}_m$ stands for either ${\rm B}_1$ (if $\mbox{Re}(\bar{z}w) > 0$), ${\rm B}_2$ (if $\mbox{Re}(\bar{z}w) < 0$), or $\{{\rm B}_1, {\rm B}_2\}$ (if $\mbox{Re}(\bar{z}w) = 0$). The latter case ($\mbox{Re}(\bar{z}w) = 0$) can only occur if $\alpha = \pi/2$. Similarly, the case ${\rm C}_n$ stands for either ${\rm C}_1$ (if $\mbox{Im}(\bar{z}w) > 0$), ${\rm C}_2$ (if $\mbox{Im}(\bar{z}w) < 0$), or $\{{\rm C}_1, {\rm C}_2\}$ (if $\mbox{Im}(\bar{z}w) = 0$). The latter case ($\mbox{Im}(\bar{z}w) = 0$) can only occur if $\alpha = 0$. The multi-element cases such as $\{{\rm A}_l, {\rm B}_m\}$ can be understood similarly. The values of $l$ and $m$ are determined by the signs of $|z|-|w|$ and $\mbox{Re}(\bar{z}w)$, and if there are two values for $l$ and one value for $m$, then the case is $\{{\rm A}_1,{\rm A}_2,{\rm B}_m\}$, a case with three MSSR curves.
\label{fig:case3graphic2}
}
\end{figure}

\begin{figure}[tb!]
  \centering
  \vspace{0.3cm}
\includegraphics[trim={3cm 1cm 3cm 1cm},width= .8\textwidth]{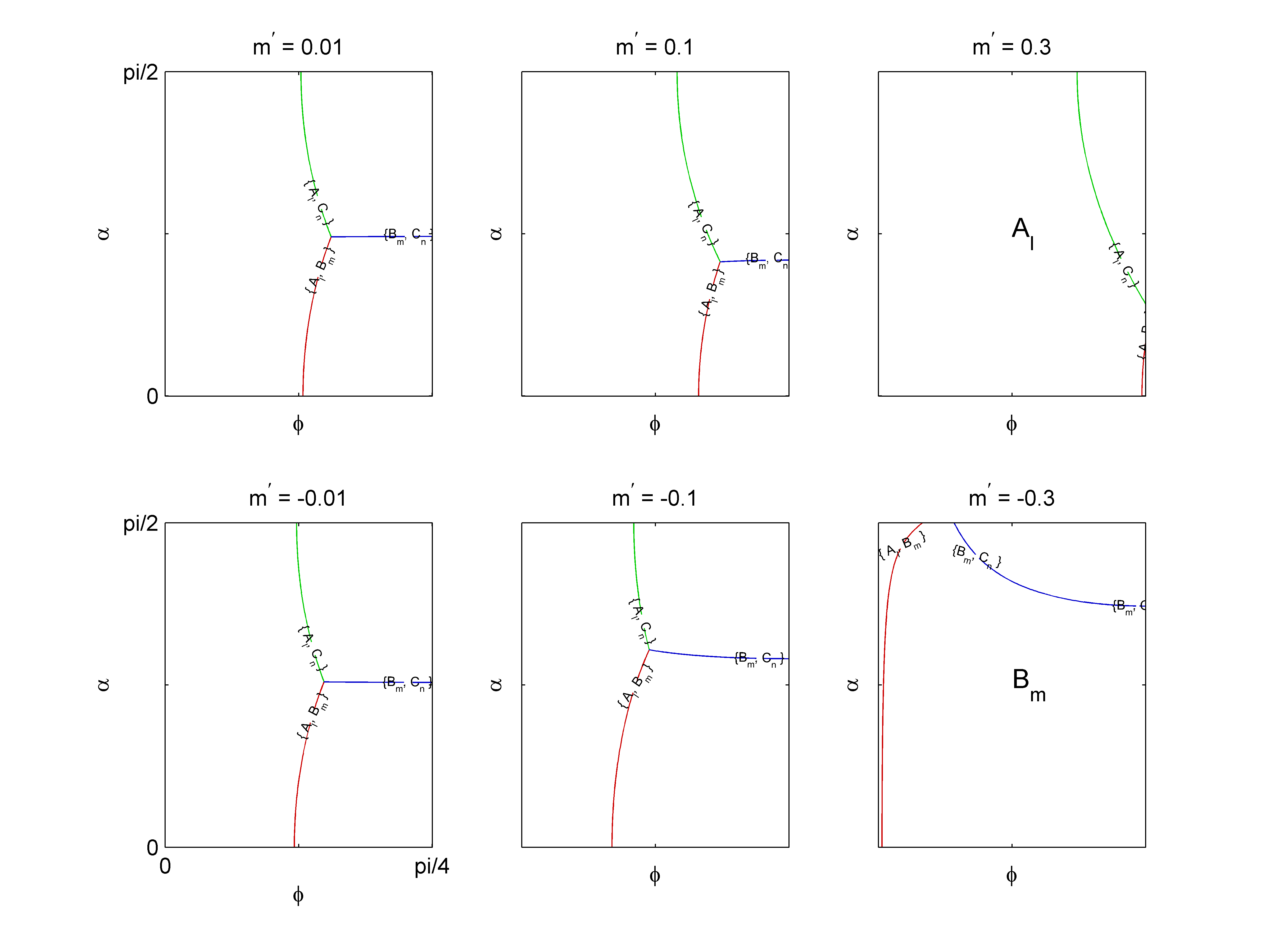}
  \caption{Unique and non-unique MSSR curves in $\symp(3)$ when $X$ has just two distinct eigenvalues, and $Y$ has three distinct eigenvalues.
\label{fig:case3graphic3}
}
\end{figure}

We also make use of the following parameters, concerning the eigenvalues of $X$ and $Y$, scaled by $k$ (where $k>0$ is as in (\ref{gm})):
\begin{align}
m_1 = \frac{2 (a- b) (c - d)}{4k},\quad m_2 = \frac{2 (a- b) (d - f)}{4k}.
\end{align}
Each of $m_1$ and $m_2$ can be either positive or negative, but they must both have the same sign. To simplify the analysis, we assume
\begin{align}
 m_1 = m_2 := m'.
\end{align}
Then we have
\begin{align}
\ell^2_{\rm id} &< \ell^2_{13}  \ \Longleftrightarrow \   \varphi^2 - \beta^2 - 2m' <0, \label{eq:4.32}\\
\ell^2_{\rm id} &< \ell^2_{12}  \ \Longleftrightarrow \   \varphi^2 - (\beta')^2 -m'<0, \label{eq:4.33}\\
\ell^2_{13} &< \ell^2_{12}  \ \Longleftrightarrow \   \beta^2 - (\beta')^2 + m' <0.\label{eq:4.34}
\end{align}
By (\ref{eq:beta and betap1}), (\ref{eq:4.32}) is equivalent to
$$
\alpha > \cos^{-1} \left(  \frac{\cos(2\sqrt{(\varphi^2 - 2 m')_+})}{\sin(2\varphi)}\right),\quad \mbox{ if } \varphi > 0,
$$
and to $m' > - \frac{\pi^2}{32} $ if $\varphi = 0$.
By (\ref{eq:beta and betap2}), (\ref{eq:4.33}) is equivalent to
$$
\alpha < \sin^{-1} \left(  \frac{\cos(2\sqrt{(\varphi^2 -  m')_+})}{\sin(2\varphi)}\right),\quad \mbox{ if } \varphi > 0,
$$
and to $m' > - \frac{\pi^2}{16} $ if $\varphi = 0$. We have not found a similarly simple inequality equivalent to
(\ref{eq:4.34}). Figure~\ref{fig:case3graphic1}, generated numerically, indicates the regions of $(\varphi, \alpha)$ corresponding to different size-orders of $\ell_{\rm id}, \ell_{13}, \ell_{12}$, for the fixed value $m' = -0.1$.
The seven cases of unique and non-unique MSSR curves summarized in Table \ref{tableS3} are graphically represented in Figs.~\ref{fig:case3graphic2} and~\ref{fig:case3graphic3}. As $|m'|$ increases, either ${\rm A}_l$ or ${\rm B}_m$ becomes the only case of MSSR curves.

In the following we provide several examples of the unique and non-unique cases of MSSR curves for the case $X \in \Smid$ and $Y \in \Stop$. We first discuss the shape classification changes of the scaling-rotation curves $\chi_{{\rm A}_l}(t)$, $\chi_{{\rm B}_m}(t)$ and $\chi_{{\rm C}_n}(t)$ $(l,m,n = 1,2)$. These depend on the sign of $m'$.
\begin{enumerate}
\item If $m' >0 $ (that is, $X$ is prolate, and $Y$ is tri-axial), then the shape-classification changes of
     $\chi_{{\rm A}_l}(t)$ are (prolate $\to$ tri-axial); for
     $\chi_{{\rm B}_m}(t)$ they are (prolate $\to$ tri-axial $\to$ oblate $\to$ tri-axial $\to$ prolate $\to$ tri-axial);  for $\chi_{{\rm C}_n}(t)$ they are (prolate $\to$ tri-axial $\to$ oblate $\to$ tri-axial).
\item If $m' <0 $ (that is, $X$ is oblate, and $Y$ is tri-axial), then the shape-classification changes of
     $\chi_{{\rm A}_l}(t)$ are (oblate $\to$ tri-axial $\to$ prolate $\to$ tri-axial $\to$ oblate $\to$ tri-axial); for
     $\chi_{{\rm B}_m}(t)$ they are (oblate $\to$ tri-axial); for
     $\chi_{{\rm C}_n}(t)$ they are (oblate $\to$ tri-axial $\to$ prolate $\to$ tri-axial).
\end{enumerate}
Two scaling-rotation curves ${\rm A}_1$ and ${\rm A}_2$ (or ${\rm B}_1$ and ${\rm B}_2$, ${\rm C}_1$ and ${\rm C}_2$) share the same scaling parameters, and also share the same rotation axis, but with different orientations (one clockwise, the other counterclockwise) and possibly different angles. The two angles differ by $\pi$.

\begin{figure}[tbp!]
  \centering
\vskip -0.1in
\includegraphics[width=1\textwidth]{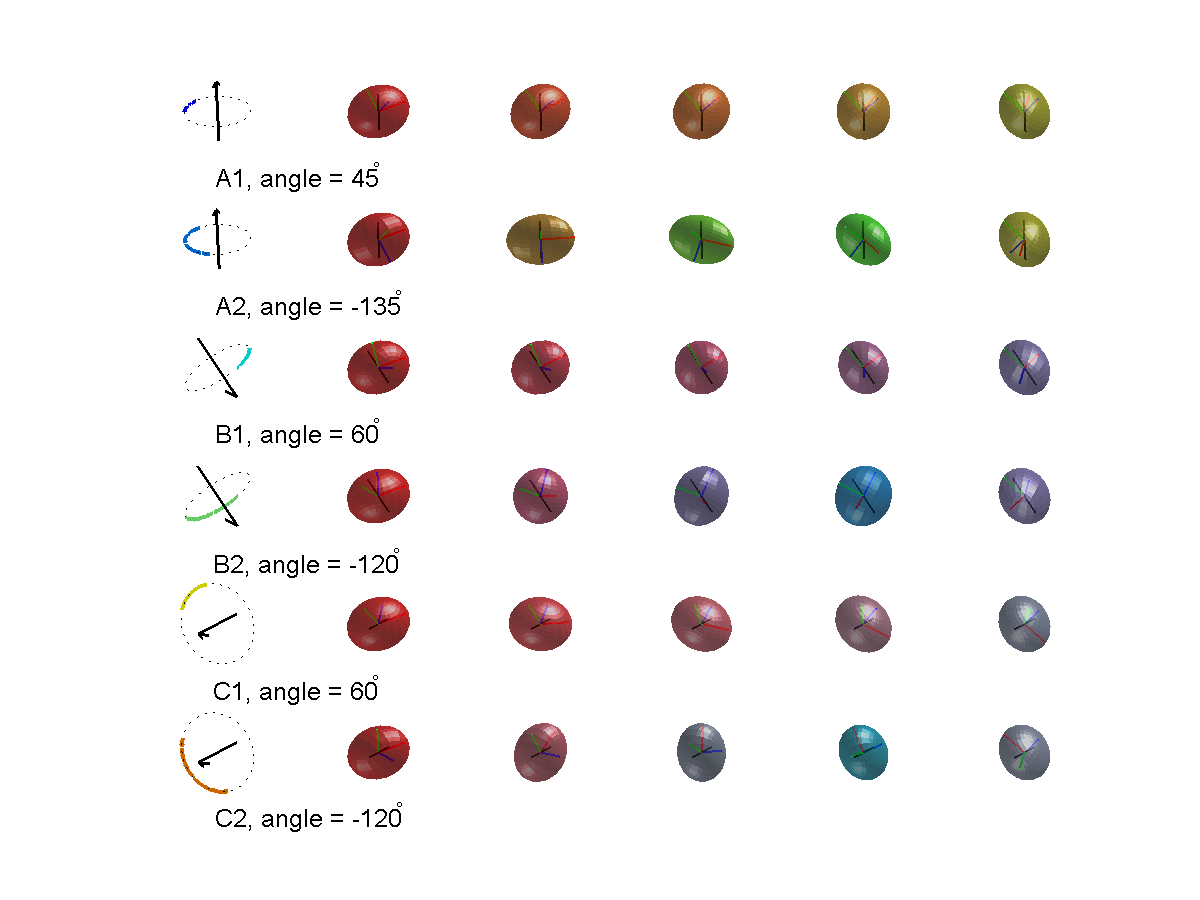}
\vskip -0.1in
\includegraphics[width=1\textwidth]{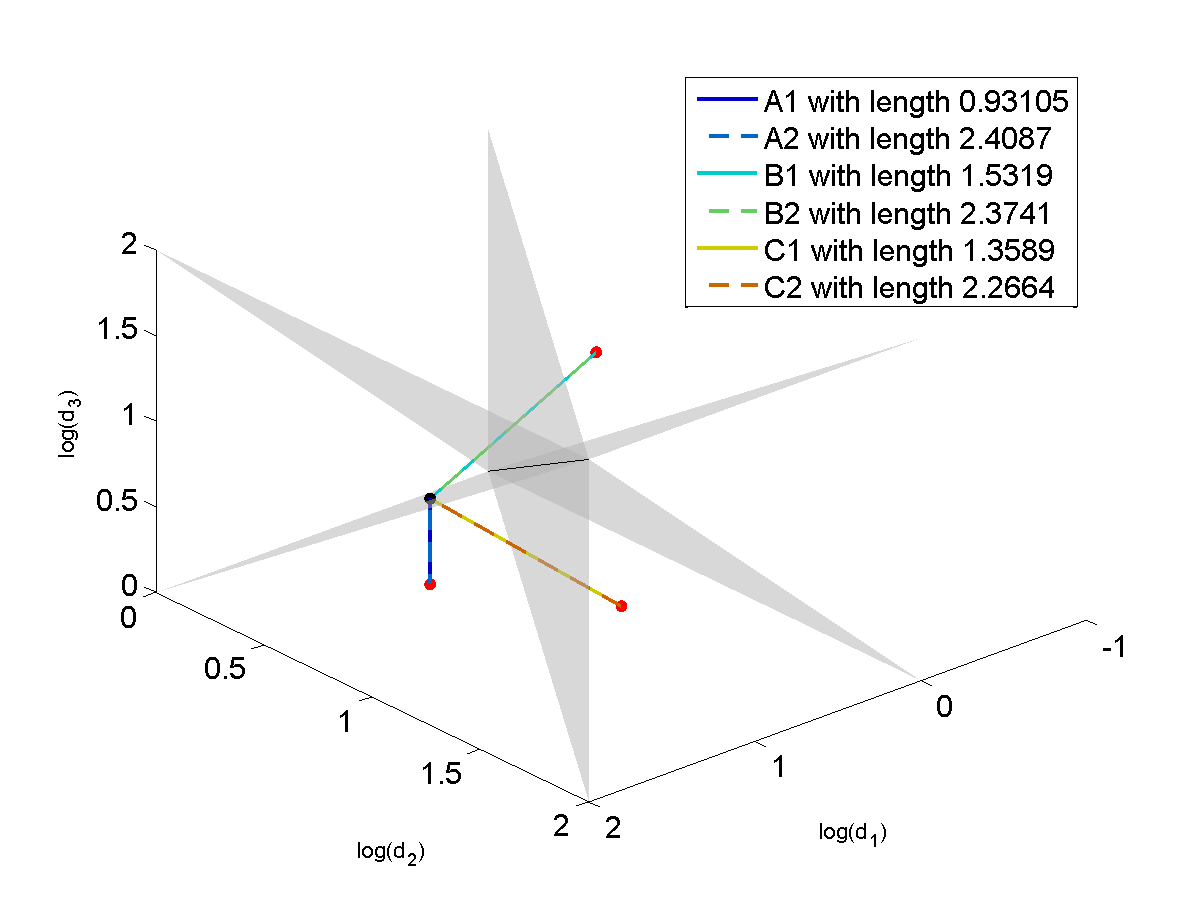}
  \caption{An example for the case $\{{\rm A}_1\}$.
   (The scaling-rotation curve corresponding to ${\rm A}_1$ is the unique MSSR curve.)
   \label{fig:A1} Each of the cases ${\rm A}_1$ - ${\rm C}_2$ represents the corresponding scaling-rotation curve defined in Table \ref{classes}, whose length can be found in the legend.
Note that $m' > 0$ in this example. The eigenvalue paths corresponding to ${\rm A}_1$ and ${\rm A}_2$ depart from a single connected component of $\Dc_{\J_1}$ (represented by the shaded open half-plane containing the black dot) and reach $\Dc_{\rm top}$.  Accordingly, the  shape-classification changes are  (prolate $\to$ tri-axial).
}
\end{figure}

\begin{figure}[tbp!]
  \centering
\includegraphics[width=1\textwidth]{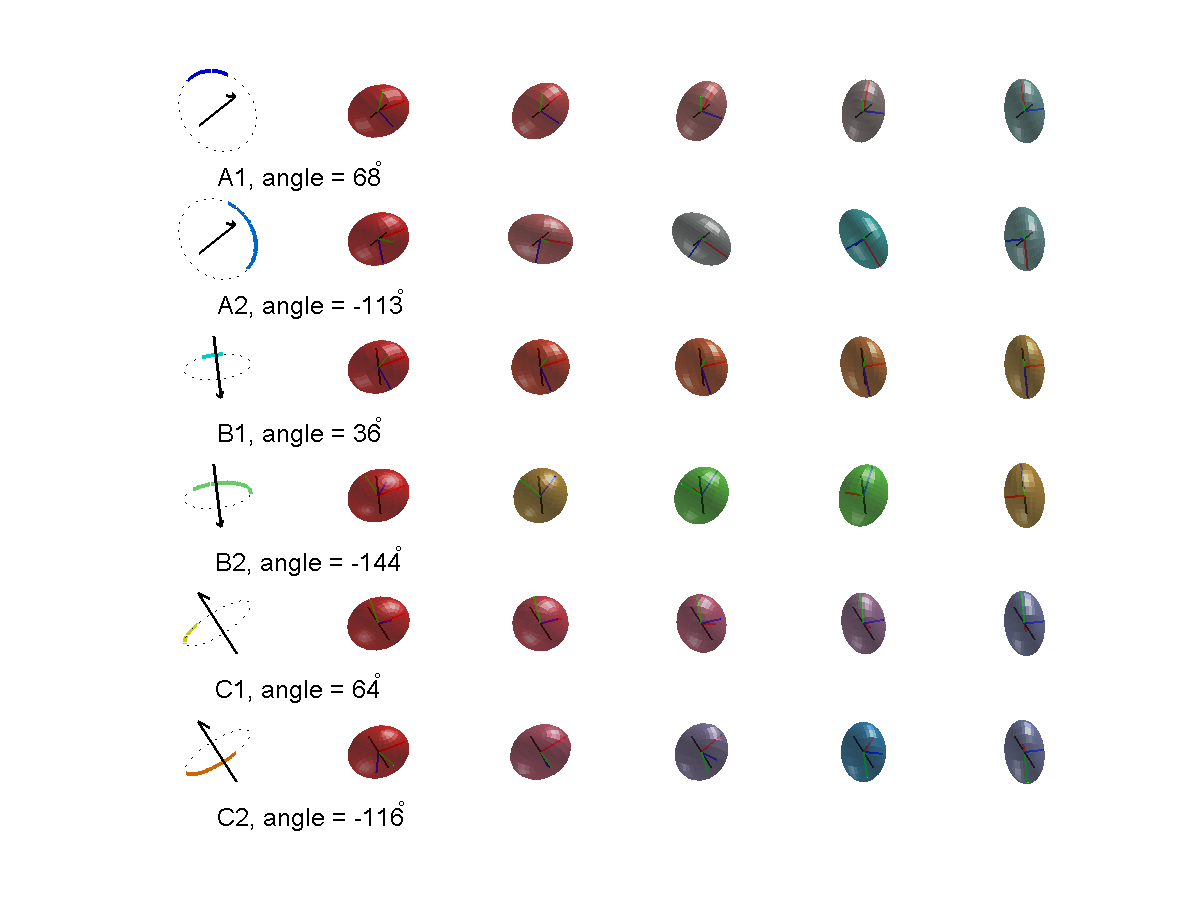}\vspace{-0.0in}
\includegraphics[width=1\textwidth]{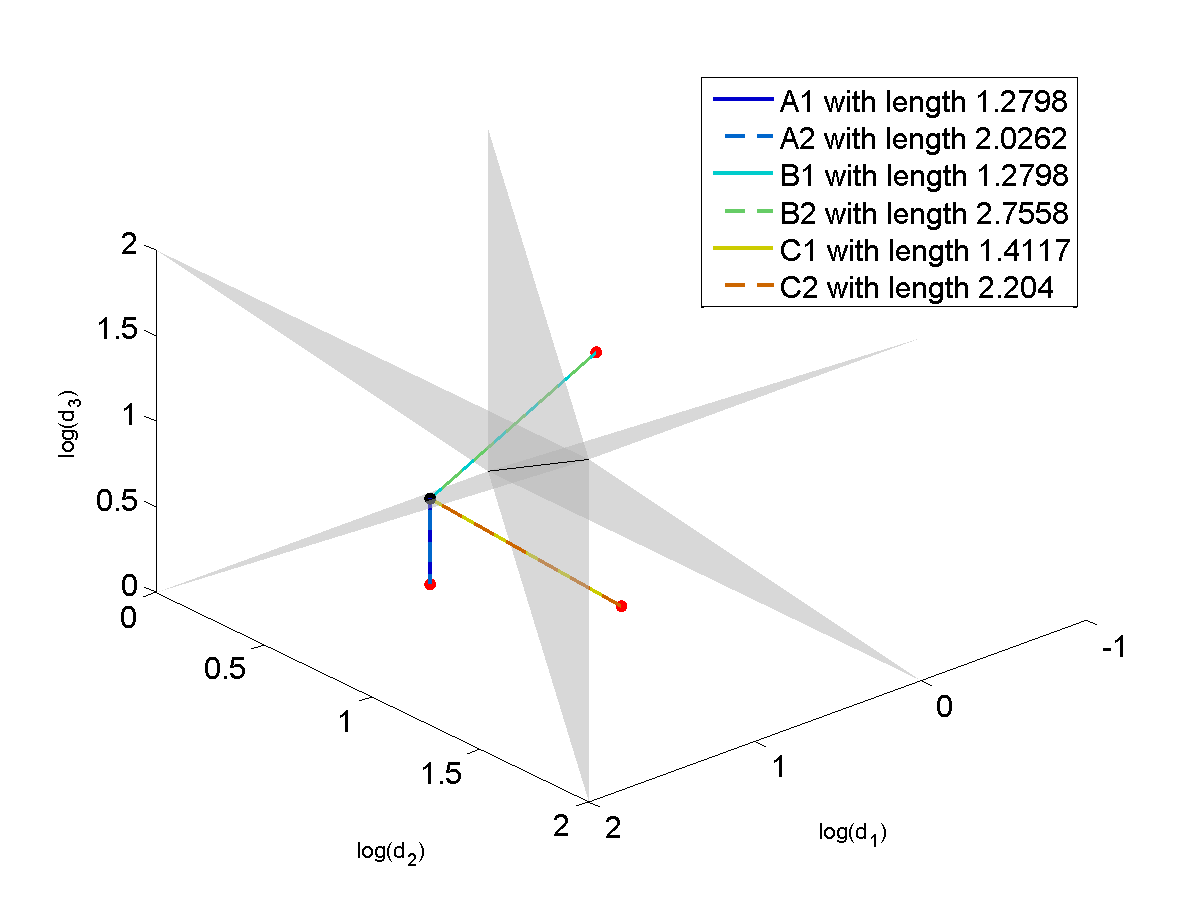}
  \caption{An example for the case $\{{\rm A}_1, {\rm B}_1\}$. \label{fig:A1B1case3}
}
\end{figure}

\begin{figure}[tbp!]
  \centering
\includegraphics[width=1\textwidth]{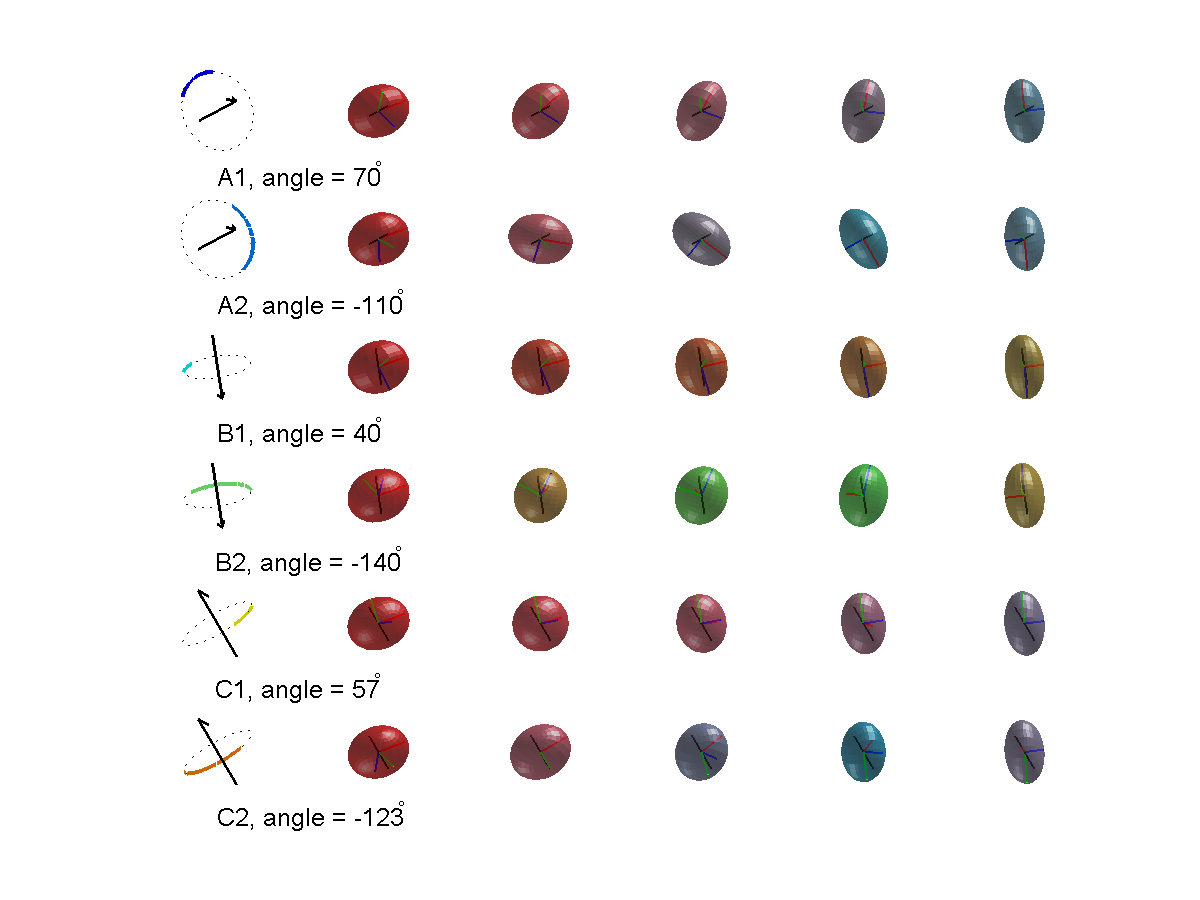}\vspace{-0.0in}
\includegraphics[width=1\textwidth]{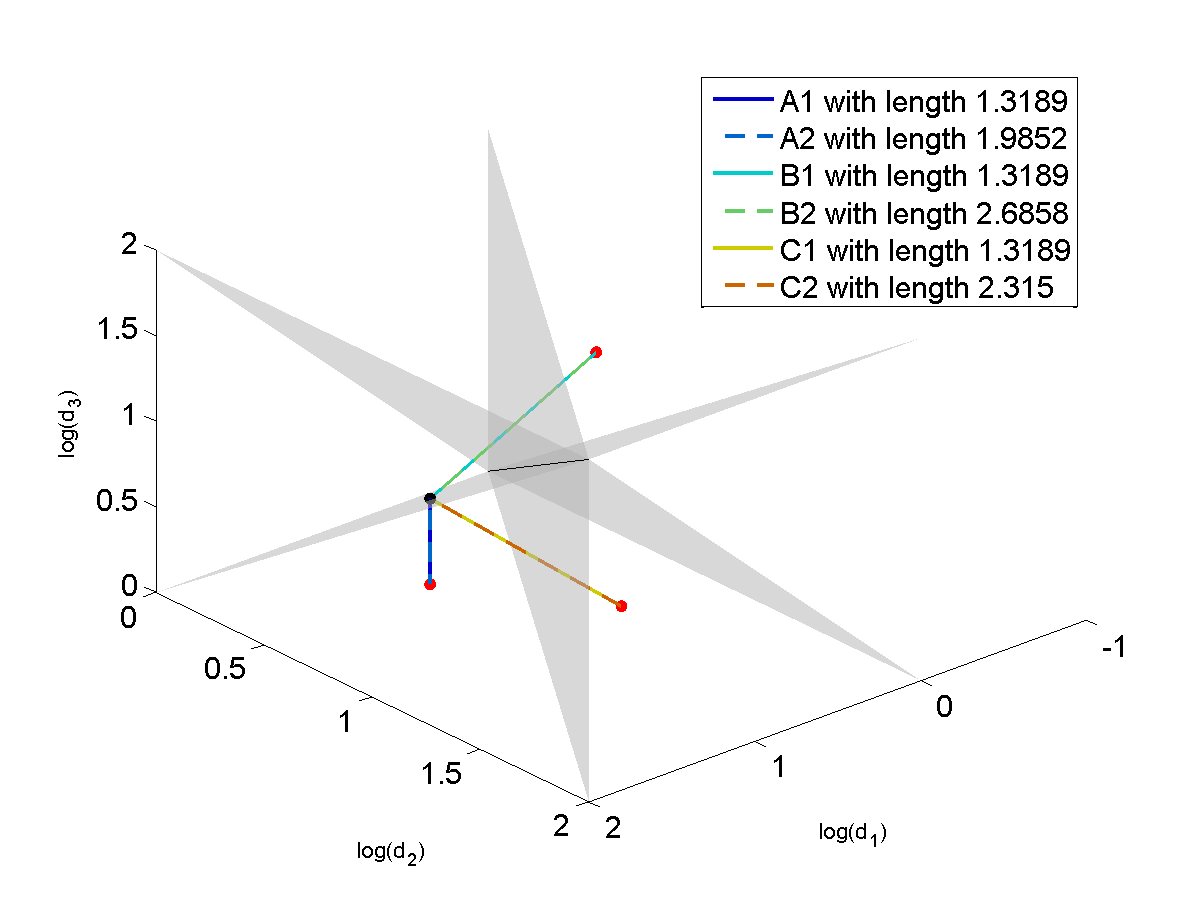}
  \caption{An example for the case $\{{\rm A}_1,{\rm B}_1,{\rm C}_1\}$. \label{fig:ABCcase3}
}
\end{figure}

In an attempt to visually illustrate examples in Fig.~\ref{fig:A1} to Fig.~\ref{fig:ABCcase3}, a scaling-rotation curve $\chi := \chi_{U,D, A, L}$ in the 6-dimensional space $\symp(3)$ is depicted by both a sequence of ellipsoids (representing the discretized scaling-rotation curve $\chi$) and the combination of the rotation parameter $A$ and changes of eigenvalues $D\exp(Lt)$.
The rotation parameters $A$ of $\chi$ are depicted as the axis-angle figure in the top left-most panels of figures below.
 The bottom panel of each figure depicts a logarithmic projection to $\diag(3)$ of the curve $D\exp(Lt) \in \diag^+(3)$ from $D$ (black dot) to $D\exp(L)$ (red dots). (See Fig.~\ref{fig:Symp3stratification} for the definition of shaded planes.)
Since the rotational degrees of freedom have been projected out in the bottom panel, the relative lengths of the straight-line-segments do not accurately reflect the relative lengths of the curves.

We collect representative visual examples in which there are a unique MSSR curve (Fig.~\ref{fig:A1}, the case $\{{\rm A}_1\}$),  exactly two MSSR curves (Fig.~\ref{fig:A1B1case3}, the case $\{{\rm A}_1,{\rm B}_1\}$), and  exactly three MSSR curves (Fig.~\ref{fig:ABCcase3}, the case $\{{\rm A}_1,{\rm B}_1,{\rm C}_1\}$). More situations, including a four-MSSR-curve case, are possible.

%
%
%
%

%
%
%
%
%
%
%
%
%
%

%


\subsection{The case in which both $X$ and $Y$ have exactly two distinct eigenvalues}\label{sec:case_mid}
For this special case, we parameterize $X$ and $Y$ with $a,b,c,d \in \Real$ (not ordered), and a unit quaternion $q = z + w j \in S_{\mathbf{H}}^3$, where $z, w \in \Complex$, $|z|^2 + |w|^2 = 1$ as
\begin{align} \label{eq:case2par}
X = U
     \begin{pmatrix}
       e^a & 0 & 0 \\
       0 & e^b & 0 \\
       0 & 0 & e^b \\
     \end{pmatrix} U^T, \quad Y = U \phi(q)
                             \begin{pmatrix}
       e^c & 0 & 0 \\
       0 & e^d & 0 \\
       0 & 0 & e^d \\
     \end{pmatrix} (U \phi(q) )^T
\end{align}
where $U \in \SO(3)$. Without loss of generality, we assume $\mbox{Re}(q) = \mbox{Re}(z) >0$ so that $\phi(q)$ is not an involution.

There are four different cases of MSSR curves arising from the parametrization of (\ref{eq:case2par}), as summarized in Table \ref{classes}. These cases are denoted ${\rm A}_1', {\rm A}_2', {\rm B}' $ and ${\rm C}'$. Our goal here is to further investigate the seven subcases of $\Mc(X,Y)$ in Table \ref{tableS2}, by partially rewriting the conditions in Table \ref{tableS2}. Note that
$0 < | z| \le 1$, $0 \le \cos^{-1} |z| < \frac{\pi}{2}$, $|w| = \sqrt{1 - |z|^2} = \cos( \frac{\pi}{2} - \cos^{-1} |z| )$, and
\begin{align}
\varphi  &= \cos^{-1}( \max\{|z|, |w| \} ) = \min\{ \cos^{-1} |z|, \cos^{-1}(\sqrt{1 - |z|^2}) \} \nonumber \\
    &= \left\{
               \begin{array}{ll}
                         \cos^{-1} |z|, & 0 \le \cos^{-1} |z| \le \frac{\pi}{4};\\
                \frac{\pi}{2}- \cos^{-1} |z|, &  \frac{\pi}{4} \le \cos^{-1} |z| < \frac{\pi}{2}.
                 \end{array}
                \right.  \label{eq:varphi_alt}
\end{align}

\begin{figure}[tb!]
  \centering
\includegraphics[width=0.7\textwidth]{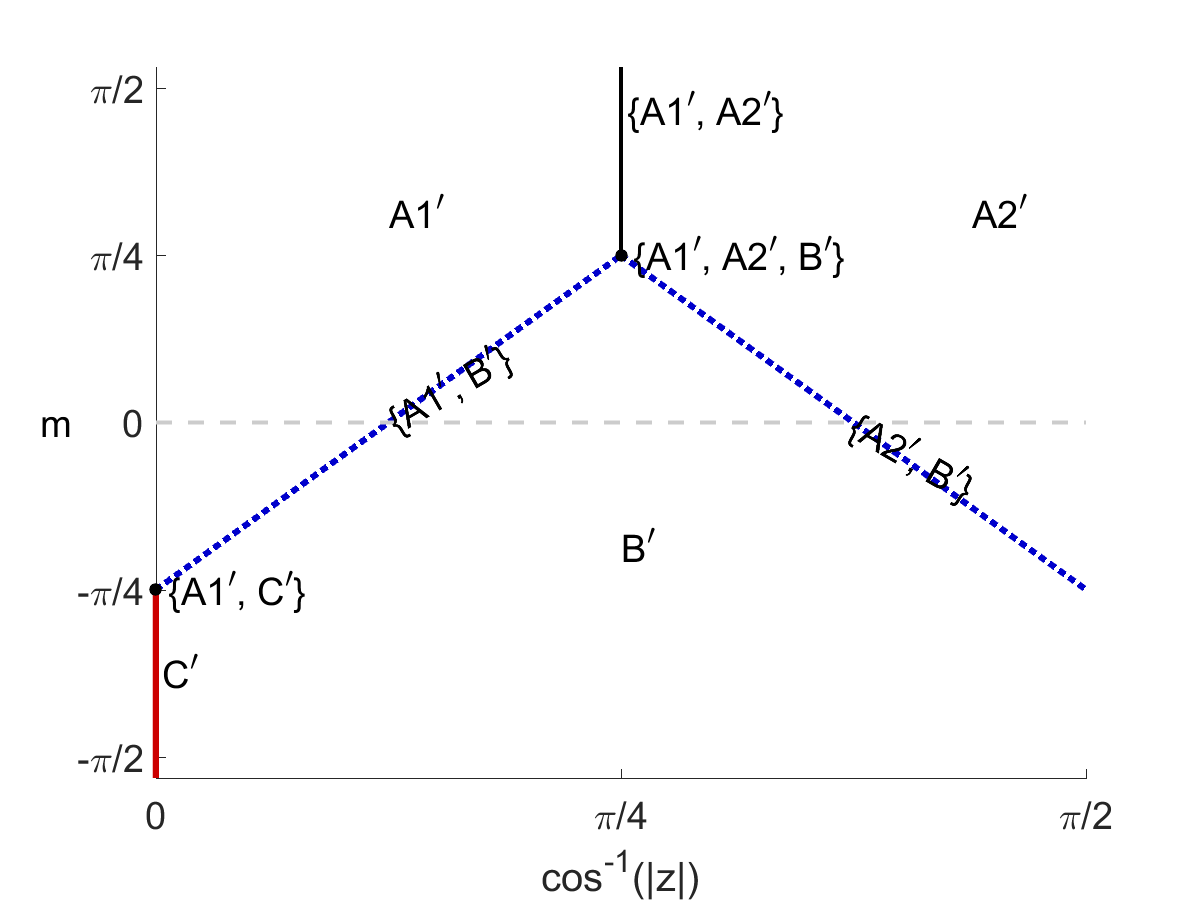}
  \caption{Unique and non-unique MSSR curves in $\symp(3)$ when both $X$ and $Y$ have exactly two distinct eigenvalues: Schematic illustration for the nine subcases. The horizontal line $m = 0$ is excluded. \label{fig:22}}
\end{figure}

\noindent Let
$m = \frac{2(a-b)(c-d)}{k \pi}.$ Unlike in the $\symp(2)$ case (cf. Section~\ref{sec:appendix-MSSR2}), this $m$ can be either positive or negative, but not zero.  If $X$ and $Y$ are both prolates (or both oblates), then $m > 0$. If $X$ is an oblate and $Y$ is a prolate (or vice versa), then $m < 0$. By Theorem~\ref{thm:dsr3},
\begin{align}
\ell_{\rm id} < \ell_{13} \ \Longleftrightarrow \  m > 2 ( \varphi - \frac{\pi}{8}).
\end{align}
Moreover,
\begin{align}
  |z| > |w| \ & \Longleftrightarrow \  \frac{1}{\sqrt{2}} < |z| \le 1 \  \Longleftrightarrow \ 0 \le \cos^{-1}|z| < \frac{\pi}{4}\\
  |z| < |w| \ & \Longleftrightarrow \  0 <  |z| < \frac{1}{\sqrt{2}} \  \Longleftrightarrow \ \frac{\pi}{4} < \cos^{-1}|z| < \frac{\pi}{2}  \label{eq:z and w}
\end{align}
Combining (\ref{eq:varphi_alt})-(\ref{eq:z and w}), the conditions for the nine subcases in Table \ref{tableS2} are represented by $m \neq 0$ and $\cos^{-1}|z| \in [0, \frac{\pi}{2})$, as shown in Fig.~\ref{fig:22}. {These subcases can be sub-divided by the sign of $m$. }

In the following we take a few representative examples of the nine subcases. Each example of MSSR curve $\chi$ is accompanied by its shape-classification changes.
Note that the seven subcases not involving ${\rm C}'$
  resemble the seven subcases of $p = 2$; see Section~\ref{sec:appendix-MSSR2}. The case ${\rm C}'$ is only defined for $|z| = 1$ (or $\cos^{-1}|z| = 0$), in which situations the cases ${\rm A}_2'$, ${\rm B}_1$ are not defined.
%

\begin{figure}[tbp]
  \centering
\includegraphics[width=1\textwidth]{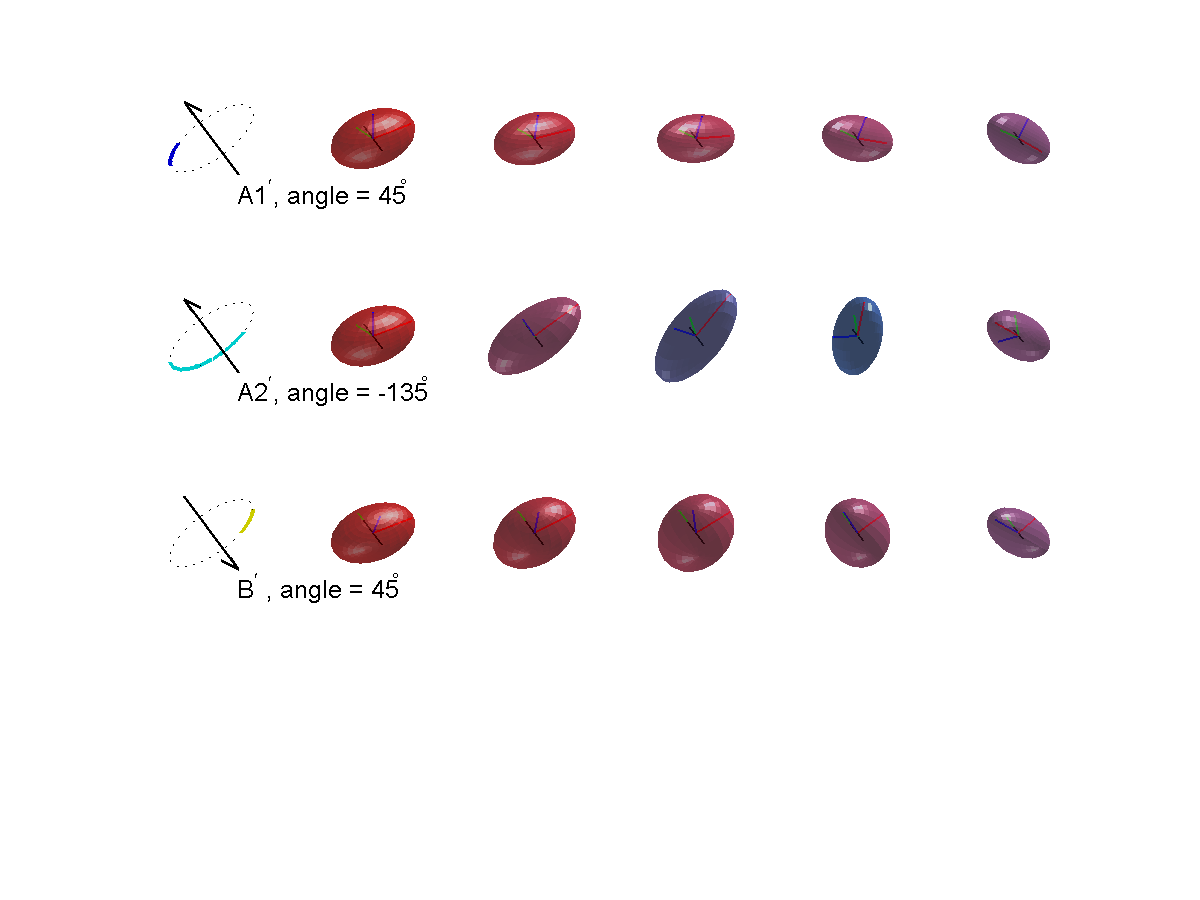}\vspace{-1.0in}
\includegraphics[width=1\textwidth]{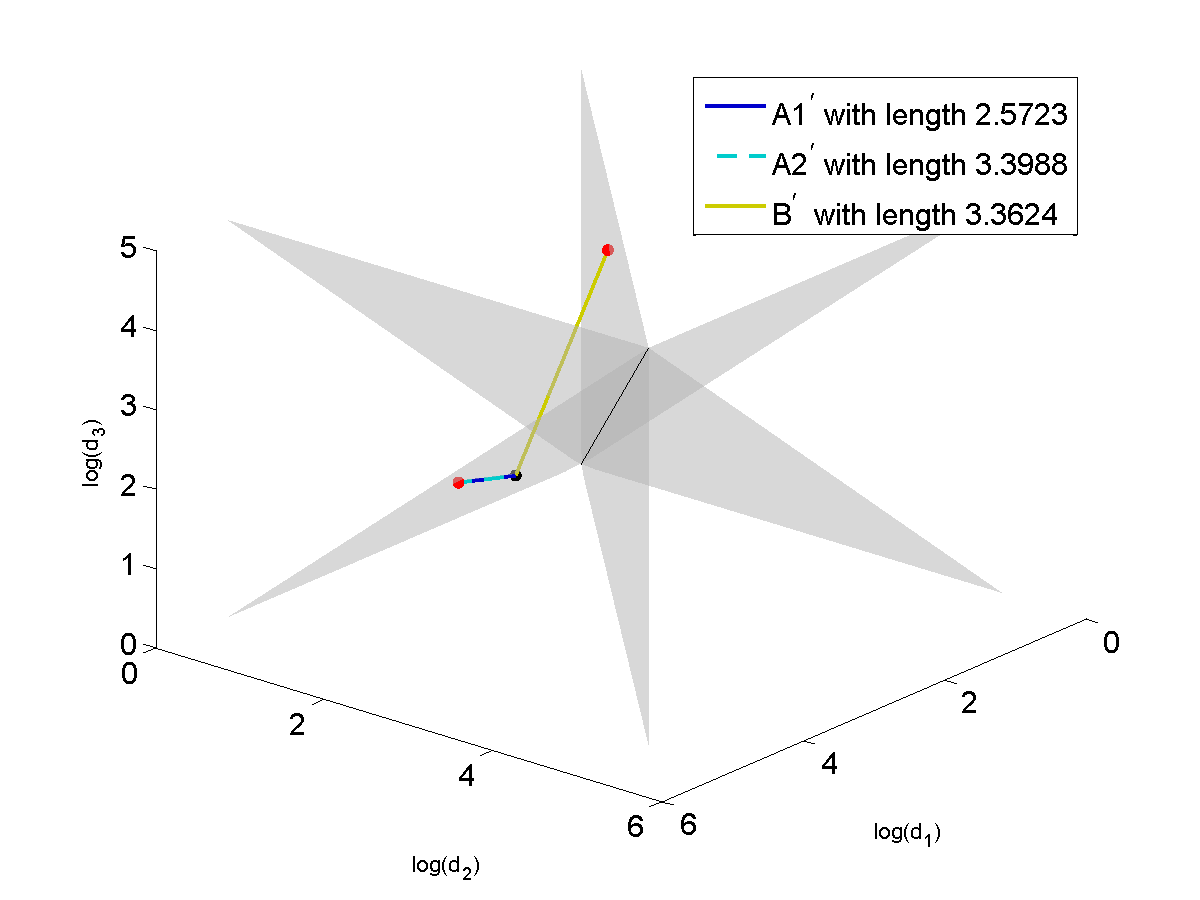}
  \caption{An example for the case $\{{\rm A}_1'\}$ (The scaling-rotation curve corresponding to ${\rm A}_1'$ is the unique MSSR curve.) \label{fig:A1pos} Each of the cases ${\rm A}_1'$, ${\rm A}_2'$, and ${\rm B}'$ represents the corresponding scaling-rotation curves defined in Table \ref{classes}, whose length can be found in the legend.  Since the rotational degrees of freedom have been projected out in the bottom panel, the relative lengths of the straight-line segments do not accurately reflect the relative lengths of the curves.
Note that $m > 0$ in this example. The eigenvalue paths corresponding to ${\rm A}_1'$ and ${\rm A}_2'$ stay in a single connected component of $\Dc_{\J_1}$ (represented by one of the shaded open half-planes) and their shapes are all prolates (if both $X$ and $Y$ are prolates) or all oblates (if both $X$ and $Y$ are oblates).  On the other hand, the eigenvalue path corresponding to ${\rm B}'$ travels through $\Dc_{\rm top}$ and another shaded plane corresponding to $\J_2$. That is, the  shape-classification changes are  (prolate $\to$ tri-axial $\to$ oblate $\to$ tri-axial $\to$ prolate).
}
\end{figure}

\begin{figure}[tbp!]
  \centering
\includegraphics[width=1\textwidth]{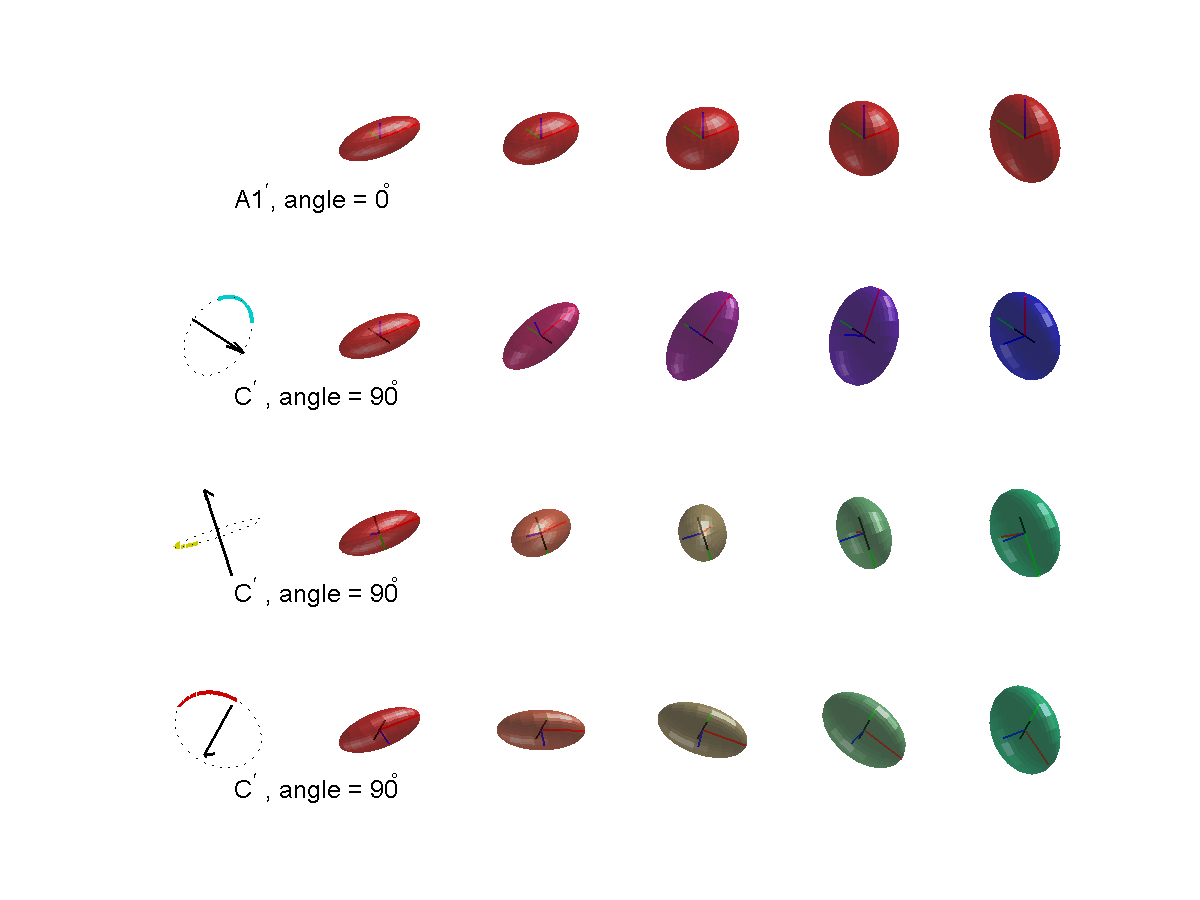}\vspace{-0.0in}
\includegraphics[width=1\textwidth]{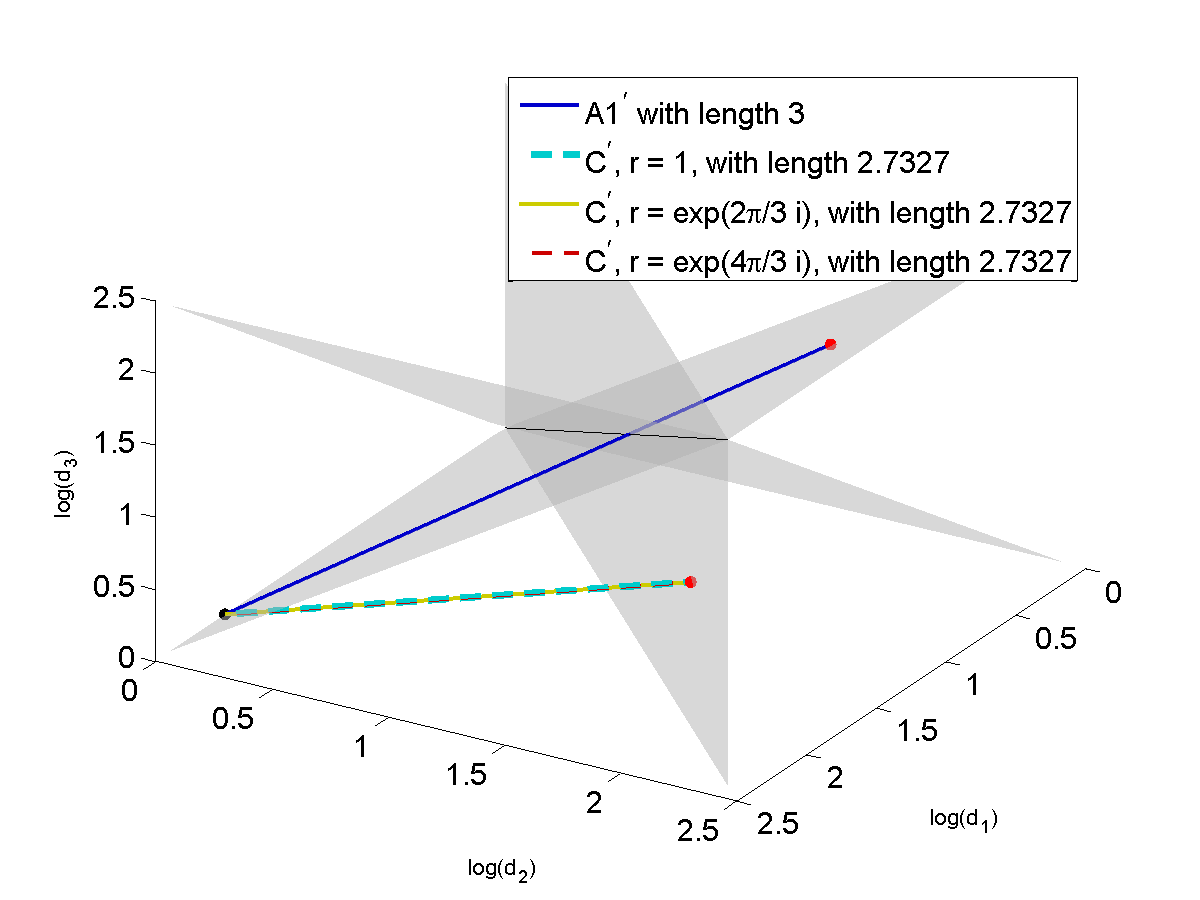}
  \caption{An example for the case $\{\rm{}C'\}$. For each of the case ${\rm C}'$, the rotation axis (depicted as the black line segment) is orthogonal to the (red) major semi-axis of $X$. It is shown in the proof of Theorem \ref{allMSR}(i) that there is one-to-one correspondence between the ``equator'' of $X$ and the family $\mathcal{M}_{{\rm C}'}$; see also Remark \ref{rem:proob}.  \label{fig:Cp}
}
\end{figure}

\begin{enumerate}
\item   $\{{\rm A}_1'\}$. Unique MSSR curve $\chi$. If $m > 0$, the shapes of $\chi(t)$ are always prolate (or oblate) for all $t \in [0,1]$; if $m<0$, the shape-classification changes of the MSSR curve $\chi(t)$ evaluated from $t= 0$ to $t= 1$ are either (oblate $\to$ sphere$\to$  prolate) or (prolate $\to$  sphere $\to$ oblate). See Fig.~\ref{fig:A1pos} for an example with $m>0$.

\item  $\{{\rm B}'\}$. Unique MSSR curve. 
        For prolate $X$, the shape-classification changes are  (prolate $\to$
tri-axial $\to$ oblate $\to$  tri-axial $\to$  prolate) if $m > 0$, (prolate $\to$
tri-axial $\to$  oblate) if $m < 0$. For oblate $X$, interchange ``oblate" and ``prolate" in these shape-classification changes.

\item   $\{{\rm A}_1',{\rm A}_2'\}$. Two MSSR curves with rotation angle $\pi/2$.  The shape-classification changes of ${\rm A}_2'$ are the same as those of ${\rm A}_1'$, and dependent on the sign of $m$; see item 1.

\item  $\{{\rm A}_1',{\rm A}_2',{\rm B}'\}$. Three MSSR curves (two with rotation angle $\pi/2$ and the other involving no rotation).

\item  $\{{\rm C}'\}$. Uncountably many MSSR curves. See Fig.~\ref{fig:Cp}. The set of MSSR curves is in natural one-to-one correspondence with $S_\mathbb{C}^1$. In this case $m < 0$, and the corresponding shape changes are the same as the case $\{B'\}$, with $m<0$: (prolate $\to$
tri-axial $\to$  oblate).
\end{enumerate}

\end{document}